\documentclass{cocv}
\usepackage{amsfonts,amssymb,amsmath}
\usepackage{graphicx,graphics,epsfig,color}
\usepackage{multicol}
\usepackage[colorinlistoftodos, textwidth=35mm, shadow]{todonotes}
\usepackage{rgsMacros}

\def\startmodifnew{\color{red}}
\def\stopmodifnew{\color{black}\normalcolor}

\let\labelindent\relax
\usepackage{enumitem}
\usepackage{balance}
\usepackage{enumitem}

\newcommand\red[1]{{\color{red}#1}}
    \definecolor{gray}{rgb}{0.33,0.4,0.47}\def\gray#1{{\color{gray}#1}}
    \definecolor{steelblue}{rgb}{0,.42,.7}\def\steelblue#1{{\color{steelblue}#1}}
    \definecolor{britishgreen}{rgb}{0,0.26,0.15}\def\britishgreen#1{{\color{britishgreen}#1}}
    \definecolor{navyblue}{rgb}{0,0,.8}\def\navyblue#1{{\color{navyblue}#1}}
    \definecolor{olivegreen}{rgb}{0.14,0.29,0}\def\olivegreen#1{{\color{olivegreen}#1}}
    \definecolor{myred}{rgb}{0.86,0.1,0.16}\def\myred#1{{\color{myred}#1}}

\newif\ifitsdraft
\def\itsdraft{\global\itsdrafttrue}
  \itsdraft  

\newtheorem{exe}{Example}
\newtheorem{corol}{Corollary}
\newtheorem{ass}{Assumption}
\newtheorem{defin}{Definition}
\newtheorem{cla}{Claim}
\newtheorem{rem}{Remark}
\newtheorem{lem}{Lemma}
\newtheorem{prop}{Proposition}
\newtheorem{thm}{Theorem}
\newtheorem{fct}{Fact}
\newtheorem{prob}{Problem}
\newenvironment{lemma}{\begin{lem}}{\hfill $\square$ \end{lem}}
\newenvironment{proposition}{\begin{prop}}{\hfill $\square$ \end{prop}}
\newenvironment{corollary}{\begin{corol}}{\hfill $\square$ \end{corol}}
\newenvironment{example}{\begin{exe}\rm }{\hfill $\square$ \end{exe}}
\newenvironment{remark}{\begin{rem}\rm }{\hfill $\bullet$ \end{rem}}
\newenvironment{assumption}{\begin{ass}}{\hfill $\bullet$ \end{ass}}
\newenvironment{theorem}{\begin{thm}}{\hfill $\square$ \end{thm}}
\newenvironment{definition}{\begin{defin}}{\hfill $\bullet$ \end{defin}}
\newenvironment{claim}{\begin{cla}}{\hfill $\bullet$ \end{cla}}

\newenvironment{problem}{\begin{prob}}{\hfill $\bullet$ \end{prob}}

\begin{document}

\title{\LARGE \bf A converse robust-safety theorem for differential inclusions}

\author{Mohamed Maghenem} \address{CNRS, Gipsa-lab,  Grenoble INP,  Universit{\'e} Grenoble Alpes,  Grenoble, France. Email: mohamed.maghenem@gipsa-lab.fr.}

\author{Masoumeh Ghanbarpour} \address{Department of Electrical and Computer Engineering, University of Colorado Bulder,  Colorado,  USA.  
\\
 Email:  masoumeh.ghanbarpour@gmail.com }


\begin{abstract}
This paper establishes the equivalence between robust safety and the existence of a barrier function certificate for differential inclusions.   More precisely,  for a robustly-safe differential inclusion, a barrier function is constructed as the 
time-to-impact function with respect to a specifically-constructed reachable set.  Using techniques from 
set-valued and nonsmooth analysis, we show that such a function,  although being possibly discontinuous,  certifies robust safety by verifying a condition involving the system's solutions.  Furthermore, we refine this construction,  using smoothing techniques from the literature of converse Lyapunov theory, to provide a smooth barrier certificate that certifies robust safety by verifying a condition involving only the barrier function and the system's dynamics. In comparison with existing converse robust-safety theorems, our results are more general as they allow the safety region to be unbounded, the dynamics to be a general continuous set-valued map,  and the solutions to be non-unique.  
\end{abstract}   
            

\maketitle

\section{Introduction}  

Safety for a dynamical system requires the solutions  starting from a given set of initial conditions to never reach a given unsafe set  \cite{prajna2007framework}.   Depending on the application,  reaching the unsafe set may correspond to 
non-applicability of a predefined feedback law,  due to saturation or a change in the dynamics or,  simply,  due to collisions with physical obstacles. Ensuring safety is in fact key in many engineering applications including traffic regulation \cite{ersal2020connected},  aerospace \cite{9656550}, and human-robot interactions \cite{9788028}. 

\subsection{Motivation}

This notion of safety is not robust in nature, as it is possible to construct safe differential equations that become unsafe when arbitrarily small perturbations are added to their right-hand side \cite[Example 1]{9683684}.   As a result, we say, roughly speaking, that a dynamical system is robustly safe if it remains safe in the presence of a perturbation term added to its dynamics.  This 
robust-safety notion was first introduced in \cite{wisniewski2016converse} for systems defined on compact manifolds.  A similar notion is studied in \cite{ratschan2018converse,  9444774} for continuous-time systems modeled by differential equations.  The same notion is considered in \cite{9683684} in the context of differential inclusions,  which generalize differential equations by allowing the right-hand side to be a general set-valued map \cite{aubin2012differential},  and thus the solutions to be non-unique.  

As the analytical expression of the solutions of a dynamical system are usually impossible to obtain, and since their precise approximation can be computationally expensive,    barrier functions are widely used to study safety and robust safety without
computing or approximating the  solutions. This is analogous to Lyapunov theory for stability.  We recall that a \textit{barrier function candidate} is a scalar function with opposite signs on the initial and the unsafe subsets. Furthermore,  it certifies safety, or robust safety, by satisfying an inequality constraint involving the barrier function candidate itself and the system's dynamics.  In which case,  the barrier function candidate becomes a \textit{barrier certificate} \cite{9705088, 9683684}.   
 Such conditions are well documented in the literature of safety under different smoothness properties of the barrier function candidate  \cite{ames2014control,
konda2019characterizing, 10.1007/978-3-642-39799-8_17} and the system's dynamics \cite{draftautomatica}. 
Furthermore, in the context of robust safety,  when an upper bound on the perturbation is known,   safety conditions involving  the (worst-case) perturbed dynamics are used in \cite{liu2020converse, JANKOVIC2018359}.  In \cite{seiler2021control},  specific classes of perturbations, solution to  some dynamical models and verifying a certain integral constraint are considered.  Perturbation-free conditions ensuring robust safety are proposed in  \cite{9683684, RubSafPI},  provided that mild regularity assumptions on the dynamics and the barrier function candidate hold.  Showing the necessity of such perturbation-free conditions is the main subject of the current paper. 

 \subsection{Background}

Converse safety and robust-safety problems pertain to show the existence of a barrier certificate for safety and robust safety, provided that the system is safe and robustly safe, respectively.

\begin{itemize}
\item In the context of safety, it is shown in \cite{9705088} that, in general, the existence of a continuous barrier certificate is not necessary for safety,   unless special cases are considered \cite{prajna2005necessity}.   Alternatively, time-varying barrier certificates are introduced and their existence is shown in \cite{9705088} to be necessary as well as sufficient,  under some assumptions on the system.    
\item In the context of robust safety,   \cite{wisniewski2016converse} solved the converse robust-safety problem by constructing a smooth barrier certificate for systems defined on smooth and compact manifolds,  provided that \blue{the system's dynamics are represented by a smooth single-valued map},  the initial and unsafe sets are compact and disjoint,  and a \textit{Meyer} function exists.     

\item In \cite{ratschan2018converse},  when the  system's dynamics are represented by a smooth 
single-valued map,  the complement of the unsafe set is bounded, and the closures of the initial and unsafe sets are disjoint,  robust safety is shown to be equivalent to the existence of a smooth barrier certificate. In particular,  to prove the converse robust-safety theorem  in \cite{ratschan2018converse}, the reachable set, denoted by $K_{\bar{\epsilon}}$, along the solutions to a perturbed version of the system starting from the initial set is introduced, where the subscript  $\bar{\epsilon}$ stands for the perturbation term added to the 
original-system's dynamics. 
After that, \blue{a barrier function candidate is defined, at any point in the state space, as the first to impact the boundary of the set $K_{\bar{\epsilon}}$ by a solution starting from that point}.  Such a function is shown to be a valid barrier function candidate. Although being only continuous, \blue{it is shown to be strictly decreasing along the solutions to the original system, when $\bar{\epsilon}$ is a robustness margin}.    
After that, a smooth barrier certificate is deduced using boundedness of the safe set \blue{and the density property of the class of smooth functions} in the space of continuous functions.   

\item A very similar converse robust-safety theorem is established in \cite{9444774} using converse Lyapunov theorems for asymptotic stability.   Indeed,  when either the reachable set $K_{\bar{\epsilon}}$ is bounded or the system's dynamics are represented by a globally Lipschitz map, and a uniform separation exists between the unsafe region and the set $K_{\bar{\epsilon}}$,  the latter set is shown to be uniformly asymptotically stable for the original system.  As a result,  existing converse Lyapunov theorems are used to show that a smooth barrier certificate exists.  
\end{itemize}
 
 \subsection{Contribution} 
 
In this paper,  we prove two  converse robust-safety theorems under mild regularity assumptions on the system's dynamics. \blue{Indeed, the latter is allowed to be represented by a set-valued map. Moreover,  we do not restrict the safety region to be bounded}. As in \cite{ratschan2018converse}, we show that blue{a specifically defined} 
time-to-impact function with respect to  the boundary of a reachable set $K_{\bar{\epsilon}}$,  when $\bar{\epsilon}$ is a robustness margin, is strictly decreasing when evaluated along the solutions to the original system lying on a neighborhood of $K_{\bar{\epsilon}}$.  However,  \blue{since the solutions are not necessarily unique},  this function is not necessarily continuous. Nonetheless, for an appropriate choice of the robustness margin $\bar{\epsilon}$,  the constructed function is shown to be a (non-smooth) barrier certificate. That is, it satisfies a sufficient condition for robust safety that is \textit{non-infinitesimal}; namely,  a condition involving the system's solutions.  
To construct a smooth barrier certificate, inspired by \cite{SUBBARAMAN201654}, we propose to smoothen the constructed nonsmooth one.  However, for the resulting smooth function to be a barrier function candidate; namely,  to have opposite signs on the initial and unsafe sets,  we need to carefully choose the set with respect to which the time-to-impact function is defined, as well as the different parameters involved in its construction.  As a consequence,  we show the existence of smooth barrier certificate provided that the system's dynamics are represented by a set-valued map that is continuous and the closures of the initial and unsafe sets are disjoint.   Finally,  we show the utility of our converse result in the context of safety for self-triggered control systems. 

 The rest of the paper is organized as follows.  Preliminaries on set-valued maps,  differential inclusions,  and  invariance and attractivity notions are in Section \ref{Sec.2}.  The problem formulation,  the main results,  and a motivational example  are in Section \ref{Sec.3}.   Preparatory materials towards the proofs of the main results are in Section \ref{Sec.4}.  
The proofs  of the main results are in Sections \ref{Sec.5} and \ref{Sec.6}. \blue{Finally, intermediate technical results are reported in the Appendix.} 

\blue{A preliminary version of this work is in \cite{9682926}, where the solutions are assumed to be forward complete. This  is not the case here. Furthermore,  proofs, detailed explanations, and the application example are not present in latter reference. }

\textbf{Notation.} 
For $x$,  $y \in \mathbb{R}^n$,  $x^{\top}$ denotes the transpose of $x$,  $|x|$ the Euclidean norm of $x$ and $\langle x, y \rangle:= x^\top y $ the inner product between $x$ and $y$.  For a set $K \subset \mathbb{R}^n$, 
we use $\mbox{int}(K)$ to denote its interior,  $\partial K$ to denote its boundary,  $U(K)$ to denote any open neighborhood of the set $K$,  and $|x|_K$ to denote the distance between $x$ and the set $K$.   Furthermore, we use $C_K(x)$ to denote the \textit{contingent} cone of $K$ at $x$, which is given by 
$$ C_K(x) := \left\{ v \in \mathbb{R}^n: \liminf_{h \rightarrow 0^+} |x + h v|_K/h = 0 \right\}.   $$ 
 For $O \subset \mathbb{R}^n$,  $K \backslash O$ denotes the subset of elements of $K$ that are not in $O$.  \blue{The sum of two sets is defined as $O + K := \{ x+y : x \in O, y \in K \}$}.    \textcolor{blue}{For a function $\phi : \dom \phi \rightarrow \mathbb{R}^m$,  $\dom \phi \subset \mathbb{R}^n$ denotes the domain of definition of $\phi$}.   By $F : \mathbb{R}^m \rightrightarrows \mathbb{R}^n $,  we denote a set-valued map associating each element $x \in \mathbb{R}^m$ \blue{with} a subset $F(x) \subset \mathbb{R}^n$.  In particular,  $\text{Proj}_{K} : \mathbb{R}^n \rightrightarrows K$ represents the projection set-valued map on $K$; namely, 
$ \text{Proj}_{K}(x) := \{y \in K : |x-y| = |x|_{K} \}$.  For a set $D \subset \mathbb{R}^m$,  $F(D) := \{ \eta \in F(x) : x \in D \}$.    For a differentiable map $B : \mathbb{R}^n \rightarrow \mathbb{R}$,  \textcolor{blue}{$\nabla_{x_i} B$ denotes  the derivative of $B$ with respect to $x_i$},  $i \in \{1,2,...,n\}$,  and $\nabla B$ denotes the gradient of $B$ with respect to $x$.   \ifitsdraft We say that $B \in \mathcal{L}^1$ if $|B|_1 := \int_{\mathbb{R}^n} |B(x)| dx$ is finite.  \fi
Finally,  we use $\mathbb{B}$ to denote the closed unit ball of appropriate dimension centered at the origin,  \blue{and $\mathcal{C}_+$ to denote the space of continuous   functions from  $\mathbb{R}^n$ to $\mathbb{R}_{>0}$. }

\section{Preliminaries} \label{Sec.2}

\subsection{Set-valued vs single-valued maps} 

Consider a set-valued map $F: K \rightrightarrows \mathbb{R}^n$, where $K \subset \mathbb{R}^m$. 

\begin{itemize}
\item $F$ is \textit{outer semicontinuous} at $x \in K$ if,  for every sequence $\left\{x_i\right\}^{\infty}_{i=0} \subset K$ and for every sequence  $\left\{ y_i \right\}^{\infty}_{i=0} \subset \mathbb{R}^n$ with $\lim_{i \rightarrow \infty} x_i = x$, $\lim_{i \rightarrow \infty} y_i = y \in \mathbb{R}^n$, and $y_i \in F(x_i)$ for all $i \in \mathbb{N}$, we have $y \in F(x)$;  see \cite{rockafellar2009variational}. 

\item  $F$ is \textit{upper semicontinuous} at $x \in K$ if,  for each $\varepsilon > 0$,  there exists a neighborhood of $x$, 
denoted by $U(x)$,  such that for each $y \in U(x) \cap K$, $F(y) \subset F(x) + \varepsilon \mathbb{B}$; see \cite[Definition 1.4.1]{aubin2009set}.

\item  $F$ is  \textit{continuous} at $x \in K$  if,  for each $\epsilon > 0$,  there exists $\delta > 0$ such that 
\begin{align} \label{eqContin}
|F(x_1) - F(x_2) |_H \leq \epsilon \qquad  \forall x_1, x_2 \in x + \delta \mathbb{B},  
\end{align}  
where $|F(x) - F(y) |_H$ stands for the Hausdorff distance between the sets $F(x)$ and $F(y)$.   
 
\item  $F$ is \textit{locally bounded} at $x \in K$ if there exists a neighborhood of $x$, denoted by $U(x)$,  and 
$\beta > 0$ such that  $|\zeta| \leq \beta$ for all $\zeta \in F(y)$ and for all $y \in U(x) \cap K$.  
\end{itemize}

Furthermore,  the map $F$ is upper,  outer semicontinuous,  continuous, or locally bounded if, respectively,  so it is for all $x \in K$.

Consider a single-valued map $B: K \rightarrow \mathbb{R}$,  where $K \subset \mathbb{R}^m$. 

\begin{itemize}
\item  $B$ is \textit{lower semicontinuous} at $x \in K$ if, for every sequence $\left\{ x_i \right\}_{i=0}^{\infty} \subset K$ such that $\lim_{i \rightarrow \infty} x_i = x$, we have $\liminf_{i \rightarrow \infty} B(x_i) \geq B(x)$. 
\item  $B$ is \textit{upper semicontinuous} at $x \in K$ if, for every sequence $\left\{ x_i \right\}_{i=0}^{\infty} \subset K$ such that $\lim_{i \rightarrow \infty} x_i = x$, we have $\limsup_{i \rightarrow \infty} B(x_i) \leq B(x)$.  
\item  $B$ is \textit{continuous} at $x \in K$ if it is both upper and lower semicontinuous at $x$. 
\end{itemize}

Furthermore,  $B$ is upper, lower semicontinuous, or continuous if, respectively,  so it is for all $x \in K$. 

\subsection{Differential inclusions} 
\textcolor{blue}{We recall the notion of a 
Carath{\'e}odory solution to a differential inclusion of the form}
\begin{align} \label{eq.1}
\Sigma : \quad  \dot x \in F(x) \qquad x \in \mathbb{R}^n.
\end{align}
\begin{definition}
A function $\phi : \dom \phi  \rightarrow \mathbb{R}^n$,  \blue{with $\dom \phi \subset \mathbb{R}$ an interval containing $\{0\}$},  is a solution to $\Sigma$ if it is locally absolutely continuous and  $\dot{\phi}(t) \in F(\phi(t))$ for almost all $t \in \dom \phi$.
\end{definition}

A solution $\phi$ to $\Sigma$ is said to start from $x$ if $\phi(0) = x$.   A solution $\phi$ to $\Sigma$ is maximal if there is no solution $\psi$ to  $\Sigma$  such that $\psi(t) = \phi(t)$ for all $t \in \dom \phi$ and $\dom \phi$ \blue{is} strictly included in $\dom \psi$.  Furthermore,  we use $\mathcal{S}_{\Sigma}(x)$ to denote the set of maximal solutions $\phi$ to $\Sigma$ starting from $x$. 

\textcolor{blue}{ We now propose to view the sets of points reached by the solutions to $\Sigma$, starting from a given initial condition and over 
a given window of time,  as set-valued maps.  Indeed,  as in \cite[Section 4.2.]{refId0} and \cite[Page 104]{aubin2012differential},  we,  respectively,    recall the set-valued maps $R_{\Sigma} : \mathbb{R} \times \mathbb{R}^n \rightrightarrows \mathbb{R}^n$ and $R^b_{\Sigma} : \mathbb{R} \times \mathbb{R}^n \rightrightarrows \mathbb{R}^n \cup \emptyset $ given by 
\begin{align*}
 R_{\Sigma}(t,x)  := \{ \phi(s): \phi \in \mathcal{S}_{\Sigma}(x), ~ 
s \in \dom \phi \cap I_t \},  \quad
 R^b_{\Sigma}(t, x) & := \left\{ \phi(t) : \phi \in \mathcal{S}_{\Sigma}(x),  ~t \in \dom \phi  \right\},
\end{align*}
where $I_t := [\min\{0,t\}, \max\{0,t\}]$.  }
 In simple words,  the set $R_{\Sigma}(t,x)$ includes all the elements reached by the solutions to $\Sigma$ starting from $x$ over the interval $I_t$.  Furthermore,  the set $R^b_{\Sigma}(t,x)$ includes the value of the solutions starting from $x$ at $t$,  when $t$ is part of their domain.  
 
Finally,  we introduce the following assumption on $F$. 
\begin{assumption} \label{ass1} 
The map $F$ is upper semicontinuous and $F(x)$ is nonempty, compact,  and convex for all $x \in \mathbb{R}^n$.
\end{assumption} 

\textcolor{blue}{Assumption \ref{ass1} guarantees  the existence of a non-trivial solution from any $x \in \mathbb{R}^n$ as well as useful structural properties for the set of solutions to $\Sigma$;  see \cite{aubin2012differential,  refId0,  filippov2013differential}.  }
 
\begin{remark} \label{remplus}
\blue{Given a set-valued map $F : K \rightrightarrows \mathbb{R}^n$, where $K \subset \mathbb{R}^m$,  
we recall,  based on \cite[Theorem 5.19]{rockafellar2009variational} and \cite[Lemma 5.15]{goebel2012hybrid},  that  the following two properties are equivalent.
\begin{itemize}
\item  $F$ is upper semicontinuous and  $F(x)$ is compact for all $x \in K$. 
\item $F$ is outer semicontinuous and locally bounded.
\end{itemize}}
\end{remark}

\section{Problem formulation and results}  \label{Sec.3}

Given a set of initial conditions $X_o \subset \mathbb{R}^n$ and an unsafe set $X_u \subset \mathbb{R}^n$  such that $X_o \cap X_u = \emptyset$, 
we recall that $\Sigma$ is safe with respect to $(X_o,X_u)$ if,  for each solution $\phi$ with $\phi(0) \in X_o$,  we have 
$\phi(\dom \phi  \cap \mathbb{R}_{\geq 0} ) \subset \mathbb{R}^n \backslash X_u$.   \blue{In words,  the solutions starting from $X_o$ never reach the set $X_u$ at any positive time}.  Note that safety with respect to  $(X_o, X_u)$ is verified if and only if there exists a set $K \subset \mathbb{R}^n$,  with $X_o \subset K$ and $K \cap X_u = \emptyset$,  that is forward invariant.   
In turn,  a set $K \subset \mathbb{R}^n$ is forward invariant if,  for each solution $\phi$ to $\Sigma$ with $ \phi(0) \in K$, $\phi(\dom \phi \cap \mathbb{R}_{\geq 0}) \subset K$.  \blue{In words,  the solutions starting from $K$ never leave the set $K$ at any positive time}.

Next, we consider the perturbed version of $\Sigma$, denoted by $\Sigma_\epsilon$, and given by 
\begin{align} \label{eq.2}
\Sigma_\epsilon : \quad \dot{x} \in F(x) + \epsilon(x) \mathbb{B} \qquad  x \in \mathbb{R}^n.
\end{align} 
\blue{Following \cite{9683684}, we recall a definition of  robust safety,  which we study in this paper. }
\begin{definition}[Robust safety]
System $\Sigma$ is robustly safe with respect to $(X_o, X_u)$ if there exists $\epsilon \in \mathcal{C}_+$ such that $\Sigma_{\epsilon}$ in \eqref{eq.2} is safe with respect to $(X_o,X_u)$. 
The function $\epsilon$ is,  in this case,  named
 \textit{robust-safety margin}.
\end{definition}
We next recall the definition of a barrier function candidate.  
\begin{definition}
A scalar function $B : \mathbb{R}^n \rightarrow \mathbb{R}$ is a barrier function candidate with respect to $(X_o,X_u)$ if
\begin{align*} 
\begin{matrix} 
B(x) > 0 & \forall x \in X_u \quad \text{and} \quad B(x) \leq 0 & \forall x \in X_o. 
\end{matrix} 
\end{align*} 
\end{definition}
Note that a barrier function candidate 
$B$ defines the zero sub-level set 
\begin{align} \label{eq.4} 
K := \left\{ x \in \mathbb{R}^n : B(x) \leq 0 \right\}, 
\end{align} 
which necessarily verifies  
$$X_o \subset K \quad \text{and}  \quad K \cap X_u = \emptyset.  $$

\subsection{Sufficient conditions for robust safety} 

Since  we can guarantee  robust safety for $\Sigma$ by guaranteeing  safety for  
$\Sigma_\epsilon$,  for some  $\epsilon \in \mathcal{C}_+$,   then robust safety is verified if 
\begin{enumerate} [label={C0)},leftmargin=*]
\item  There exists a barrier function candidate $B : \mathbb{R}^n \rightarrow \mathbb{R}$ such that the set $K$ in \eqref{eq.4} is forward invariant for  $\Sigma_\epsilon$,  for some   $\epsilon \in \mathcal{C}_+$.  
\end{enumerate}

\blue{ The latter allows us to recall the following solution-dependent sufficient condition for robust safety \cite{9705088}. }

\begin{proposition}
$\Sigma$ is robustly safe with respect to $(X_o,X_u) \in \mathbb{R}^n \times \mathbb{R}^n$ provided that   
\begin{enumerate} [label={C1)},leftmargin=*]
\item \label{item:C2bis}  
\blue{There exists a barrier function candidate $B$ such that the set $K$ in \eqref{eq.4} is closed,  there exists $U(\partial K)$ an open neighborhood of $\partial K$,  and  there exists 
$\epsilon \in \mathcal{C}_+$ such that,  along every solution (not necessarily maximal)  
$\phi$  to $\Sigma_{\epsilon}$ satisfying $\phi(\dom \phi) \subset U(\partial K)$,    the map 
$t \mapsto B(\phi(t))$ is non increasing. }
\end{enumerate} 
\end{proposition}

Note that \ref{item:C2bis} does not require from $B$ to be smooth neither $F$ to satisfy Assumption \ref{ass1}. However, it involves the solutions to $\Sigma_\epsilon$ (hence, it also involves the perturbation term $\epsilon$).   

\blue{We now recall from \cite{9683684,RubSafPI} a sufficient condition for robust safety that uses only the barrier function candidate $B$ and the nominal dynamics $F$.}

\begin{proposition}
Consider system $\Sigma$ such that  Assumption \ref{ass1} holds. $\Sigma$ is robustly safe with respect to $(X_o,X_u) \in \mathbb{R}^n \times \mathbb{R}^n$ provided that 
\begin{enumerate} [label={C2)},leftmargin=*]
\item \label{item:C3} There exists 
a continuously differentiable barrier function candidate $B$ such that
\begin{align} 
\langle \nabla B (x), \eta \rangle < 0 \qquad \forall \eta \in F(x), \quad \forall x \in \partial K.   \label{eq.2cbis}
\end{align}
\end{enumerate} 
\end{proposition}
\blue{Note that \eqref{eq.2cbis} involves only the barrier function candidate $B$,  which is now required to be continuously differentiable,  and the nominal dynamics $F$, which is required to verify Assumption \ref{ass1}.} 

\subsection{Converse robust-safety theorems}

We start introducing the following two assumptions. 

\begin{assumption} \label{ass3} 
The map $F$ is continuous and $F(x)$ is nonempty, compact, and convex for all $x \in \mathbb{R}^n$.
\end{assumption}

\begin{assumption} \label{ass4-} 
$\cl(X_o) \cap X_u = \emptyset$. 
\end{assumption}

\blue{Using the latter two assumptions, we can formulate our first converse robust-safety theorem.}

\begin{theorem} \label{thm3}
Consider system $\Sigma$ that is robustly safe with respect to $(X_o,X_u)$ and such that Assumptions  \ref{ass3} and \ref{ass4-} hold.  Then,   \ref{item:C2bis} holds. 
\end{theorem}

In the next converse robust-safety theorem,   instead of Assumption \ref{ass4-},  we use the following relatively stronger assumption.  

\begin{assumption} \label{ass4} 
$\cl(X_o) \cap \cl (X_u) = \emptyset$. 
\end{assumption}   

\begin{theorem} \label{thm4}
Consider system $\Sigma$ that is robustly safe with respect to $(X_o,X_u)$ and such that Assumptions \ref{ass3}  and \ref{ass4} hold.  Then,   \ref{item:C3} holds.
\end{theorem}

\begin{remark} \label{remthm2}
Our proof of Theorem \ref{thm4} actually allows us to conclude that, for any $k \in \{1,2,...\}$, their exists a continuously-differentiable barrier certificate $B$ such that
$$ \langle \nabla B (x), \eta \rangle < -1 \qquad \forall \eta \in F(x),  \quad \forall x \in \partial K.  $$   
\end{remark}

\begin{remark}
\blue{ We believe that relaxing Assumption 2, and using Assumption 1 instead, would require a totally different approach to prove Theorems \ref{thm3} and \ref{thm4}.  Indeed,  a key property allowing our barrier-function construction is established in Lemma \ref{lem1} below.   This property does not hold when only Assumption \ref{ass1} is verified; 
see Example \ref{exp1}.   Similarly,  relaxing Assumptions \ref{ass4-} and \ref{ass4} would require using a completely different approach to prove Theorems \ref{thm3} and \ref{thm4},  respectively.  Indeed,  such separations between the two sets $X_o$ and $X_u$ are necessary to squeeze the zero-level set of the constructed barrier function between the two sets; see the forthcoming Section \ref{Sec.Separ} for more details. }
\end{remark}

\subsection{Application: Safety for self-triggered control systems}

Consider the control system $\Sigma_u$ given by
$$ \Sigma_u : \dot{x} = f(x,u)  \qquad  (x,u) \in  \mathbb{R}^n \times  \mathbb{R}^{m},  $$  
 where $f : \mathbb{R}^n \times  \mathbb{R}^{m} \rightarrow \mathbb{R}^n$ is a continuous function. Furthermore, we let   $\kappa: \mathbb{R}^n \rightarrow \mathbb{R}^{m}$ be a continuous feedback law such that the resulting closed-loop system 
 \begin{align} \label{eqSigmaST} 
 \Sigma : \dot{x} =  F(x) := f(x,\kappa(x))  \qquad x \in \mathbb{R}^n 
 \end{align}
 is safe with respect to $(X_o,X_u) \subset \mathbb{R}^n \times \mathbb{R}^n$.  
 
  In a self-triggered (ST) control framework \cite{di_benedetto_digital_2013},  we construct a monotonically increasing sequence  
  $\{t_i\}^{\infty}_{i=0} \subset \mathbb{R}_{\geq 0}$ and we ask the controller to   remain constant between each two time samples $t_i$ and $t_{i+1}$,  i.e.,  
  $$ u(t) = \kappa(x(t_i)) \qquad  \forall t \in [t_i, t_{i+1}), \qquad i \in \mathbb{N}. $$  
  Hence,  a solution 
$\phi$ to the ST closed-loop system must satisfy 
\begin{align*} 
\dot{\phi}(t) =  F(\phi(t))  + \Gamma(\phi(t),\phi(t_i))   \quad \forall t \in [t_i, t_{i+1}), \quad \Gamma(\phi(t),\phi(t_i)) :=  f(\phi(t),\kappa(\phi(t_i))) - f(\phi(t),\kappa(\phi(t))).
\end{align*}

Our goal here,  knowing that system $\Sigma$ in \eqref{eqSigmaST} is robustly safe,   is to address the following problem.  

\begin{problem} \label{probap}
\blue{Prove the existence of a sequence $\{t_i\}^{\infty}_{i=0} \subset \mathbb{R}_{\geq 0}$ and $T > 0$ such that $t_{i+1} - t_i  > T$ for all $i \in \mathbb{N}$ and the resulting ST closed-loop system is safe with respect to $(X_o,X_u)$.}  
\end{problem}

We will show that Theorem \ref{thm4} is key to solve Problem \ref{probap}.  Indeed,  Theorem \ref{thm4} allows us to formulate the following corollary. 

\begin{corollary} \label{coroo}
Consider the control system $\Sigma_u$ and a feedback law $\kappa$ such that the resulting closed-loop system $\Sigma$ in \eqref{eqSigmaST} is robustly safe with respect to $(X_o,X_u)$. 

Then, there exist a continuously-differentiable barrier certificate $B$ and
continuously-differentiable functions $\alpha : \mathbb{R}^n \rightarrow \mathbb{R}$ and $\gamma : \mathbb{R}^n \times \mathbb{R}^n \rightarrow \mathbb{R}$ such that
\begin{align} \label{eqConC-} 
\alpha(x) \geq 3/4 \quad \forall x \in \partial K,  \qquad 
\gamma (x,x) \leq 1/8 \quad \forall x \in \mathbb{R}^n, 
\end{align}
and
\begin{align} \label{eqConC}
\hspace{-0.6cm}
\langle \nabla B(x),  f(x,\kappa(y)) \rangle   \leq - \alpha(x) + \gamma(x,y) \qquad  \forall (x,y) \in \mathbb{R}^n \times \mathbb{R}^n.  
\end{align}
\end{corollary}

\begin{proof}
\blue{The existence of the 
continuously-differentiable barrier certificate $B$ is guaranteed by Theorem \ref{thm4}. Furthermore,  according to Remark \ref{remthm2}, we can choose the barrier certificate $B$ to satisfy
$$ \hat{\alpha}(x) := - \langle \nabla B(x) ,  F(x)  \rangle > 1 \quad \forall x \in \partial K \quad \text{and} \quad   \hat{\gamma}(x,x) := \langle \nabla B(x) ,   \Gamma(x,x)  \rangle  = 0 \quad \forall x \in \mathbb{R}^n. $$
As a result, using Whitney approximation theorem for continuous functions, we conclude that we can always find  continuously-differentiable functions $\gamma$ and $\alpha$ such that \eqref{eqConC-} holds and at the same time
 \begin{align} \label{eqprop31}
 \alpha(x) \leq  \hat{\alpha}(x) \quad  \forall x \in \mathbb{R}^n
 \quad \text{and}  \quad 
 \gamma(x,y)  \geq \hat{\gamma}(x,y) \quad \forall (x,y) \in \mathbb{R}^n \times \mathbb{R}^n.  
\end{align}
The choice in \eqref{eqprop31} allows us to verify 
\eqref{eqConC}. }
\end{proof}

The proposed solution to  Problem \ref{probap}
requires the following assumption.

\begin{assumption} \label{assexp}
\blue{The set $\mathbb{R}^n \backslash X_u$ is bounded, and there exists $\tau>0$ such that the solutions to $\Sigma_y : \dot{x} = f(x,\kappa(y))$ starting from $y \in \mathbb{R}^n \backslash X_u$ cannot blow up on the time  interval $[0,\tau]$.}  
\end{assumption}

Under the first part of Assumption \ref{assexp}, we conclude that the zero sub-level set $K$ of the barrier certificate  
 is compact. Hence, there exist $\beta>0$ and $T_1>0$ such that the following properties hold. 
 
\blue{\begin{itemize}
\item[Pr1)] For each $y \in G:=  \{x \in K: |x|_{\partial K} \leq \beta\}$, we have $\alpha(y) - \gamma(y,y)  \geq \alpha(y) - 1/8  \geq 1/4$. 
This is true because $\alpha$ is continuously differentiable. 
\item[Pr2)] The solutions  to $\Sigma_y : \dot{x} = f(x,\kappa(y))$ starting from  $y \in L :=   \{x \in K: |x|_{\partial K} \geq \beta\}$
remain in $K$ on the interval $[0,T_1]$.
This is because $f$ and $\kappa$ are continuous and $K$ is compact.
\end{itemize}}

The following result  addresses Problem \ref{probap}. 

\begin{proposition}
\blue{Consider the control system $\Sigma_u$ whose dynamics $f$ is continuous. Consider $(X_o,X_u) \subset \mathbb{R}^n \times \mathbb{R}^n$  and a continuous feedback law $\kappa$ such that 
the resulting closed-loop system $\Sigma$ in \eqref{eqSigmaST} is robustly safe with respect to $(X_o,X_u)$ and Assumption \ref{assexp} holds.

 Then,  a solution to Problem \ref{probap} is given by
$$ t_{i+1} := 
\left\{ 
\begin{matrix} 
t_i +  \max \{ T_1,  T_r(\phi(t_i)) \}  & \text{if} ~ \phi(t_i) \in L
\\
t_i +  T_r(\phi(t_i))    &  \text{if} ~ \phi(t_i) \in G,
\end{matrix} 
\right.
$$ 
where $G$ and $(T_1,L)$ are introduced in Pr1) and Pr2), respectively. Furthermore, the map $y \mapsto T_r(y)$ is defined, for all $y \in K$, as
\begin{align*}
T_r(y) & := 0 &  \text{if} ~ \alpha(y) - \gamma(y,y) \leq 0, 
\\
T_r(y) & := \left\{ 
\begin{matrix}
\tau &  \text{if} ~ M_r(\tau,y) \leq 0
\\ 
\min \left\{ \tau, \frac{2(\alpha(y) - \gamma(y,y))}{M_r(\tau,y)} \right\} & \mbox{otherwise}
\end{matrix}
\right.   &  \text{otherwise}, 
\end{align*}
where $\tau$ is introduced in Assumption \ref{assexp} and 
\begin{equation}
\label{eqMreq}
\begin{aligned}
M_r(\tau, y) & :=   \sup \{\langle \nabla_x \gamma(x, y), f(x, \kappa(y)) \rangle:  x \in R_{\Sigma_y}(\tau, y) \} 
+ \sup \{\langle - \nabla \alpha(x) , f(x, \kappa(y)) \rangle: x \in R_{\Sigma_y}(\tau, y)\},
\end{aligned}
\end{equation}
$R_{\Sigma_y}(\tau, y)$ is the set of points reached by the solutions to $\Sigma_y$ starting from $y$ over the window of time $[0,\tau]$, 
and the functions $\alpha$ and $\gamma$ are introduced in Corollary \ref{coroo}. }
\end{proposition}

\begin{proof}
Using Pr1), we conclude that $T_r(y)>0$ for all $y \in G$. Furthermore, using Pr2), we conclude that $t_{i+1} - t_i > 0$ for all $i \in \mathbb{N}$. 

Next, we show that any solution $\phi$ to $\dot{x} = f(x,\kappa(y))$,  starting from $y \in K$,  satisfies 
$ \phi([0,T_r(y)]) \subset K$. This is enough to conclude that the self-triggered closed-loop system is safe with respect $(X_o,X_u)$.
To this end, we note that $\phi$ is locally absolutely continuous,  and  since $\alpha$ and $\gamma$ are continuously differentiable, it follows that  $t \mapsto \alpha(\phi(t))$ and $t \mapsto \gamma(\phi(t),y)$ are also locally absolutely continuous.   Hence,  for each $t \in \dom \phi$,  there exists a sequence 
 $\{\tau_n\}_{n=0}^N\subset [0,t]$,  with  $N \in \mathbb{N}^* \cup \{\infty\}$, such that 
 $$ \lim_{n\rightarrow N} \tau_n = t, \qquad   \tau_n - \tau_{n-1} > 0, $$  
 and the maps $t \mapsto \alpha(\phi(t))$ and $t\mapsto \gamma(\phi(t), y)$ are  differentiable on each $(\tau_{n-1}, \tau_{n})$. 

 \textcolor{blue}{Consider the map} $\bar{\gamma}(\cdot) := \gamma(\cdot,y)$ 
 and note that
$$ \bar{\gamma}(\phi(t)) - \bar{\gamma}(y) = \sum_{n=1}^N  \left( \bar{\gamma}(\phi(\tau_n)) - \bar{\gamma}(\phi(\tau_{n-1})) \right), \quad \alpha(\phi(t)) - \alpha(y)  = \sum_{n=1}^N \left[\alpha(\phi(\tau_n)) - \alpha(\phi(\tau_{n-1}))\right].  $$
As a result, using the classical mean-value theorem,  we conclude that, for each $n\in\{1,2,...,N\}$, there exist $c_n, d_n \in (\tau_{n-1}, \tau_{n})$ such that
\begin{align*}
\bar{\gamma}(\phi(t)) - \bar{\gamma}(y)  = \sum_{n=1}^{N} \left( \frac{d}{dt} \bar{\gamma}(\phi(t)) \Big\vert_{t=c_n} (\tau_{n} - \tau_{n-1}) \right),
\quad 
\alpha(\phi(t)) - \alpha(y)  = \sum_{n=0}^{N}
\left(\frac{d}{dt} \alpha(\phi(t)) \Big\vert_{t=d_n}  (\tau_{n+1} - \tau_n)\right). 
\end{align*}
As a result,  when $t \in [0,\tau]$,  we conclude that
\begin{equation}
    \label{eqciteit}
    \begin{aligned} 
    \bar{\gamma}(\phi(t))  - \bar{\gamma}(y) & \leq  t \sup\{\langle \nabla\bar{\gamma}(x),f(x,\kappa(y)) \rangle:  x \in R_{\Sigma_y}(\tau,y) \},
\\
-\alpha(\phi(t)) + \alpha(y)    
   & \leq t \sup \{\langle - \nabla \alpha(x), f(x,\kappa(y)) \rangle :   x \in R_{\Sigma_y}(\tau,y) \}.
\end{aligned}
\end{equation}
Next, we note that 
 $$ \frac{d}{dt} B(\phi(t)) = \langle \nabla B(\phi(t)), \dot{\phi}(t) \rangle \qquad   \text{for almost all} ~ t \in \dom \phi.  $$  
Integrating the previous equality 
from $0$ to $t \leq \tau$, we obtain
\begin{align*} 
 B(\phi(t)) -B(y) & \leq \int_0^t [-\alpha(\phi(s)) + \bar{\gamma} (\phi(s))]  ds  \leq \int_0^t [-\alpha(y) + \bar{\gamma}(y) + s M_r(\tau, y)] ds    
\\ & 
\leq  - t (\alpha(y) - \bar{\gamma}(y))  + \frac{t^2}{2} M_r(\tau, y) \qquad \forall t \in [0,\tau].
\end{align*} 
To obtain the latter inequalities, we used \eqref{eqConC}, \eqref{eqciteit}, and \eqref{eqMreq}. 

Now, we note that, when $\alpha(y) - \bar{\gamma}(y) > 0$ and $M_r(\tau,y) > 0$, it follows that
$$ B(\phi(t)) - B(y) \leq 0 \qquad  \forall t \in \left[ 0, 2 \frac{ \left(\alpha(y) - \bar{\gamma}(y) \right)}{M_r(\tau,y)} \right] \cap [0,\tau]. $$ 
 Otherwise,  when $\alpha(y) - \bar{\gamma}(y) > 0$ and $M_r(\tau,y) \leq 0$, we conclude that
$$ B(\phi(t)) - B(y) \leq 0 \qquad  \forall t \in [0,\tau]. $$ 
The latter is enough to conclude that the proposed triggering sequence guarantees safety for the resulting self-triggered 
closed-loop system.

In the rest of the proof, 
we show the existence of $T > 0$ such that $ t_{i+1} - t_i  > T$ for all  $i \in \mathbb{N}$. To do so,  it is enough to show that the maps $y \mapsto T_r(y)$ is lower semicontinuous on $G$. Indeed, since it is already positive on $G$ and $G$ is compact,  it would follow using \cite[Theorem B.2]{puterman2014markov} that $T_r$ reaches its minimum on $G$. As a result, we can take $T := \min \{ \min \{T_r(y): y \in G \}, T_1 \}$. 

To show that $T_r$ is lower semicontinuous on $G$, we start noting that the set-valued map $y \mapsto R_{\Sigma_y}(\tau,y)$ is upper semicontinuous with compact images on $G$; see Lemma \ref{lem1-} in the Appendix. Next,  we use \cite[Theorem 1.4.16]{aubin2009set},  under smoothness properties of $\gamma$ and $\alpha$,  to conclude that the single-valued map $y \mapsto M_r(\tau,y)$ is upper semicontinuous on $G$. \blue{ It is also locally bounded on $G$ since so is the set-valued map $y \mapsto R_{\Sigma_y}(\tau,y)$.} 

Now,  since  $\alpha$ is positive and continuous on $G$,  we conclude that $y \mapsto \frac{\alpha(y) - \bar{\gamma}(y)}{M_r(\tau,y)}$ is lower semicontinuous and positive on $G$. To complete the proof, we consider a sequence $ \{x_i\}^{\infty}_{i=0} \subset G$ that converges to $x_o \in G$. Since $y \mapsto M_r(\tau,y)$ is upper semicontinuous, we conclude that $\limsup_{i \rightarrow \infty} M_r(\tau,x_i) \leq M_r(\tau,x_o)$.   \textcolor{blue}{Moreover, by selecting an appropriate subsequence, we can assume, without loss of generality, that} $ \liminf_{i \rightarrow \infty} T_r(x_i)  =    \lim_{i \rightarrow \infty} T_r(x_i)  = a > 0$. Also, since $M_r$ is locally bounded,  \textcolor{blue}{one can select another subsequence to conclude the existence of} $\beta \in \mathbb{R}$ such that 
$  \lim_{i \rightarrow \infty} M_r(\tau,x_i) = \beta \leq M_r(\tau,x_o)$. 
 \textcolor{blue}{To finalize the proof, we} distinguish the following three scenarios:  
\begin{enumerate}
\item  When $\beta < 0$, we conclude that 
$M_r(\tau,x_i) < 0$ for all $i \in \mathbb{N}$ sufficiently large. In this case,  $T_r(x_i) = \tau$ for all $i \in \mathbb{N}$ sufficiently large, thus, $\lim_{i \rightarrow \infty} T_r(x_i) =  
\tau \geq T_r(x_o)$.   
\item  When $\beta > 0$, we conclude that $M_r(\tau,x_i) > 0$ for all $i \in \mathbb{N}$ sufficiently large.   In this case, we conclude that 
$T_r(x_i) = \min \left\{ \tau,  \frac{\alpha(x_i) - \bar{\gamma}(x_i)}{M_r(\tau,x_i)}  \right\}$  for all $i \in \mathbb{N}$ large. Hence,   
$
\lim_{i \rightarrow \infty} T_r(x_i)  =  
\min \left\{\tau, \lim_{i \rightarrow \infty} \frac{ \alpha(x_i) - \bar{\gamma}(x_i)}{ M_r(\tau,x_i)} \right\} 
\geq \min \left\{ \tau,   \frac{\alpha(x_o) - \bar{\gamma}(x_o)}{M_r(\tau,x_o)} \right\} \geq T_r(x_o)$.
\item   When $\beta = 0$, we conclude that 
$| \frac{\alpha(x_i) - \bar{\gamma}(x_i)}{ M_r(\tau,x_i)} | \geq \tau$ for all $i \in \mathbb{N}$ sufficiently large. Hence, 
$T_r(x_i) = \tau$ for all $i \in \mathbb{N}$ sufficiently large.
\end{enumerate}
\end{proof} 

\section{Preparatory material} \label{Sec.4}

In this section, we present key intermediate results that allow us to construct a smooth barrier certificate for robustly-safe systems. 

\subsection{Perturbed differential inclusions} 

\blue{We will show that when  Assumption \ref{ass3} is verified,  given $x \in \mathbb{R}^n$ and $\epsilon \in \mathcal{C}_+$,  we can find $T > 0$ such that  $R^b_{\Sigma}(T,x) \subset  \text{int}\left( R^b_{\Sigma_\epsilon}(T,x) \right)$}.   A similar result requiring $F$ to be locally Lipschitz can be found in \cite{puri1995varepsilon} and \cite{ratschan2018converse}.

\begin{lemma} \label{lem1}
Consider system $\Sigma$  such that Assumption  \ref{ass3} holds.  Then,  for each $x \in \mathbb{R}^n$ and for each $\epsilon \in \mathcal{C}_+$,  
there exists $T>0$ such that,  for each 
$t \in (0,T]$,  there exists $\delta > 0$ such that 
\begin{align}
y + \delta \mathbb{B} & \subset R^b_{\Sigma_\epsilon}(t,x) \qquad \forall y \in R^b_{\Sigma}(t,x) \backslash \{x\}, 
\label{eqreachcover} \\
x + \delta \mathbb{B} & \subset 
 R^b_{\Sigma_\epsilon}(-t,y) = R^b_{\Sigma^-_\epsilon}(t,y) \qquad \forall y \in R^b_{\Sigma}(t,x) \backslash \{x\}, \label{eqreachcover1} 
\end{align}
where $\Sigma^- :   \dot{x} \in - F(x)$ with $x \in \mathbb{R}^n$.  
\end{lemma}

\begin{proof}
Given $x \in \mathbb{R}^n$ and $\epsilon \in \mathcal{C}_+$,  we pick $\Delta \in (0,1)$ such that,  for each $x_1,x_2 \in x + \Delta \mathbb{B}$ and for each $f_1 \in F(x_1)$, there exists $f_2 \in F(x_2)$ such that
\begin{align} \label{eqcontin}
|f_1 - f_2| \leq \underline{\epsilon}/2, \quad \underline{\epsilon} := \min \{\epsilon(y) : y \in x + \mathbb{B} \}. 
\end{align}  
The latter is possible since the set-valued map $F$ is assumed to be continuous. 

Now,  using Lemma \ref{lemprepre} in the Appendix,  we conclude the existence of $b$ and $\bar{T}>0$ such that $R_{\Sigma} (\bar{T},  (x + b \mathbb{B}))$ and 
$R_{\Sigma^-} (\bar{T},  (x + b \mathbb{B}))$
are bounded. As a result,  we invoke Lemma \ref{lem1-} in the appendix to conclude that  $R_\Sigma$ and $R_{\Sigma^-}$ are outer semicontinuous and locally bounded  on 
$ [0,\bar{T}] \times (x + b \mathbb{B})$,  which implies,  according to Remark \ref{remplus},  that $R_\Sigma$ and $R_{\Sigma^-}$ are upper semicontinuous and have  compact images on 
$[0,\bar{T}] \times (x + b \mathbb{B})$. 
Hence,   we can find 
$T \in (0,\bar{T}]$ such that 
\begin{align} 
R_{\Sigma}(T,x) & \subset 
\left(x + \frac{\Delta}{2} \mathbb{B} \right), 
\label{eqcontin1a}
\\
R_{\Sigma^{-}}(T,R_{\Sigma}(T,x)) & \subset 
\left(x + \frac{\Delta}{2} \mathbb{B} \right).
\label{eqcontin1b}
\end{align}

To prove \eqref{eqreachcover}, we let $\delta \in (0, \Delta/2]$ and we let $y \in R^b_{\Sigma}(t,x) \backslash \{ x \}$ for some $t \in (0,T]$.  Hence,  there exists a solution $\phi \in \mathcal{S}_{\Sigma}(x)$ such that $\phi(0) = x$ and $\phi(t) = y$.

Now,  given $z \in y + \delta \mathbb{B}$, we consider the function $\eta : [0,t] \rightarrow \Delta \mathbb{B}$ given by
$$ \eta(s) := \phi(s) - \frac{s}{t} (y-z) \qquad \forall s \in [0,t]. $$
Note that
$$ \dot{\eta}(s) = \dot{\phi}(s) - \frac{1}{t} (y-z) \qquad \text{for almost all} \quad s \in [0,t]. $$
Next, for almost all $s \in [0,t]$, we let 
$f_\eta(s) \in F(\eta(s))$ be such that 
$ |f_\eta(s) - \dot{\phi}(s)| \leq \underline{\epsilon}/2. 
$
This is possible using \eqref{eqcontin}.  Hence, we conclude that 
\begin{align*} 
\dot{\eta}(s) & = f_\eta(s) + (\dot{\phi}(s) - f_\eta(s)) - \frac{1}{t} (y-z)  \in F(\eta(s)) + 
\left( \frac{\underline{\epsilon}}{2} + \frac{\delta}{t} \right) \mathbb{B} \qquad \quad \text{for almost all} ~~ s \in [0,t]. 
\end{align*} 
As a result, by taking  $\delta := \frac{\min \{ \underline{\epsilon}t, \Delta \}}{2}$,  we conclude that
$$
\dot{\eta}(s)  \in F(\eta(s)) + 
\underline{\epsilon} \mathbb{B} \subset F(\eta(s)) + 
\epsilon(\eta(s)) \mathbb{B} \qquad \text{for almost all} \quad s \in [0,t]. 
$$ 
Hence,  $\eta : [0,t] \rightarrow x + \Delta \mathbb{B}$ is a solution to 
$\Sigma_{\epsilon}$ verifying $\eta(0) = x$ and $\eta(t) = z$,  which proves \eqref{eqreachcover} since $z$ is arbitrary within $y + \delta \mathbb{B}$. 

The proof of \eqref{eqreachcover1} follows using the same steps used to prove \eqref{eqreachcover},  while invoking \eqref{eqcontin1b} instead of \eqref{eqcontin1a}. 
\end{proof}

It is important to note that,  when only  
Assumption \ref{ass1} is verified;  namely,  when $F$ is not required to be continuous,  then we can find a system $\Sigma$,  $x \in \mathbb{R}^n$,  and  $\epsilon \in \mathcal{C}_+$  such that 
$$R^b_{\Sigma}(t,x) \backslash \{x\}   \not\subset  \text{int}(R^b_{\Sigma_\epsilon}(t,x)) \qquad \forall t > 0.  $$ 

\begin{example} \label{exp1}
Consider system $\Sigma$ with $n = 2$ and 
$$ F(x) := \left\{ 
\begin{matrix} 
[0 \quad 1]^\top  & \text{if}~ x_2 < 0
\\
\co \{[0 \quad 1]^\top, [-1 \quad 0]^\top\} & \text{otherwise}.   
\end{matrix} 
\right. $$
Let the constant function $\epsilon : \mathbb{R}^2 \rightarrow \{1/2\}$.

Note that $\Sigma$ admits a solution $\phi$ starting from $x_o = 0$ that is given by 
$ \phi(t) := [-t \quad 0]^\top$ for all $t \geq 0$.  
Note that,  for each $t>0$,  we have 
$ \phi(t) \in R^b_{\Sigma}(t,x_o) \backslash \{x_o\} $. 

Now, given $t \geq 0$,  we show that
\begin{align*} 
 \phi(t) - [0 \quad 1/i]^\top \notin R_{\Sigma_\epsilon}(t,x_o) \qquad  \forall i \in \{1,2,... \}.  
 \end{align*}
To do so, it is enough to show that the set 
$K := \{x \in \mathbb{R}^2 : x_2 \geq 0 \}$ is forward invariant for $\Sigma_{\epsilon := \frac{1}{2}}$. 
We show the latter by applying  Lemma \ref{lemA14} in the Appendix,  via verifying 
\begin{align} \label{eqteng} 
F(x) + \mathbb{B}/2 \subset C_K(y)  \qquad  \forall y \in \text{Proj}_{K}(x) \qquad \forall  x \in \mathbb{R}^2 \backslash K. 
\end{align}
To show \eqref{eqteng},  
we note that 
$\mathbb{R}^2 \backslash K = \{ x \in \mathbb{R}^2 : x_2 < 0 \}$ and that the projection of $x := [x_1 \quad x_2]^\top \in \mathbb{R}^2 \backslash K$ on $K$ is $y := [x_1 \quad 0]^\top$. Furthermore,   using 
\cite[Lemma 3]{draftautomatica},  we conclude that 
$$ C_K(y) = \{v \in \mathbb{R}^2 : v_2 \geq 0 \} \qquad \forall y \in \partial K.  $$  
Now,  for each $x \in \mathbb{R}^2 \backslash K$,  we have  
$$ F(x) + \mathbb{B}/2 = 
\{[ \alpha/2 \quad (\beta/2) + 1 ]^\top : \sqrt{\alpha^2 + \beta^2} \leq 1 \}. $$ 
As a result, for each $v := [v_1 \quad v_2]^\top \in F(x) + 
\mathbb{B}/2$, we conclude that $v_2 \geq 0$, which means that $v \in C_K(y)$ for all $y \in \text{Proj}_{K}(x)$.   Hence, \eqref{eqteng} follows.
\end{example}
 
\subsection{The time to impact contractive and recurrent sets}

We start recalling a definition of forward contractivity of a closed subset $K$ for a differential inclusion $\Sigma$ \cite{Blanchini:1999:SPS:2235754.2236030}.  

\begin{definition}[Forward contractivity] 
A closed subset $K \subset \mathbb{R}^n$ is forward contractive for $\Sigma$ if it is forward invariant for $\Sigma$ and,  for every $x_o \in \partial K$ and for every  $\phi \in \mathcal{S}_{\Sigma}(x_o)$,  there exists $T > 0$ such that 
$$ \phi(t) \in \mbox{int}(K) \qquad  \forall t \in \dom \phi \cap (0,T].  $$
\end{definition}

\blue{Next, we recall a definition of global reccurence 
of a subset $K \subset \mathbb{R}^n$ for a differential inclusion $\Sigma$ }
\cite{SUBBARAMAN201654}. 
\begin{definition}[Global recurrence] 
\blue{A set $K \subset \mathbb{R}^n$ is globally recurrent for $\Sigma$ if,  for each $\phi \in \mathcal{S}_{\Sigma}(x)$ with  
$x \in \mathbb{R}^n$,  there exists $t \in \dom \phi \cap \mathbb{R}_{\geq 0}$ such that $\phi(t) \in \mbox{int}(K)$.} 
\end{definition}

\blue{Finally, we deduce a definition   of local recurrence,  which we will use in our work.}

\begin{definition}[Local recurrence]
A set $K \subset \mathbb{R}^n$ is locally recurrent for $\Sigma$ on a neighborhood of $\partial K$,  denoted by $U(\partial K)$,  
if, for each $\phi \in \mathcal{S}_{\Sigma}(x)$ with $x \in U(\partial K)$, there exists $t \in \dom \phi \cap \mathbb{R}_{\geq 0}$ such that 
$\phi(t) \in \mbox{int}(K)$. 
\end{definition}

\begin{remark}
Note that local (respectively, global) recurrence of the set $K \subset \mathbb{R}^n$ implies,  respectively,   
local (resp. global) recurrence of the set $\cl(K)$ and vice versa.
\end{remark}

\blue{ We recall here a variant of   the time-to-impact function (also known as the hitting-time function \cite{Aubin:1991:VT:120830}) with respect to a closed subset $K \subset \mathbb{R}^n$  along the solutions to $\Sigma$,  which we denote by $B_K : \mathbb{R}^n \rightarrow \mathbb{R}$.   We will show that such a function enjoys some regularity properties when $K$ is locally recurrent and forward contractive for $\Sigma$.   
We additionally show that $B_K$ is strictly decreasing along the solutions that lie within a neighborhood of 
$\partial K$.  }  

Before defining $B_K$,  we start introducing the   map 
$\hat{B}_{K} : U_1 \rightarrow  \mathbb{R} \cup \{\pm \infty\}$ 
given by 
\begin{align} \label{eqhatB}
\hat{B}_K(x) := 
\left\{ 
\begin{matrix}
\inf \{ T_{K} (\phi) : \phi \in \mathcal{S}_{\Sigma}(x) \} 
& \text{if $x \in \text{int}(K)$}
\\
0 & \text{if $x \in \partial K$}
\\
\sup \{ T_K (\phi) : \phi \in \mathcal{S}_{\Sigma}(x) \} & \text{otherwise},
\end{matrix} \right.
\end{align}
where,  for each $\phi \in \mathcal{S}_{\Sigma}(U_1)$,  
the functional $T_K : \mathcal{S}_{\Sigma}(U_1) \rightarrow \mathbb{R} \cup \{  \pm \infty \} $ is given by
$$ T_K (\phi) :=
\argmin \{ |t| : t \in \dom \phi,~ \phi(t) \in \partial K \}.  $$

\blue{ Roughly speaking,  $T_K(\phi)$ associates to each $\phi \in \mathcal{S}_{\Sigma}(x)$ the time, with the smallest norm,  at which $\phi$ hits $\partial K$.   
Then, $\hat{B}_K(x)$ is defined to be the smallest 
(respectively,   the largest) value among such times,  over all the $\phi$s in $\mathcal{S}_{\Sigma}(x)$,  when $x \in \text{int}(K)$ (respectively,  when $x \in \mathbb{R}^n \backslash K$).  }

\blue{ In \cite{Aubin:1991:VT:120830},  functionals similar to $T_K$ are introduced under the name hitting- and exit-time functionals,  denoted by $\theta_K$ and 
$\tau_K$, respectively.
Note that,  under Assumption \ref{assnew} below,  the latter two functionals coincide, and are related to $T_K$ through the relationship: 
$$ T_{K} (\phi) =  \theta_{\mathbb{R}^n \backslash K}(\phi) = \tau_{\mathbb{R}^n \backslash K}(\phi) \qquad \forall \phi \in \mathcal{S}_\Sigma(x) ~ \text{with} ~ x \in \mathbb{R}^n \backslash K.  $$
}

\begin{assumption} \label{assnew}
There exists a closed subset 
$U_{1} \subset \mathbb{R}^n$ such that $\partial K \subset \text{int}(U_1)$ and
\begin{enumerate} [label={A\ref{assnew}\arabic*)},leftmargin=*]
\item \label{item:71--} The set $K$ is locally recurrent for $\Sigma$ on $U_1$.
\item \label{item:72--} The set $\mathbb{R}^n\backslash K$ is locally recurrent for $\Sigma^-$ on $U_1$.
\item \label{item:73--} The set $K$ is forward contractive for $\Sigma$.
\item \label{item:74--} The set $\mathbb{R}^n \backslash \mbox{int}(K)$ is forward contractive for $\Sigma^-$.
\end{enumerate}
\end{assumption}

\blue{
Using \cite[Proposition 4.2.4]{Aubin:1991:VT:120830},  we can guarantee that 
$\hat{B}_K(x) = \sup\{T_K(\phi) : \phi \in \mathcal{S}_\Sigma(x)\}$
is upper semicontinuous on $\mathbb{R}^n \backslash K$ provided that Assumption \ref{ass1} holds and $F$ is strict Marchaud.  In the following lemma,  we establish the same conclusion,  among others,  using Assumptions \ref{ass1} and \ref{assnew} ($F$ is not required to be strict Marchaud). }

\begin{lemma} \label{lem4-}
Consider system $\Sigma$ such that Assumption \ref{ass1} holds.  \blue{Furthermore,  consider a closed subset $K \subset \mathbb{R}^n$ for which Assumption \ref{assnew} holds.}  
Then,  the following properties are satisfied.
\begin{enumerate} [label={P\ref{lem4-}\arabic*)},leftmargin=*]
\item \label{item:71-}   $\hat{B}_{K}$ is well defined on $U_1$; i.e.,  for each $x \in U_1$,  $\hat{B}_{K}(x)$ exists and is finite.
\item \label{item:72-}  $\hat{B}_K$ is upper 
semicontinuous on $U_1 \backslash \text{int}(K)$.  
\item \label{item:73-} $\hat{B}_K$ is lower semicontinuous on $K \cap U_1$.  
\end{enumerate}
\end{lemma}     

\begin{proof}
We first use \ref{item:71--} and \ref{item:72--} to conclude that each maximal solution 
$\phi \in \mathcal{S}_{\Sigma}(x)$ with $x \in U_1$ reaches $\partial K$ at some time, which can be either positive or negative.  Hence,  $T_{K}(\phi)$ exists and it is finite for all 
$\phi \in \mathcal{S}_{\Sigma}(U_1)$.   

Next,  we will show that the map $x \mapsto T_K(\mathcal{S}_\Sigma(x))$ is locally bounded on $U_1 \backslash K$.   
To find a contradiction,   we pick $x_o \in U_1 \backslash K$, a sequence 
$\{ x_i \}^{\infty}_{i=1} \subset  U_1 \backslash K$ such that  $ \lim_{i \rightarrow \infty} x_i = x_o$, and a sequence 
$\{ \phi_i \}^{\infty}_{i=1}$ of solutions (not necessarily maximal) to $\Sigma$ such that $ \phi_i(0) = x_i$,   $\dom \phi_i = [0,T_{K}(\phi_i)]$ for all $i \in \{1,2,...\}$, and $ \lim_{i \rightarrow \infty} T_{K} (\phi_i) = +\infty$.   Without loss of generality,  we also assume that $i \mapsto T_K(\phi_i)$ is strictly increasing. 
Since $F$ is locally bounded,  using Lemma \ref{lemprepre} in the Appendix,  we conclude the existence of $\bar{T} \in \mathbb{R}_{>0} \cup \{+ \infty\}$ the largest time  such that,  on any interval $[0,T] \subset [0,\bar{T})$,   the sequence $\{ \phi_i \}^{\infty}_{i=0}$ is uniformly bounded.   As a result,  by passing to an appropriate subsequence and \blue{using \cite[Theorem 5.29]{goebel2012hybrid}} recursively on each closed interval $[0,T] \subset [0,\bar{T})$,  we conclude the existence of $\phi : [0,\bar{T}) \rightarrow \mathbb{R}^n$   solution  to $\Sigma$ starting from $x_o$ such that 
\begin{align} \label{eqlimit-}
\text{lim}_{i \rightarrow \infty} \phi_i(t) = \phi(t) \qquad \forall t \in [0,\bar{T}). 
\end{align}
Now,  since $x_o \in U_1 \backslash K$,  we conclude the existence of $\alpha \in (0, \bar{T})$ such that 
$T_{K} (\phi) < \alpha$. \blue{Indeed, $\alpha$ must be smaller than  $\bar{T}$ because  $T_K(\phi)$ exists and is finite, and when $\bar{T}$ is finite, we necessarily have  $\lim_{t \rightarrow \bar{T}} \phi(t) = + \infty$.}
Next, by forward contractivity of 
the set $K$ for $\Sigma$,  we conclude that 
$\phi(\alpha) \in \mbox{int}(K)$; thus, there exists 
$\beta>0$ such that 
$ | \phi(\alpha) |_{\mathbb{R}^n \backslash K}  \geq \beta$.
However, there exists $i^* \in \mathbb{N}$ such that 
$ \phi_i (\alpha) \notin K$ for all $i \geq i^*$,  
which implies that $ |\phi(\alpha) - \phi_i(\alpha)| \geq \beta$ for all $i \geq i^*$.    The latter contradicts \eqref{eqlimit-}; hence, the map $x \mapsto T_K(\mathcal{S}_\Sigma(x))$ is locally bounded and thus \ref{item:71-} follows as a direct consequence. 

To  prove \ref{item:72-},  we start re-expressing $\hat{B}$ as
$$  \hat{B}(x) :=
 \sup \{ t_\phi : t_\phi \in T_{K} (\mathcal{S}_{\Sigma}(x)) \} \qquad \forall x \in U_1 \backslash \text{int}(K). 
$$
Furthermore, we propose to show that the set-valued map 
$x \mapsto T_{K} (\mathcal{S}_{\Sigma}(x))$ is outer 
semicontinuous on $U_1 \backslash \text{int}(K)$. For this, we consider a sequence $\{x_i\}^\infty_{i=1} \subset U_1 \backslash \text{int}(K)$ that converges to $x_o \in U_1 \backslash \text{int}(K)$ and a sequence $\{t_i\}^\infty_{i=1} \subset \mathbb{R}_{\geq 0}$ that converges to $t_o \in \mathbb{R}_{\geq 0}$ such that 
$ t_i \in T_{K} (\mathcal{S}_{\Sigma}(x_i))$  for all i $\in \{1,2, ...\}$,   and we show that $t_o \in T_{K} (\mathcal{S}_{\Sigma}(x_o))$.  
Indeed,   let a sequence of solutions $ \{ \phi_i \}^{\infty}_{i =1}$,  with $\dom \phi_i = [0,  t_i]$, $t_i = T_K(\phi_i)$,  and $\phi_i(0) = x_i$ for all $i \in \{1,2,...\}$.  Here,  we distinguish between two cases.   
\begin{itemize}
\item When the sequence  $\{ \phi_i \}^{\infty}_{i =1}$ is uniformly bounded,  using \cite[Theorem 5.29.]{goebel2012hybrid}, we conclude the existence of a solution $\phi : [0,t_o] \rightarrow \mathbb{R}^n$ that is the graphical limit of an appropriate subsequence of  $\{ \phi_i \}^{\infty}_{i =1}$.  Next,  we show that 
$t_o = T_K(\phi) $.  To find a contradiction,  we assume that  $t_o \neq  T_K(\phi)$. Hence,  using \ref{item:71--} and \ref{item:73--},  it follows that there exists $\delta > 0$ such that either $|\phi(t_o)|_K \geq \delta$ or 
$|\phi(t_o)|_{\mathbb{R}^n \backslash K} \geq \delta$. 
On the other hand,  \blue{by graphical convergence and continuity of solutions}, we conclude that
$ \lim_{i \rightarrow \infty} \phi_i(t_i) = \phi(t_o)$.  Hence,  there exists $i^\star \in \{1,2,...\}$ such that $|\phi_i(t_i) - \phi(t_o)| \leq \delta/2$ for all $i \geq i^\star$.  The latter implies that 
$|\phi_i(t_i)|_{\partial K} \geq \delta/2$  for all $i \geq i^\star$,   and the contradiction follows. 

\item If the sequence is not uniformly bounded, we   exploit local boundedness of $F$ and Lemma \ref{lemprepre} in the Appendix,  to conclude the existence of $\bar{T} \in (0,  t_o]$ the largest time such that the sequence $\{ \phi_i \}^{\infty}_{i =1}$ is uniformly bounded on any $[0, T] \subset [0,\bar{T})$.  
\blue{Using \cite[Theorem 5.29]{goebel2012hybrid} recursively on each interval $[0, T] \subset [0, \bar{T})$},  we conclude the existence of an unbounded  solution $\phi : [0,\bar{T}) \rightarrow \mathbb{R}^n$  such that,  after passing to an appropriate subsequence,  we obtain  
 \begin{align} \label{eqlimit}  
 \lim_{i \rightarrow \infty} \phi_i(t)  = \phi(t) \qquad \forall t \in [0,\bar{T}).   
 \end{align}
 In this case,  we conclude that $\phi$ must reach $\partial K$ at some $ \alpha := T_K(\phi) \in [0, \bar{T})$.  Since $K$ is contractive,  then we can find 
$\beta \in (\alpha,  t_o)$ such that 
$|\phi(\beta)|_{\mathbb{R}^n \backslash K} = \gamma$ for some $\gamma > 0$.  Now,  using \eqref{eqlimit}, we conclude that $ \lim_{i \rightarrow \infty} |\phi_i(\beta)|_{\mathbb{R}^n \backslash K}  = \gamma  > 0$ with 
$\beta < t_o$; hence, $\beta < t_i$ for $i$ large enough. This yields to a contradiction implying that $\alpha = t_o$.   
\end{itemize}

Having the map $x \mapsto T_K(\mathcal{S}_\Sigma(x))$ locally bounded and  outer semicontinuous on $U_1 \backslash \text{int}(K)$ implies that the map $x \mapsto T_K(\mathcal{S}_\Sigma(x))$ is upper semicontinuous with compact images on $U_1 \backslash \text{int}(K)$; see Remark \ref{remplus}.   Therefore,  using  \cite[Theorem 1.4.16]{aubin2009set},  we conclude that $\hat{B}_{K}$ is upper semicontinuous on $ U_1 \backslash \text{int}(K)$.

Finally,  to prove \ref{item:73-},  it is enough to 
notice that   
\begin{align*}
\inf \{ T_{K} (\phi) : \phi \in \mathcal{S}_{\Sigma}(x) \}  & =  - \sup\{ - T_{K} (\phi) : \phi \in \mathcal{S}_{\Sigma}(x) \}  =  - \sup\{  T_{K} (\phi) : \phi \in \mathcal{S}_{\Sigma^-}(x) \}.
\end{align*}
Therefore,  the previous arguments apply,  under \ref{item:72--} and \ref{item:74--},   to conclude that $-\hat{B}_K$ is upper semicontinuous on $K \cap U_1$; thus, 
 $\hat{B}_K$ is lower semicontinuous on $K \cap U_1$.
\end{proof}

Next, we introduce the function $B_K : \mathbb{R}^n \rightarrow \mathbb{R}$, which extends $\hat{B}_K$ to the entire $\mathbb{R}^n$ and is given by
\begin{equation}
\label{eqbarcand}
\begin{aligned} 
\hspace{-0.3cm} B_K(x) := 
\left\{
\begin{matrix}
\inf \{ \hat{B}_K(y) : y \in \text{Proj}_{U_1}(x)  \}   & \hspace{-0.2cm} \text{if} ~ x \in \text{int}(K)
\\
0 & \hspace{-0.2cm} \text{if $x \in \partial K$}
\\
\sup \{ \hat{B}_K(y) : y \in \text{Proj}_{U_1}(x) \} & \hspace{-0.2cm} \text{otherwise}.
\end{matrix}
\right. 
\end{aligned}
\end{equation}

\begin{lemma} \label{lem5-}
Consider system $\Sigma$ such that Assumption \ref{ass1} holds.  \blue{Furthermore,  consider $(X_o,X_u) \subset \mathbb{R}^n \times \mathbb{R}^n$ and a closed subset $K \subset \mathbb{R}^n$  such that $X_o \subset K$, $X_u \cap K = \emptyset$, and 
Assumption \ref{assnew} holds}.   

Then, the following properties are satisfied.  
\begin{enumerate} [label={P\ref{lem5-}\arabic*)},leftmargin=*]
\item \label{item:B1--}  $B_K$ is a barrier function candidate with respect to $(X_o,X_u)$ and $K = \{ x \in \mathbb{R}^n : B_K(x) \leq 0\}$. 
\item \label{item:B3--} $B_K$ is upper semicontinuous on $\mathbb{R}^n \backslash \text{int}(K)$.
\item \label{item:B4--} $B_K$ is lower semicontinuous on $K$.
\item \label{item:B5--} Along each solution $\phi$ to $\Sigma$ (not necessarily maximal) such that $\phi(\dom \phi) \subset U_1$,  the map $t \mapsto B_{K}(\phi(t))$ is  decreasing. 
In particular, we have 
$$ B_{K}(\phi(t')) - B_{K}(\phi(t)) \leq  t-t'  \quad \forall t, t' \in \dom \phi  ~~ \text{with} ~~ t \leq t'.  $$
\end{enumerate}
\end{lemma}

\begin{proof}
To prove \ref{item:B1--},  we use 
the fact  that the set $K$ is forward contractive for $\Sigma$. Hence, the solutions to $\Sigma$ starting outside the set $K$ cannot reach $K$ at negative times. This implies that $ B_K(x) > 0$ for all $x \in \mathbb{R}^n \backslash K$.   Similarly,  the set $\mathbb{R}^n \backslash 
\mbox{int}(K)$ is forward contractive for 
$\Sigma^-$. Hence, the solutions to $\Sigma$ starting from $K$ cannot reach $\mathbb{R}^n \backslash \mbox{int}(K)$ at  positive times. This implies that $B_K(x) \leq 0$ for all $x \in K$. 

To prove \ref{item:B3--}, we start noting that 
\begin{align} \label{eqBBhat}
B_{K}(x) = \hat{B}_{K}(x)  \qquad \forall x \in U_1, 
\end{align}

which implies, using \ref{item:72--}, that $B_K$ is upper semicontinuous on the set $U_1 \backslash K$. Furthermore,  using \cite[Example 5.23]{rockafellar2009variational},  we conclude that $\text{Proj}_{U_1}$ is outer semicontinuous and locally bounded; hence, admits compact images. 
Thus,   using \cite[Theorem 1.4.16]{aubin2009set},  \ref{item:B3--} follows.   

Similarly,  to prove \ref{item:B4--}, we use \eqref{eqBBhat} and 
\ref{item:73-} to conclude that $B_K$ is lower semicontinuous on $U_1 \cap K$.    Furthermore, for each $x \in K \backslash U_1$,  we have 
$ B_K(x) := - \sup \{ -\hat{B}_K(y) : y \in \text{Proj}_{U_1}(x) \}.   $
Hence,  using \cite[Example 5.23]{rockafellar2009variational},  we know that $\text{Proj}_{U_1}$ is outer semicontinuous and locally bounded; thus,  admits compact images.  Thus,  using \cite[Theorem 1.4.16]{aubin2009set},  we conclude that $-B_K$ is upper semicontinuous; thus, $B_K$ is lower semicontinuous on $K$.

Finally, to prove \ref{item:B5--}, we let $\phi$ be a solution to $\Sigma$ (not necessarily maximal) starting from $x_o \in U_1 $ and such that $\phi(\dom \phi) \subset U_1$. Furthermore, we let
$t \in \dom \phi$ and $\sigma>0$ such that $[t,t+ \sigma] \subset \dom \phi$.  Hence,
$$ B_{K}(\phi(s)) =  \hat{B}_{K}(\phi(s)) \qquad \forall s \in [t,t+\sigma]. $$  
Next, we distinguish three situations. 
\begin{enumerate} [label={\arabic*)},leftmargin=*]
\item  When $\phi([t,t+\sigma]) \subset 
U_1 \backslash K$, there exists a solution $\psi \in \mathcal{S}_{\Sigma}(\phi(t+
\sigma))$ such that
$ B_K(\phi(t+\sigma))=T_K(\psi)$. Additionally, we introduce the solution 
$\hat{\psi} \in \mathcal{S}_{\Sigma}(\phi(t))$ satisfying
\begin{align*} 
\hat{\psi}(s) = 
\left\{
\begin{matrix}
\psi (s-\sigma) & \forall s \geq \sigma, \\
\phi(t+s) &  \forall s \in [0, \sigma]. \end{matrix} \right.
\end{align*}
Then, $T_{K}(\hat{\psi}) = \sigma + T_{K}(\psi)$. 
Hence, \eqref{eqhatB} implies that
$$ B_{K}(\phi(t)) \geq T_{K}(\hat{\psi}) = \sigma + B_{K}(\phi(t+\sigma)). $$

\item  When $\phi([t,t+\sigma]) \subset U_1 \cap \mbox{int}(K)$,
 there exists a solution 
 $\psi \in \mathcal{S}_{\Sigma}(\phi(t))$ 
such that $B_{K}(\phi(t))=T_{K}(\psi)$. 
Additionally, we introduce the solution $\hat{\psi} \in \mathcal{S}_{\Sigma}(\phi(t+
\sigma))$ satisfying
\begin{align*} 
\hat{\psi}(s) = 
\left\{ 
\begin{matrix}
\psi(s+\sigma) & \forall s \leq - \sigma, \\
 \phi(t+\sigma+s) & \forall s \in [- \sigma, 0].
\end{matrix} \right.
\end{align*}
As a result, we have $T_{K}(\hat{\psi}) = -\sigma + T_{K}(\psi)$ and \eqref{eqhatB} implies that
$$ B_{K}(\phi(t+\sigma)) 
\leq T_{K}(\hat{\psi}) = - \sigma + 
B_{K}(\phi(t)). $$

\item Finally, we consider a solution 
$\phi$ satisfying 
$ \phi(t) \in U_1  \backslash K$,  $\phi(t+\sigma) \in  U_1 \cap \mbox{int}(K)$,  
and such that there exists a unique $t_1 \in (t, t+\sigma)$ for which  
$ \phi(t_1) \in \partial K$.   Indeed, the aforementioned scenario complement the scenarios in the previous two items since, due to \ref{item:73--} and \ref{item:74--}, the solutions cannot slide on $\partial K$.  

Now, as in the previous steps, we conclude that 
$$ B_{K}(\phi(t+t'_1)) - B_{K}(\phi(t)) \leq -  t'_1 \qquad \forall t'_1 \in [0,t_1) $$
and
$$
B_{K}(\phi(t+\sigma)) - B_{K}(\phi(t + t'_1)) 
\leq  - (\sigma -  t'_1) \qquad \forall t'_1 \in (t_1,\sigma]. 
$$
Next, we note that
$$ B_{K}(\phi(t+t'_1)) > 0 = B_{K}(\phi(t+t_1))  \qquad \forall t'_1 \in [0,t_1) $$
and 
$$
 B_{K}(\phi(t+t'_1)) < 0 = 
 B_{K}(\phi(t+t_1)) \qquad \forall t'_1 \in (t_1,\sigma]. $$
Hence,  
$$ B_K(\phi(t+t_1)) - B_K(\phi(t)) \leq - t_1 \quad \text{and} 
\quad 
B_K(\phi(t+\sigma)) - B_K(\phi(t + t_1))  
\leq  - (\sigma -  t_1),
$$
 implying that  
$ B_K(\phi(t+\sigma)) - B_K(\phi(t))  
\leq  - \sigma. $  
\end{enumerate}
\end{proof}

\section{Proof of Theorem \ref{thm3}} 
\label{Sec.5}

\subsection{Proof steps} 

We propose to prove Theorem \ref{thm3}
by following three steps. 

\begin{itemize}
\item[\underline{Step 1}] Given $\bar{\epsilon} \in \mathcal{C}_+$,  we introduce the set 
\begin{align} \label{eqinfreachbis}
K_{\bar{\epsilon}} := \bigcup_{t \geq 0} \bigcup_{x \in X_o} R_{\Sigma_{\bar{\epsilon}}} (t, x).
\end{align}
 When $\Sigma$ is robustly safe with respect to $(X_o,X_u)$ and Assumption \ref{ass4-} holds,  we show that
\begin{align} \label{eqdisconnect-}
X_u \cap \cl(K_{\bar{\epsilon}}) = \emptyset, 
\end{align}
for appropriately chosen robustness margin $\bar{\epsilon}$.  In particular,   when additionally Assumption \ref{ass4} holds, we show that 
\begin{align} \label{eqinter}
\cl(X_u) \cap \cl(K_{\bar{\epsilon}}) = \emptyset.
\end{align}

\begin{lemma}  \label{lem2bis}
Consider system $\Sigma$ such that Assumption \ref{ass3} holds.  Consider $(X_o,X_u) \subset \mathbb{R}^n \times \mathbb{R}^n$ such that $\Sigma$ is robustly safe with respect to $(X_o,X_u)$. Then, given a robustness margin $\bar{\epsilon}_o$, we conclude that, for each  $\bar{\epsilon} \in \mathcal{C}_+$ satisfying 
\begin{align*} 
\bar{\epsilon}(x) < \bar{\epsilon}_o(x) \qquad \forall x \in \mathbb{R}^n,
\end{align*}
 the following properties hold. 
\begin{enumerate} [label={\arabic*)},leftmargin=*]
\item \label{item:P52b} If Assumption \ref{ass4} holds, then \eqref{eqinter} holds. 
\item \label{item:P52a} If Assumption \ref{ass4-} holds, then \eqref{eqdisconnect-} holds.
\end{enumerate} 
\end{lemma}

\item[\underline{Step 2}] It is shown in \cite[Lemma 2]{ratschan2018converse} that,  given $\bar{\epsilon} \in \mathcal{C}_+$,  when  
$\cl(K_{\bar{\epsilon}})$ is bounded and $F$ is single valued and smooth,  every maximal solution to $\Sigma$  starting from $R_{\Sigma}(t,  \partial K_{\bar{\epsilon}})$,   for some $t \in \mathbb{R}$,   must cross $\partial K_{\bar{\epsilon}}$ only one time. Motivated by this fact, we establish contractivity of  
$\cl(K_{\bar{\epsilon}})$ for $\Sigma_\epsilon$, for an appropriate choice of the perturbation $\epsilon$. 

\begin{lemma} [Contarctivity of the set 
$\cl(K_{\bar{\epsilon}})$]  \label{lem2}
Consider system $\Sigma$  such that Assumption \ref{ass3} holds.   
Consider $X_o \subset\mathbb{R}^n$,   $\bar{\epsilon} \in \mathcal{C}_+$,  and the set $K_{\bar{\epsilon}}$ introduced in \eqref{eqinfreachbis}.  Then,  for each  $\epsilon \in \mathcal{C}_+$ satisfying  
\begin{align} \label{eqepsineq}
\epsilon(x) < \bar{\epsilon}(x) \qquad \forall x \in \mathbb{R}^n,
\end{align}
the following properties hold.
\begin{enumerate} 
[label={P\ref{lem2}\arabic*)},leftmargin=*]
\item \label{item:P3a} The set 
$\cl(K_{\bar{\epsilon}})$ is forward contractive for $\Sigma_{\epsilon}$.
\item \label{item:P3b} The set $\mathbb{R}^n \backslash \mbox{int}(K_{\bar{\epsilon}})$ is forward contractive for $\Sigma^-_{\epsilon}$.
\item \label{item:P3c}  
The solutions 
to $\Sigma_\epsilon$ starting from $\text{int}(K_{\bar{\epsilon}})$ never reach 
$\partial K_{\bar{\epsilon}}$ for positive times.
\item \label{item:P3d} 
The solutions to $\Sigma_\epsilon$ starting from $\mathbb{R}^n \backslash \cl(K_{\bar{\epsilon}})$ never reach 
$\partial K_{\bar{\epsilon}}$ for negative times.
\end{enumerate}
\end{lemma}

\item[\underline{Step 3}] 
It is shown in  \cite[Lemma 3]{ratschan2018converse} that,  when the set $\cl(K_{\bar{\epsilon}})$ is bounded and when $F$ is single-valued and smooth,  there exists a neighborhood of $\partial K_{\bar{\epsilon}}$,  denoted  $U \subset \mathbb{R}^n$,  such that every maximal solution to $\Sigma$ starting from $U$ reaches $\partial K_{\bar{\epsilon}}$ in finite (positive or negative) time.   
Inspired by this observation and using Lemma \ref{lem2},  we establish local recurrence of the set $K_{\bar{\epsilon}}$ on a neighborhood of $\partial K_{\bar{\epsilon}}$.

\begin{lemma} [Recurrence of the set 
$\cl(K_{\bar{\epsilon}})$] \label{lem3}
Consider system $\Sigma$ such that Assumption  \ref{ass3} holds.  
Consider $X_o \subset \mathbb{R}^n$,  $\epsilon_1, \bar{\epsilon} \in \mathcal{C}_+$ such that 
\begin{align} \label{eqepseps1}
\epsilon_1(x) < \bar{\epsilon}(x) \qquad \forall x \in \mathbb{R}^n.
\end{align}
Then, there exists a closed subset $U_{1} \subset \mathbb{R}^n$ such that 
\begin{align} \label{eqinter0}
\partial K_{\bar{\epsilon}} \subset \text{int} (U_1),
\end{align}
where the set $K_{\bar{\epsilon}}$ is introduced in \eqref{eqinfreachbis}.  Moreover,   for each 
$\epsilon \in \mathcal{C}_+$ satisfying  
\begin{align} \label{eqepsilon2}
\epsilon(x) \leq \epsilon_1(x) \qquad \forall x \in \mathbb{R}^n, 
\end{align}
the following properties hold:

\begin{enumerate} [label={P\ref{lem3}\arabic*)},leftmargin=*]
\item \label{item:R1} The set $K_{\bar{\epsilon}}$ is locally recurrent for $\Sigma_{\epsilon}$ on $U_1$.
\item \label{item:R2} The set 
$\mathbb{R}^n\backslash K_{\bar{\epsilon}}$ is locally recurrent for $\Sigma^-_{\epsilon}$ on $U_1$.
\end{enumerate}
\end{lemma}
\end{itemize}

\blue{ As a consequence,  when $\Sigma$ is robustly safe with respect to $(X_o,X_u)$ and the assumptions in Theorem \ref{thm3} hold,  using Lemmas \ref{lem2bis}, \ref{lem2}, and \ref{lem3},  we are able to prove the existence of robust-safety margins $\bar{\epsilon}$ and $\epsilon_1$ satisfying \eqref{eqepseps1} such that,  for each $\epsilon  \in \mathcal{C}_+$ satisfying \eqref{eqepsilon2}, 
 Assumption \ref{assnew}  is verified while replacing 
 $(\Sigma,K)$ therein by 
 $(\Sigma_\epsilon,\cl(K_{\bar{\epsilon}}))$,  
 $X_o \subset \cl(K_{\bar{\epsilon}})$,  and 
$X_u \cap \cl(K_{\bar{\epsilon}}) = \emptyset$.  
Hence,   applying Lemma \ref{lem5-},  we conclude that \ref{item:C2bis} is verified for $B = B_{\cl(K_{\bar{\epsilon}})}$ as defined in \eqref{eqhatB} and \eqref{eqbarcand} while replacing 
$(\Sigma, K)$ therein by 
$(\Sigma_{\epsilon_1},  \cl(K_{\bar{\epsilon}}))$ and for any $\epsilon  \in \mathcal{C}_+$ satisfying \eqref{eqepsineq}.  }

\subsection{Proof of Lemma \ref{lem2bis} }  \label{Sec.Separ}

We prove \ref{item:P52b} using contradiction. That is, when \eqref{eqinter} is not satisfied, we conclude that there exists $x \in \partial X_u \cap \partial K_{\bar{\epsilon}}$ and such that $x \in \partial K_{\bar{\epsilon}_o}$. The latter is the only possibility since $\bar{\epsilon}_o$ is a robustness margin and since 
$K_{\bar{\epsilon}} \subset K_{\bar{\epsilon}_o}$.  Furthermore, using Assumption \ref{ass4}, we conclude that $x \notin \cl(X_o)$. Next, we consider a sequence $\{x_i\}^{\infty}_{i = 0} \subset K_{\bar{\epsilon}}$ 
that converges to $x$.  For 
$ \underline{\epsilon} := 
\min \{ \bar{\epsilon}_o(y)-\bar{\epsilon}(y) : 
y \in x + \mathbb{B} \}, $
we pick $\Delta \in (0,1)$ such that 
$ (x + \Delta \mathbb{B}) \cap \cl (X_o) = \emptyset $
and, for each 
$(x_1,x_2) \in (x + \Delta \mathbb{B}) \times (x + \Delta \mathbb{B})$ and for each $\eta_1 \in F(x_1)$, there exists $\eta_2 \in F(x_2)$ such that
\begin{align} \label{eqAdded}
|\eta_1 - \eta_2| \leq \frac{\underline{\epsilon}}{2}. 
\end{align}  
The latter is possible since the set-valued map $F$ is assumed to be continuous.
Without loss of generality, we assume that 
\begin{align} \label{eqstep1}
x_i \in x + \frac{\Delta}{4} \mathbb{B} \qquad \forall i \in \mathbb{N}. 
\end{align}
Furthermore, we consider a sequence of solutions 
$\{\phi_i\}^{\infty}_{i = 0}$ to $\Sigma_{\bar{\epsilon}}$ such that each solution $\phi_i$ starts from $x_{oi} \in \cl(X_o)$. Moreover, for each $i \in \mathbb{N}$, there exist $t^1_i > t^o_i > 0$ such that
$ \phi_i([t^o_i,t^1_i]) \subset x + \frac{\Delta}{2} \mathbb{B}$,   $\phi_i([t^o_i,t^1_i]) \subset K_{\bar{\epsilon}}, $ 
$ 
|\phi_i(t^o_i) - x| = \frac{\Delta}{2}$,  and
$\phi_i(t^1_i) = x_i$.   
Next, using \eqref{eqstep1} and local boundedness of the map $F + \bar{\epsilon} \mathbb{B}$, we conclude the existence of $t > 0$ such that 
\begin{align} \label{eqDelta+i}
t^1_i - t^o_i \geq t \qquad \forall i \in \mathbb{N}. 
\end{align}
Now, we let $\tilde{\Sigma} := \Sigma_{\bar{\epsilon}}$; which implies that  
$\tilde{\Sigma}_{\epsilon := \bar{\epsilon}_o - \bar{\epsilon}} = \Sigma_{\bar{\epsilon}_o}$. 
Furthermore, using a similar approach as in  the proof of Lemma \ref{lem1},  we  will show the existence of $\delta>0$ such that
\begin{align} \label{eqreachreach}
 \phi_i(t^1_i) + \delta \mathbb{B} \subset R^b_{\tilde{\Sigma}_{\epsilon}}(t^1_i - t^o_i,\phi_i(t^o_i)) \quad \forall i \in \mathbb{N}. 
\end{align}    
Indeed, let $\delta \in (0, \Delta/4]$ and note that 
$ x_i \in 
R^b_{\tilde{\Sigma}}(t^1_i - t^o_i,y_i) 
\backslash \{ y_i \}$,   where $y_i := \phi_i(t^o_i)$. 

Now, given $z \in x_i + \delta \mathbb{B}$, we consider the function 
$$ \eta_{zi}(s) := \phi_i(s+t^o_{i}) - \frac{s}{t^1_i-t^o_i} (x_i-z) \quad \forall s \in [0,t^1_i-t^o_i]. $$
Note that 
$ \eta_{zi}(s) \in x + \Delta \mathbb{B}$ for all 
$s \in [0,t^1_i-t^o_i]$, and for all $i \in \mathbb{N}$.   Furthermore, for almost all $s \in [0,t^1_i-t^o_i]$, we have 
$$ \dot{\eta}_{zi}(s) := \dot{\phi}_i(s + t^o_i) - \frac{1}{t^1_i-t^o_i} (x_i-z). $$
Next, for almost all 
$s \in [0, t^1_i-t^o_i]$, we let 
$\eta_{1i}(s) \in F(\eta_{zi}(s))$ be such that 
$|\eta_{1i}(s) - \dot{\phi}_i(s + t^o_i)| \leq \underline{\epsilon}/2$.  The latter is possible using \eqref{eqAdded}. Hence, for almost all $s \in [0,t^1_i-t^o_i]$, we have  
\begin{align*} 
\dot{\eta}_{zi}(s) & := \eta_{1i}(s) + (\dot{\phi}_i(s + t^o_i) - \eta_{1i}(s)) - \frac{1}{t^1_i-t^o_i} (x_i-z) \in F(\eta_{zi}(s)) + 
\left( \frac{\underline{\epsilon}}{2} + \frac{\delta}{t^1_i-t^o_i} \right) \mathbb{B}. 
\end{align*} 
Hence, by taking  $\delta := \min \{ \underline{\epsilon} t, \Delta \}/4,$
where $t$ comes from \eqref{eqDelta+i}, we conclude that, for almost all 
$s \in [0, t^1_i-t^o_i]$,
\begin{align*} 
\dot{\eta}_{zi}(s) & \in F(\eta_{zi}(s)) + 
\underline{\epsilon} \mathbb{B} \subset F(\eta_{zi}(s)) + 
\epsilon(\eta_{zi}(s)) \mathbb{B}. 
\end{align*}
Hence, $\eta_{zi} : 
[0,t^1_i-t^o_i] \rightarrow x + \Delta \mathbb{B}$ is a solution to 
$\tilde{\Sigma}_{\epsilon}$  with 
$\eta_{zi}(0) = y_i$ and $\eta_{zi}(t^1_i-t^o_i) = z$, which proves \eqref{eqreachreach}. 

Next, since $\delta$ is uniform in $i$, we conclude the existence of $i \in \mathbb{N}$ sufficiently large such that 
$$ x \in \mbox{int} \left( R^b_{\tilde{\Sigma}_{\epsilon}}(t^1_i - t^o_i,y_i) \right). $$    Hence, since $x \in \partial K_{\bar{\epsilon}_o}$, we conclude that there exists 
$y \notin \cl(K_{\bar{\epsilon}_o})$ such that $ y \in  R^b_{\tilde{\Sigma}_{\epsilon}}(t^1_i - t^o_i,y_i)$.  
However, $y_i = \phi_i(t^o_i) \in K_{\bar{\epsilon}} \subset K_{\bar{\epsilon}_o}$ and the latter set is, by definition, forward invariant for $\tilde{\Sigma}_{\epsilon} = \Sigma_{\bar{\epsilon}_o}$, which yields to a contradiction.  

Finally, we prove \ref{item:P52a}) using contradiction. That is, when \eqref{eqdisconnect-} is not satisfied, we conclude that there exists 
$x \in \partial X_u \cap 
\partial K_{\bar{\epsilon}} \cap X_u$ such that $x \in \partial K_{\bar{\epsilon}_o}$. The latter is the only possibility since $\bar{\epsilon}_o$ is a robustness margin and $K_{\bar{\epsilon}} \subset K_{\bar{\epsilon}_o}$ by definition. 
Furthermore, we distinguish two situations:

\begin{itemize}
\item When $x \notin \cl(X_o)$, the contradiction follows using the same steps as in the proof of \ref{item:P52b}).

\item When $x \in \cl(X_o)$, the contradiction follows using Assumption \ref{ass4-}.
\end{itemize}

\subsection{Proof of Lemma \ref{lem2}} 
 
We prove \ref{item:P3a} by first proving that the set 
$\cl(K_{\bar{\epsilon}})$ is forward invariant for $\Sigma_{\epsilon}$. To find a contradiction, we assume the existence of $x \in \partial K_{\bar{\epsilon}}$ and  $\phi \in \mathcal{S}_{\Sigma_{\epsilon}}(x)$ such that, for some $T_1 >0$,
$$ \phi(s) \in \mathbb{R}^n \backslash \cl(K_{\bar{\epsilon}}) \qquad  \forall s \in (0,T_1].  $$   Next, we take $t \in (0, \min\{T_1,T\}]$,   \blue{where $T>0$ is such that, for each $t \in (0,T]$,  there exists $\delta>0$ such that
$ x + \delta \mathbb{B} \subset 
 R^b_{\Sigma_{\bar{\epsilon}}}(-t, \phi(t))$; see Lemma \ref{lem1}}.  Hence, there exists 
$y \in \mbox{int} (K_{\bar{\epsilon}})$ and a solution $\psi$ to $\Sigma_{\bar{\epsilon}}$ such that $\psi(0) = y$ and  $\psi(t) = \phi(t)$. 
The latter implies that the set 
$K_{\bar{\epsilon}}$ is not forward invariant for 
$\Sigma_{\bar{\epsilon}}$, which yields to a contradiction. 

To complete the proof, we show that the solutions to $\Sigma_{\epsilon}$ cannot slide on $\partial K_{\bar{\epsilon}}$; namely, for each $x \in \partial K_{\bar{\epsilon}}$ and for each 
$\phi \in \mathcal{S}_{\Sigma_{\epsilon}}(x)$, there exists $T>0$ such that 
\begin{align*} 
\phi(t) \in \mbox{int} (K_{\bar{\epsilon}}) \qquad \forall t \in (0,T].  
\end{align*}
To find a contradiction,  we assume the existence of a solution $\phi$ to $\Sigma_{\epsilon}$ starting from $x_o \in \mathbb{R}^n$ and an interval $(t_1, t_2) \subset \dom \phi$,  with $t_2 > t_1$,  such that 
$ \phi(s) \in \partial K_{\bar{\epsilon}}$ for all $s \in (t_1, t_2)$. 
Next, we pick $t_3 \in (t_1, t_2)$ and we take $t \in (0,\min\{(t_2 - t_3),T\})$, 
\blue{where $T>0$ is such that, for each $t \in (0,T]$,  there exists $\delta>0$ such that
\begin{align} \label{eqball+} 
\phi( t + t_3) + \delta \mathbb{B} \subset 
 R^b_{\Sigma^-_{\bar{\epsilon}}}(t, \phi(t_3));
 \end{align}
  see Lemma \ref{lem1}}. 
Similarly,   we take $t \in (0,\min\{(t_3 - t_1),T\})$, 
\blue{where $T>0$ is such that, for each $t \in (0,T]$,  there exists $\delta>0$ such that
\begin{align} \label{eqball-} 
\phi( -t + t_3) + \delta \mathbb{B} \subset 
R^b_{\Sigma^-_{\bar{\epsilon}}}(t, \phi(t_3));
\end{align}
see Lemma \ref{lem1}}.   Finally, combining \eqref{eqball+} and \eqref{eqball-}, we conclude the existence of a solution to $\Sigma_{\bar{\epsilon}}$ starting from $\mbox{int} (K_{\bar{\epsilon}})$ and that leaves the set $K_{\bar{\epsilon}}$,  which yields to a contradiction.

To prove \ref{item:P3b}, we recall that the solutions to $\Sigma_{\epsilon}$ cannot slide on $\partial K_{\bar{\epsilon}}$ along positive time intervals. Hence, the same property must hold for the solutions to $\Sigma^{-}_{\epsilon}$. As a result, if \ref{item:P3b} is not satisfied, then there exists a maximal solution $\phi \in \mathcal{S}_{\Sigma^{-}_{\epsilon}}(x_o)$, for some $x_o \in \partial K_{\bar{\epsilon}}$, and $T_o>0$ such that
$\phi((0,T_o]) \subset \mbox{int} (K_{\bar{\epsilon}})$.
Using Lemma \ref{lem1},  we conclude the existence of $T \in (0, T_o]$ and $\delta > 0$ such that,  for $y := \phi(T) \in K_{\bar{\epsilon}}$,  we have 
$$ x_o + \delta \mathbb{B}  \subset 
 R^b_{\Sigma^-_{\bar{\epsilon}}}(-T,y) = R^b_{\Sigma_{\bar{\epsilon}}}(T,y).  $$   
 However, the latter contradicts forward invariance of the set $K_{\bar{\epsilon}}$ for $\Sigma_{\bar{\epsilon}}$.

To prove \ref{item:P3c} using contradiction,   we assume the existence of a solution $\phi$  to $\Sigma_\epsilon$ starting from 
$x \in \text{int}(K_{\bar{\epsilon}})$  such that, for some $t_1>0$, we have 
 $ \phi(t_1) \in 
\partial K_{\bar{\epsilon}}$ and $\phi([0,t_1)) \in 
\text{int} (K_{\bar{\epsilon}}). $
Next, we take 
$ y := \phi(t_1) \in R^b_{\Sigma_{\epsilon}}(t_1,x) \backslash 
\{x\} $
and 
$$ z := \phi(t_1-t) \in R^b_{\Sigma_{\epsilon}}
(-t,y) \backslash \{y\} = R^b_{\Sigma^{-}_{\epsilon}}(t,y) \backslash \{y\}, $$
for some $t \in (0,t_1]$ to be determined. 
Hence, $z \in \text{int}(K_{\bar{\epsilon}})$. Next, using Lemma \ref{lem1} while 
replacing $(x,y,\Sigma,\Sigma_{\epsilon})$ by 
$(y,z,\Sigma^{-}_{\epsilon}, 
\Sigma^{-}_{(\epsilon + \bar{\epsilon})/2})$, we conclude the existence of $T>0$ such that, 
for $t := \min\{T,t_1\}$, there exists 
$\delta>0$ such that 
$  y + \delta \mathbb{B} \subset 
R^b_{\Sigma_{(\bar{\epsilon} + \epsilon)/2}}(t,z). $
However, according to \ref{item:P3a}, the 
set $\cl (K_{\bar{\epsilon}})$ must be forward invariant for $\Sigma_{(\bar{\epsilon} + \epsilon)/2}$, since 
\begin{align} \label{eqeps}
0 < (\bar{\epsilon}(x) + \epsilon(x))/2 < \bar{\epsilon}(x) \qquad \forall x \in \mathbb{R}^n, 
\end{align}
which yields to a contradiction.

Similarly, to prove \ref{item:P3d} using contradiction, we assume the existence of a solution $\phi$ to $\Sigma_\epsilon$ starting from $x \in \mathbb{R}^n \backslash \cl(K_{\bar{\epsilon}})$ 
such that, for some $t_1>0$, we have 
 $ \phi(-t_1) \in \partial K_{\bar{\epsilon}}$ and $\phi([-t_1,0)) \subset \text{int} (K_{\bar{\epsilon}})$.   
Next, we take 
$$ y := \phi(-t_1) \in R^b_{\Sigma_{\epsilon}}(-t_1,x) \backslash 
\{x\} $$
and 
$$ z := \phi(t-t_1) \in R^b_{\Sigma_{\epsilon}}
(t,y) \backslash \{y\} = R^b_{\Sigma^{-}_{\epsilon}}(-t,y) \backslash \{y\}, $$
for some $t \in (0,t_1]$ to be 
determined. Hence, $z \in \mathbb{R}^n \backslash \cl(K_{\bar{\epsilon}})$. Next, using Lemma \ref{lem1} while 
replacing $(x,y,\Sigma,\Sigma_{\epsilon})$ by 
$(y,z,\Sigma_{\epsilon}, 
\Sigma_{(\epsilon + \bar{\epsilon})/2})$, we conclude the existence of $T>0$ such that, 
for $t := \min\{T,t_1\}$, there exists 
$\delta>0$ such that 
$$  y + \delta \mathbb{B} \subset 
R^b_{\Sigma^{-}_{(\bar{\epsilon} + \epsilon)/2}}(t,z) = R^b_{\Sigma_{(\bar{\epsilon} + \epsilon)/2}}(-t,z). $$
However, according to \ref{item:P3b}, the 
set $\cl (K_{\bar{\epsilon}})$ must be forward invariant for $\Sigma_{(\bar{\epsilon} + \epsilon)/2}$,  since  \eqref{eqeps} holds, 
which yields to a contradiction.

\subsection{Proof of Lemma \ref{lem3}} 
 
To prove \ref{item:R1},  we will show that,  for each $x_o \in \partial K_{\bar{\epsilon}}$,  there exists $\delta > 0$ such that,  for each $x \in (x_o + \delta \mathbb{B}) \backslash \cl(K_{\bar{\epsilon}})$ and for each solution $\phi \in \mathcal{S}_{\Sigma_{\epsilon}}(x)$, there exists $T_\phi > 0$ such that $ \phi(T_\phi) \in \mbox{int} (K_{\bar{\epsilon}}) $.  First,  using Lemma \ref{lem2} and  according to \ref{item:P3a} and \ref{item:P3c},  
we conclude that $ R^b_{\Sigma_\epsilon}(T,x_o) \subset \mbox{int}(K_{\bar{\epsilon}})$ for all $T>0$.  Furthermore,  since $F$ is locally bounded,  using Lemmas \ref{lemprepre} and \ref{lem1-},  we conclude that,  for each $T>0$ small,  the set $R^b_{\Sigma_\epsilon}(T,U(x_o))$ is compact,  for $U(x_o)$ a sufficiently small neighborhood of $x_o$.  
Hence,  there exists $\alpha>0$ such that 
$$ \min \{|y|_{\partial K_{\bar{\epsilon}}} : y \in  R^b_{\Sigma_\epsilon}(T,x_o)\}  \geq \alpha.  $$   
Now,  using Lemma \ref{lem1-} in the Appendix,  we conclude that there exists $\delta>0$ such that, for each $x \in (x_o + \delta \mathbb{B}) \backslash \cl(K_{\bar{\epsilon}})$, we have 
$$ R^b_{\Sigma_\epsilon}(T, x) \subset R^b_{\Sigma_\epsilon}(T, x_o) + \alpha/2 \mathbb{B} \subset \mbox{int}(K_{\bar{\epsilon}}),  $$ 
which proves \ref{item:R1}.  Finally, \ref{item:R2} can be proved following the exact steps, while using \ref{item:P3b} and \ref{item:P3d} instead of \ref{item:P3a} and \ref{item:P3c}.

\section{Proof of Theorem \ref{thm4}} \label{Sec.6}
 
Given  $\delta, ~\bar{\epsilon} \in \mathcal{C}_+$,  we introduce the set $K_{\bar{\epsilon},\delta}$ given by
\begin{align} \label{Kinter1}
K_{\bar{\epsilon},\delta} := \bigcup_{x \in K_{\bar{\epsilon}}} (x + \delta(x) \mathbb{B}). 
\end{align}

Furthermore,  given $\rho_o, \epsilon_1 \in \mathcal{C}_+$,  we introduce the set $K_{\bar{\epsilon},\rho_o,\epsilon_1}$ given by
\begin{equation} 
\label{Kinter1+}
\begin{aligned} 
K_{\bar{\epsilon},\rho_o,\epsilon_1}   := \bigcup_{t \geq 0} \bigcup_{x \in K_{\bar{\epsilon},\rho_o}} R_{\Sigma_{\epsilon_1}} (t, x),
\quad 
K_{\bar{\epsilon},\rho_o}  := \bigcup_{x \in K_{\bar{\epsilon}}} (x + \rho_o(x) \mathbb{B}).
\end{aligned}
\end{equation}

The proof of Theorem \ref{thm4} follows in five steps.

\begin{itemize}
\item [\underline{Step 1}] 
\blue{ The next statement establishes recurrence of the set $K_{\bar{\epsilon},\delta}$.   This result is similar to \cite[Theorem 2]{SUBBARAMAN201654} which,  although  formulated for general hybrid systems,  studies  recurrence of bounded sets.    }

\begin{lemma} \label{lem6}
Consider system $\Sigma$ such that Assumption \ref{ass3} holds. Consider two subsets $(X_o,X_u) \in \mathbb{R}^n \times \mathbb{R}^n$ such that $\cl(X_o) \cap \cl(X_u) = \emptyset$ and  
$\bar{\epsilon}, \epsilon_1 \in \mathcal{C}_+$  such that 
\begin{align*} 
\epsilon_1(x) < \bar{\epsilon}(x) \qquad \forall x \in \mathbb{R}^n.
\end{align*}
Consider the set $K_{\bar{\epsilon}}$ introduced in \eqref{eqinfreachbis} and a closed set $U_1 \subset \mathbb{R}^n$ such that the conclusions of Lemmas \ref{lem2} and \ref{lem3} hold. 
Then, there exists
$\delta \in \mathcal{C}_+$  such that the set $K_{\bar{\epsilon},\delta}$ in \eqref{Kinter1} satisfies 
\begin{align}   \label{eqpropadd}
\cl(K_{\bar{\epsilon}, \delta}) \backslash \text{int} (K_{\bar{\epsilon}}) \subset \text{int} (U_1),  \qquad  \cl(K_{\bar{\epsilon}, \delta}) \cap \cl(X_u) = \emptyset
\end{align}
and,  for each 
$\epsilon \in \mathcal{C}_+$  satisfying 
\begin{align} \label{eqepsilon2n}
\epsilon(x) \leq \epsilon_1(x) \qquad \forall x \in \mathbb{R}^n, 
\end{align}
the following properties hold. 
\begin{enumerate} [label={P\ref{lem6}\arabic*)},leftmargin=*]
\item \label{item:I2} The set $K_{\bar{\epsilon},\delta}$ is locally recurrent for $\Sigma_{\epsilon}$ on $U_1$.
\item \label{item:I3} The set $\mathbb{R}^n \backslash K_{\bar{\epsilon},\delta}$ is locally recurrent for 
$\Sigma^{-}_{\epsilon}$ on $U_1$.
\end{enumerate}
\end{lemma}

\item[\underline{Step 2}] 
\blue{ Since we fail to establish contractivity of the set $K_{\delta,\bar{\epsilon}}$,  we consider the set  $K_{\bar{\epsilon},\rho_o,\epsilon_1}$ which we show to be both recurrent and contractive.   }

\begin{lemma} \label{lem6+}
Consider system $\Sigma$ such that Assumption  \ref{ass3} holds.   Consider two subsets $(X_o,X_u) \subset \mathbb{R}^n \times \mathbb{R}^n$ such that $\cl(X_o) \cap \cl(X_u) = \emptyset$ and  
$\bar{\epsilon}, \epsilon_1, \epsilon_2 \in \mathcal{C}_+$ such that 
\begin{align*} 
\epsilon_2 (x) < \epsilon_1(x) < \bar{\epsilon}(x) \qquad \quad 
\forall x \in \mathbb{R}^n.
\end{align*}
Consider the set $K_{\bar{\epsilon}}$ introduced in \eqref{eqinfreachbis},
 a closed subset $U_1 \subset \mathbb{R}^n$, and  $\delta \in \mathcal{C}_+$ such that the conclusions in Lemmas \ref{lem2}, \ref{lem3},  and \ref{lem6} hold.
 Then, there exists   
$\rho_o \in \mathcal{C}_+$ such that:

\begin{align} 
K_{\bar{\epsilon}} & \subset  \text{int}(K_{\bar{\epsilon}, \rho_o, \epsilon_1}) \label{eqpropadd1-}
\\
  \cl (K_{\bar{\epsilon}, \rho_o, \epsilon_1}) & \subset \text{int}( K_{\bar{\epsilon}, \delta}),  \label{eqpropadd1}
\\
\bigcup_{x \in \mathbb{R}^n \backslash \text{int}(K_{\bar{\epsilon},\delta})}  \hspace{-0.4cm} (x + \rho_o(x) \mathbb{B})  \cap \cl(K_{\bar{\epsilon},\rho_o,\epsilon_1}) & = \emptyset,  \label{eqpropadd2}
\end{align}

and the following properties hold.  

\begin{enumerate} [label={P\ref{lem6+}\arabic*)},leftmargin=*]
\item \label{item:I8+} There exists a closed subset $\hat{U}_1 \subset \text{int}(U_1)$ such that
\begin{align}     
\cl(K_{\bar{\epsilon}, \delta}) \backslash \text{int} (K_{\bar{\epsilon}}) & \subset \text{int} (\hat{U}_1), \label{eqlem102}  \\
 \bigcup_{x \in \hat{U}_1}  (x + \rho_o(x) \mathbb{B}) & \subset U_1, \label{eqlem101} 
\end{align}

\item  For each
$\epsilon \in \mathcal{C}_+$ satisfying
\begin{align*}
\epsilon (x) \leq  \epsilon_2(x) \qquad \forall x \in \mathbb{R}^n,
\end{align*}
\begin{enumerate}
\item \label{item:I2+a}   $K_{\bar{\epsilon},\rho_o,\epsilon_1}$ is locally recurrent for $\Sigma_{\epsilon}$ on $U_1$.
\item \label{item:I2+b} $\mathbb{R}^n \backslash K_{\bar{\epsilon},\rho_o,\epsilon_1}$ is locally recurrent for 
$\Sigma^{-}_{\epsilon}$ on $U_1$.
\item \label{item:I2+c} 
$\cl(K_{\bar{\epsilon},\rho_o,\epsilon_1})$ is forward contractive for $\Sigma_{\epsilon}$.

\item \label{item:I2+d}  $\mathbb{R}^n \backslash \mbox{int}(K_{\bar{\epsilon},\rho_o,\epsilon_1})$ is forward contractive for $\Sigma^-_{\epsilon}$.
\end{enumerate}
\end{enumerate}
\end{lemma}

\item[\underline{Step 3}]  
We assume that $\cl(X_o) \cap \cl(X_u) = \emptyset$
and we let $\bar{\epsilon}_o$,  
$\bar{\epsilon}$,  $\epsilon_1$,  and $\epsilon_2$ be robust-safety margins satisfying
\begin{align*} 
\epsilon_2 (x) < \epsilon_1(x) < \bar{\epsilon}(x)  < \bar{\epsilon}_o(x)
\qquad \quad 
\forall x \in \mathbb{R}^n.
\end{align*} 
Moreover,  we consider $\delta, \rho_o \in \mathcal{C}_+$, 
 the sets $(K_{\bar{\epsilon}}, K_{\bar{\epsilon},\delta}, K_{\bar{\epsilon},\rho_o,\epsilon_1})$  defined in \eqref{eqinfreachbis},  \eqref{Kinter1},  and  \eqref{Kinter1+}, respectively,  and closed subsets 
$U_1 \subset \mathbb{R}^n$ and $\hat{U}_1 \subset \mathbb{R}^n$  such that the conclusions of Lemmas \ref{lem3},  \ref{lem6},  and \ref{lem6+} hold.

We introduce the map 
$B_{\cl(K_{\bar{\epsilon},\rho_o,\epsilon_1})} : \mathbb{R}^n \rightarrow \mathbb{R} \cup 
\{\pm \infty \}$ defined in \eqref{eqbarcand} while replacing $(K, \Sigma)$ therein by $(\cl (K_{\bar{\epsilon},\rho_o,
\epsilon_1}),\Sigma_{\epsilon_2})$, and we show that 
\begin{enumerate} [label={S3\arabic*)},leftmargin=*]
\item \label{item:B1--+} 
$\cl(K_{\bar{\epsilon},\rho_o,\epsilon_1}) = \{ x \in \mathbb{R}^n : B_{\cl(K_{\bar{\epsilon},\rho_o,\epsilon_1})}(x) \leq 0 \}$ and  $B_{\cl(K_{\bar{\epsilon},\rho_o,\epsilon_1})}$ is a barrier candidate with respect to $(X_o,X_u)$. 
\item \label{item:B3--+}  $B_{\cl(K_{\bar{\epsilon},\rho_o,\epsilon_1})}$ is upper semicontinuous on $\mathbb{R}^n \backslash \text{int} (K_{\bar{\epsilon},\rho_o,\epsilon_1})$.

\item \label{item:B4--+}  $B_{\cl(K_{\bar{\epsilon},\rho_o,\epsilon_1})}$ is lower semicontinuous on $\cl(K_{\bar{\epsilon},\rho_o,\epsilon_1})$.

\item \label{item:B5--+} 
Along each solution $\phi$ (not necessarily maximal) to 
$\Sigma_{\epsilon_2}$ with
$\phi(\dom \phi) \subset U_1$,  we have 
$$ B_{\cl(K_{\bar{\epsilon},\rho_o,\epsilon_1})}
(\phi(t')) - B_{\cl(K_{\bar{\epsilon},\rho_o,\epsilon_1})}(\phi(t)) \leq - (t'-t) \qquad \forall (t,t') \in \dom \phi \times \dom \phi  \text{  s.t.  }  t' \geq t.  
$$
\end{enumerate}
We establish the proof by verifying the conditions in  Lemma \ref{lem5-} when $(K, \Sigma) = (\cl(K_{\bar{\epsilon}, \rho_o, \epsilon_1}), \Sigma_{\epsilon_2})$.  Indeed,  we start using \eqref{eqpropadd},  \eqref{eqpropadd1}-\eqref{eqpropadd2}, and the definition of the set $K_{\bar{\epsilon}, \rho_o, \epsilon_1}$,  to conclude that  $\cl(K_{\bar{\epsilon}, \rho_o, \epsilon_1}) \cap \cl(X_u) = \emptyset$,  $ \partial K_{\bar{\epsilon}, \rho_o, \epsilon_1} \subset \text{int} (U_1)$,  and  $X_o \subset \text{int}(K_{\bar{\epsilon}, \rho_o, \epsilon_1})$.  
Hence,  \ref{item:B1--+} is verified.   Next,  we show that \ref{item:71--}-\ref{item:74--} are verified for $(K, \Sigma)= (\cl(K_{\bar{\epsilon}, \rho_o, \epsilon_1}), \Sigma_{\epsilon_2})$.   Indeed, using \ref{item:I2+a}, we conclude that \ref{item:71--} holds,  using \ref{item:I2+b}, we conclude that \ref{item:72--} holds,  using \ref{item:I2+c}, we conclude that \ref{item:73--} holds,  and, finally, using \ref{item:I2+d}, we conclude that \ref{item:74--} holds.

\item[\underline{Step 4}] 

Given a smooth function $\rho_2 : \mathbb{R}^n \rightarrow \mathbb{R}_{>0}$  and $v \in \mathbb{B}$,  we introduce the function
$B^v_{\cl(K_{\bar{\epsilon},
\rho_o,\epsilon_1})} : \mathbb{R}^n \rightarrow \mathbb{R}$
given by 
\begin{align} \label{eq.hatB}
B^{v}_{\cl(K_{\bar{\epsilon},
\rho_o,\epsilon_1})}(x) :=  B_{\cl(K_{\bar{\epsilon},
\rho_o,\epsilon_1})}(x + \rho_2(x) v),
\end{align}
and we show that,  for each  $\rho_2 \in \mathcal{C}_+$ satisfying  
\begin{align} \label{eqrho}
 \rho_2(x) \leq \rho_1(x) := \min \{ \epsilon_2(x), \rho_o(x) \} \qquad \forall x \in \mathbb{R}^n, 
\end{align}
 the following properties hold.

\begin{enumerate} [label={S4\arabic*)},leftmargin=*]
\item \label{item:Bhat1}  $B^{v}_{\cl(K_{\bar{\epsilon},\rho_o,\epsilon_1})}$ is a barrier function candidate with respect to 
$(X_o,X_u)$.
\item \label{item:Bhat2}  $B^{v}_{\cl(K_{\bar{\epsilon},\rho_o,\epsilon_1})}$ is upper semicontinuous on the set $K_{vo} := \{x \in \mathbb{R}^n : B^{v}_{\cl(K_{\bar{\epsilon},\rho_o,\epsilon_1})}(x) \geq 0 \}$.

\item \label{item:Bhat3} $B^{v}_{\cl(K_{\bar{\epsilon},\rho_o,\epsilon_1})}$ is lower semicontinuous on the set 
$K_{vi} := \{x \in \mathbb{R}^n : B^{v}_{\cl(K_{\bar{\epsilon},\rho_o,\epsilon_1})}(x) \leq 0 \}$.
\end{enumerate}
Furthermore,  we show  the existence of a smooth function $\rho_2 : \mathbb{R}^n \rightarrow \mathbb{R}_{>0}$ satisfying \eqref{eqrho} such that, for each $v \in \mathbb{B}$, the function $B^{v}_{\cl(K_{\bar{\epsilon},\rho_o,\epsilon_1})}$ in \eqref{eq.hatB} satisfies the property:
\begin{enumerate} 
[label={S45)},leftmargin=*]
\item \label{item:Bhhat5} Along each solution $\phi$ (not necessarily maximal) to $\Sigma$ with $\phi(\dom \phi) \subset \hat{U}_1$,  we have  
$$ B^{v}_{\cl(K_{\bar{\epsilon},\rho_o,\epsilon_1})}(\phi(t')) - B^{v}_{\cl(K_{\bar{\epsilon},\rho_o,\epsilon_1})}(\phi(t)) \leq - (t'-t) \qquad \forall 
(t',t) \in \dom \phi \times \dom \phi \text{  s.t.   } t' \geq t.   $$
\end{enumerate}
Indeed,  we already showed  that $B_{\cl(K_{\bar{\epsilon}, \rho_o, \epsilon_1})}$ is well defined on $\mathbb{R}^n$. Therefore, $B^v_{\cl(K_{\bar{\epsilon}, \rho_o, \epsilon_1})}$  is also well defined on $\mathbb{R}^n$.  
Next,  by definition,  we know that 
$$ B_{\cl(K_{\bar{\epsilon}, \rho_o, \epsilon_1})}(y) \leq 0 \qquad  \forall y \in \cl(K_{\bar{\epsilon}, \rho_o, \epsilon_1}). $$
Then, given $v \in \mathbb{B}$, we conclude that
$$ B^{v}_{\cl(K_{\bar{\epsilon}, \rho_o, \epsilon_1})}(x) \leq 0 \qquad  \forall (x+\rho_2(x)v) \in \cl(K_{\bar{\epsilon}, \rho_o, \epsilon_1}). $$  

To show that 
\begin{align} \label{eqtoprove}
B^{v}_{\cl(K_{\bar{\epsilon}, \rho_o, \epsilon_1})}(x) \leq 0 \qquad \forall x \in X_o,
\end{align}
it is enough to use \eqref{eqrho} and the definition of the set $K_{\bar{\epsilon}, \rho_o, \epsilon_1}$ in  \eqref{Kinter1+},  to conclude that
$$\rho_2(x) \leq \rho_o(x) \qquad \forall x \in \mathbb{R}^n
\qquad \quad   \text{and} \quad  \qquad  x + \rho_o(x) v \in K_{\bar{\epsilon}, \rho_o, \epsilon_1} \qquad  \forall x \in K_{\bar{\epsilon}}.  $$
Next,  since $X_o \subset K_{\bar{\epsilon}}$, we conclude that
$$ x + \rho_2(x) v \in K_{\bar{\epsilon}, \rho_o, \epsilon_1} \qquad \forall x \in X_o,  $$ 
which proves \eqref{eqtoprove}. 

Next, to show that 
\begin{align} \label{eqproofstep}
B^{v}_{\cl(K_{\bar{\epsilon}, \rho_o, \epsilon_1})}(x)>0 \qquad \forall x \in X_u, 
\end{align}
we start noting that 
$$B_{\cl(K_{\bar{\epsilon}, \rho_o, \epsilon_1})}(y)>0 \qquad  \forall y \in \mathbb{R}^n \backslash \cl(K_{\bar{\epsilon}, \rho_o, \epsilon_1}). $$
Next, using Lemma \ref{lem6} and \eqref{eqpropadd},  we conclude that $ \cl(K_{\bar{\epsilon},\delta}) \cap \cl(X_u) = \emptyset. $
Hence,  
$$ x \in \mathbb{R}^n \backslash \cl(K_{\bar{\epsilon}, \delta}) \qquad  \forall x \in X_u.   $$ 
Moreover, using \eqref{eqpropadd1}-\eqref{eqpropadd2},  we conclude that
$$ x+\rho_o(x) v \notin \cl(K_{\bar{\epsilon}, \rho_o,\epsilon_1}) \qquad  \forall x \in \mathbb{R}^n \backslash \cl(K_{\bar{\epsilon}, \delta}).  $$   
Therefore,  \eqref{eqproofstep} is satisfied and $B^{v}_{K_{\bar{\epsilon}, \rho_o,\epsilon_1}}$ is a barrier function candidate with respect to 
$(X_o,X_u)$. 

Now,   since the function $\rho_2$ is  continuous,  the upper and lower semicontinuity properties of $B^v_{\cl(K_{\bar{\epsilon}, \rho_o,\epsilon_1})}$ are inherited from the those  of $B_{\cl(K_{\bar{\epsilon}, \rho_o,\epsilon_1})}$. 
Indeed, for each $x \in K_{vo}$, 
$y := x + \rho_2(x) v \in \mathbb{R}^n \backslash \text{int} (K_{\bar{\epsilon},\rho_o,\epsilon_1})$. 
According to item \ref{item:B3--+}, the function $B_{\cl(K_{\bar{\epsilon},\rho_o,\epsilon_1})}$ is upper semicontinuous at $y$. Hence, $B^v_{\cl(K_{\bar{\epsilon},\rho_o,\epsilon_1})}$ is upper semicontinuous on $K_{vo}$. Similarly, for each $x \in K_{vi}$,  $y := x + \rho_2(x) v \in  \cl(K_{\bar{\epsilon},\rho_o,\epsilon_1})$. 
Using item \ref{item:B4--+}, the function $B_{\cl(K_{\bar{\epsilon},\rho_o,\epsilon_1})}$ is lower semicontinuous at $y$.  Hence, $B^v_{\cl(K_{\bar{\epsilon},\rho_o,\epsilon_1})}$ is lower semicontinuous on $K_{vi}$.

Finally,   using Lemma \ref{lem8}, we conclude the existence of a smooth function $\rho_2 : \mathbb{R}^n \rightarrow \mathbb{R}_{>0}$ satisfying \eqref{eqrho}, such that \ref{item:star} holds. That is, for this choice of the function $\rho_2$, for each $v \in \mathbb{B}$, for each $\phi$ solution to $\Sigma$ starting from $x_o \in \hat{U}_1$ and remaining in $\hat{U}_1$, and given $(t,t') \in \dom \phi \times \dom \phi$ with $t' \geq t$, we conclude the existence of $\psi : [t,t'] \rightarrow \mathbb{R}^n$ solution to $\Sigma_{\rho_1}$ such that $ \psi(s) := \phi(s) + \rho_2(\phi(s)) v$ for all $s \in [t,t']$.  Thus, the function $\psi$ is solution to $\Sigma_{\epsilon_2}$.   Next, using the last item in Lemma \ref{lem6+},  we  conclude that $ \psi(s) := \phi(s) + \rho_2(\phi(s)) v \in U_1$ for all $s \in [t,t']$.   Finally, using \ref{item:B5--+}, we conclude that
$$ 
B^{v}_{\cl(K_{\bar{\epsilon},\rho_o,\epsilon_1})}(\phi(t')) - B^{v}_{\cl(K_{\bar{\epsilon},\rho_o,\epsilon_1})}(\phi(t)) 
=
B_{\cl(K_{\bar{\epsilon},\rho_o,\epsilon_1})}
(\psi(t')) - B_{\cl(K_{\bar{\epsilon},\rho_o,\epsilon_1})}(\psi(t)) \leq - (t'-t), $$
which implies that \ref{item:Bhhat5} is satisfied.

\item[\underline{Step 5}] 

We consider  the function $B : \mathbb{R}^n \rightarrow \mathbb{R}$ given by 
\begin{equation*}
\begin{aligned} 
B(x) & := 
\int_{\mathbb{R}^n}  B_{\cl(K_{\bar{\epsilon},
\rho_o,\epsilon_1})}(x + \rho_2(x) v) \Psi(v) dv,  
\end{aligned}
\end{equation*}
where $\rho_o : \mathbb{R}^n \rightarrow \mathbb{R}_{>0}$ and $\Psi : \mathbb{R}^n \rightarrow [0,1]$ are smooth functions such that
\begin{equation*}
\begin{aligned} 
\Psi(v) = 0  \quad \forall v \in \mathbb{R}^n \backslash \mathbb{B} \quad \text{and} \quad
\int_{\mathbb{R}^n}  \Psi(v) dv = 1, 
\end{aligned}
\end{equation*}
 is continuously differentiable.  We show that the following properties hold.
\begin{enumerate} [label={S5\arabic*)},leftmargin=*]
\item \label{item:B2-+}  $B$ is a barrier candidate with respect to $(X_o,X_u)$. 
\item \label{item:B3-+}  $B$ is continuously differentiable.
\item \label{item:B4-+} 
$\langle \nabla B(x), \eta \rangle 
\leq -1$ for all $\eta \in F(x)$ and for all $x \in \hat{U}_1$.

\item \label{item:B5-+} 
$\partial K \subset  \text{int}(\hat{U}_1)$, where $K := \{ x \in \mathbb{R}^n : B(x) \leq 0 \}$.
\end{enumerate}

To prove \ref{item:B2-+}, 
we use \ref{item:Bhat1} to conclude that
$$B^{v}_{\cl(K_{\bar{\epsilon},\rho_o,\epsilon_1})}(x) \leq 0 \qquad \forall v \in \mathbb{B},  \quad \forall x \in X_o. $$
Hence, since $\Psi$ is nonnegative, we conclude that
$$ \int_{\mathbb{R}^n} B^{v}_{\cl(K_{\bar{\epsilon},\rho_o,\epsilon_1})}(x) \Psi(v) dv \leq 0 \qquad  \forall x \in X_o.  $$ 
Furthermore,  we have that
$$ B^{v}_{\cl(K_{\bar{\epsilon},\rho_o,\epsilon_1})}(x) > 0 \qquad 
\forall v \in \mathbb{B},  \qquad \forall x \in X_u. $$   
Hence, using the fact that $\psi$ is nonnegative and there exist points in $\mathbb{B}$ where $\psi$ not zero, we conclude that 
$$ \int_{\mathbb{R}^n} B^{v}_{\cl(K_{\bar{\epsilon},\rho_o,\epsilon_1})}(x) \Psi(v) dv > 0 \qquad  \forall x \in X_u. $$ 

To prove \ref{item:B3-+}, we use Lemma \ref{lemsmoothing} under the fact that $B_{\cl(K_{\bar{\epsilon}, \rho_o, \epsilon_1})}$ is locally bounded and,  at every $x \in \mathbb{R}^n$,   either upper or lower semicontinuous.   

To prove \ref{item:B4-+}, we consider a solution $\phi$ to $\Sigma$ (not necessarily maximal) such that $\phi(\dom(\phi))\subset \hat{U}_1$.   For each $(t,t') \in \dom \phi \times \dom \phi$ with $t'\geq t$, we have 
\begin{equation*}
\begin{aligned}  
B(\phi(t')) & = \int_{\mathbb{B}} B^{v}_{\cl(K_{\bar{\epsilon},\rho_o,\epsilon_1})}(\phi(t')) \Psi(v) dv
\\ &
\leq \int_{\mathbb{B}} \left[ B^{v}_{\cl(K_{\bar{\epsilon},\rho_o,\epsilon_1})}(\phi(t)) - (t'-t) \right] \Psi(v) dv 
\\ &
= \int_{\mathbb{B}} B^{v}_{\cl(K_{\bar{\epsilon},\rho_o,\epsilon_1})}(\phi(t)) \Psi(v) dv - (t'-t) 
\\ & = B(\phi(t))-(t'-t),
\end{aligned}
\end{equation*}
where the first inequality is obtained  using \ref{item:Bhhat5}. Next, using Lemma \ref{lemA9bis+} in the Appendix,  under the fact that $B$ is continuously differentiable,   \ref{item:B4-+} follows.

To prove \ref{item:B5-+},  we propose to show that 
\begin{align} \label{eqlast} 
K_{\bar{\epsilon}} \subset K := \{ x \in \mathbb{R}^n : B(x) \leq 0 \} \subset K_{\bar{\epsilon},\delta}.  
\end{align}
Indeed,  we use \eqref{eqrho} and \eqref{Kinter1+} to conclude that
$$ x + \rho_2(x) \mathbb{B} \subset x+ \rho_o(x) \mathbb{B} \subset K_{\bar{\epsilon}, \rho_o, \epsilon_1} \qquad 
\forall x \in K_{\bar{\epsilon}}.  $$ 
Hence, for each $v \in \mathbb{B}$, we have 
$$ B^{v}_{\cl(K_{\bar{\epsilon}, \rho_o, \epsilon_1})}(x)  = B_{\cl(K_{\bar{\epsilon}, \rho_o, \epsilon_1})}(x + \rho_2(x) v) \leq 0 \qquad  \forall x \in K_{\bar{\epsilon}}. $$
The latter implies that 
\begin{align} \label{eqtouse} 
B(x) \leq 0 \qquad \forall x \in K_{\bar{\epsilon}}. 
\end{align}

Next,  using \eqref{eqpropadd1}-\eqref{eqpropadd2} and \eqref{eqrho},  we obtain
$$
x + \rho_2(x) \mathbb{B}  \subset
x + \rho_o(x) \mathbb{B} \subset \mathbb{R}^n \backslash 
\cl(K_{\bar{\epsilon}, \rho_o, \epsilon_1}) \qquad \forall x \in \mathbb{R}^n \backslash \text{int}(K_{\bar{\epsilon},\delta}).  
$$
Hence, for each 
$v \in \mathbb{B}$, 
$$ B^{v}_{\cl(K_{\bar{\epsilon}, \rho_o, \epsilon_1})}(x) = B_{\cl(K_{\bar{\epsilon}, \rho_o, \epsilon_1})}(x + \rho_2(x)v)  > 0 \qquad \forall x \in \mathbb{R}^n \backslash \text{int}(K_{\bar{\epsilon}, \delta}).  $$
The latter yields
\begin{align} \label{eqtouse1}
B(x) > 0 \qquad \forall x \in \mathbb{R}^n \backslash \text{int}(K_{\bar{\epsilon}, \delta}).
\end{align}
Finally,  using \eqref{eqtouse} and \eqref{eqtouse1},  we conclude that \eqref{eqlast} follows,  which implies \ref{item:B5-+} also follows using \eqref{eqpropadd}. 
\end{itemize} 

\subsection{ Proof of Lemma \ref{lem6}}
\blue{We start using 
Lemma \ref{lem2bis} to conclude that $\cl(X_u) \cap \cl(K_{\bar{\epsilon}}) = \emptyset$. 
We next use 
Lemma \ref{lem3} to conclude the existence of a closed subset $U_1 \subset \mathbb{R}^n$ such that \eqref{eqinter0} holds and \ref{item:R1}-\ref{item:R2} hold 
for each   
$\epsilon \in \mathcal{C}_+$ 
 satisfying \eqref{eqepsilon2n}.
As a result, to prove the existence of  $\delta \in \mathcal{C}_+$ for which \eqref{eqpropadd} holds,  we combine \eqref{eqinter0} and the fact that $\cl(X_u) \cap \cl(K_{\bar{\epsilon}}) = \emptyset$ to conclude that \eqref{eqpropadd}  is verified for
$$ \delta(x) := \frac{1}{2} \min \left\{ |x|_{\mathbb{R}^n \backslash U_1},  |x|_{\text{cl}(X_u)} \right\}.  $$

To prove \ref{item:I2}, we use \ref{item:R1}, which shows that the set $K_{\bar{\epsilon}}$ is locally recurrent for $\Sigma_{\epsilon}$ on $U_1$. Hence, \ref{item:I2} follows since, by definition, 
$K_{\bar{\epsilon}} \subset K_{\bar{\epsilon},\delta}$. 

To prove \ref{item:I3}, we grid the set 
$\partial K_{\bar{\epsilon}}$ using a sequence of nonempty compact subsets 
$\{\partial K_i\}^{N}_{i=1}$, where $N \in \{1,2,...,\infty\}$. That is, we assume that 
 $\bigcup^{\infty}_{i = 1} \partial K_i = \partial K_{\bar{\epsilon}}$.   Furthermore, for each $i \in \{1,2,...,N\}$, there exists a finite set  $\mathcal{N}_i \subset \{1,2,...,N \}$ and $\beta_i>0$ such that  
 \begin{align} \label{eqhelp}
 (\partial K_i + \beta_i \mathbb{B}) \cap \partial K_j = \emptyset \quad  
 \forall j \notin \mathcal{N}_i  \quad  \text{and} \quad  
\mbox{int}_{\partial K_{\bar{\epsilon}}} (\partial K_i \cap \partial K_j) = \emptyset \quad  
\forall j \in \mathcal{N}_i, 
\end{align}
where $\mbox{int}_{\partial K_{\bar{\epsilon}}}(\cdot)$ denotes the  interior of $(\cdot)$ relative $\partial K_{\bar{\epsilon}}$.   Such a decomposition always exists according to the Whitney Covering Lemma \cite{10.2307/1989708}.

Next, we make the following claim that we prove later. 
\begin{claim} \label{clm1}
For each $i \in 
\{1,2,...,N\}$,  we can find $t_i>0$ and $\alpha_i>0$ such that
\begin{align} \label{eqtoshow}
 \bigcup_{i \in 
 \{1,2,...,N\}} R^b_{\Sigma^{-}_\epsilon}(t_i, \partial K_i)        
\bigcap  
\bigcup_{i \in \{1,2,...,N\}} \bigcup_{x \in \partial K_i} (x + \alpha_i \mathbb{B})
= \emptyset. 
\end{align}
\end{claim}
Under Claim \ref{clm1}, we introduce the function $\alpha : \partial K_{\bar{\epsilon}} \rightarrow \mathbb{R}_{>0}$ given by 
$$ \alpha(x) := \min \{\alpha_i/2 : 
i \in \{1,2,...,N\} ~ \text{s.t.} ~ x \in \partial K_i \}. $$
The function $\alpha$ is lower semicontinuous and allows, in view of \eqref{eqtoshow}, to conclude that
$$ \bigcup_{i \in \{1,2,...,N\}} R^b_{\Sigma^{-}_\epsilon}(t_i, \partial K_i)        
\bigcap  
 \bigcup_{x \in \partial K_{\bar{\epsilon}}} (x + \alpha(x) \mathbb{B})
= \emptyset. $$
Hence,  using   \cite[Theorem 1]{katvetov1951real},  we conclude the existence of a continuous function $\bar{\alpha} : \partial K_{\bar{\epsilon}} \rightarrow \mathbb{R}_{>0}$ 
such that
$$  \bar{\alpha}(x) \leq \alpha(x) \qquad  
\forall x \in \partial K_{\bar{\epsilon}}. $$  Hence,  
\begin{align} \label{eqhelp2}  
\bigcup_{i \in \{1,2,...,N\}} R^b_{\Sigma^{-}_\epsilon}(t_i, \partial K_i)        
\bigcap  
 \bigcup_{x \in \partial K_{\bar{\epsilon}}} (x + \bar{\alpha}(x) \mathbb{B})
= \emptyset. 
\end{align}
That is, every solution to $\Sigma^{-}_{\epsilon}$ starting from $\partial K_{\bar{\epsilon}}$ leaves the set
 $K_{\bar{\epsilon}} \cup \bigcup_{x \in \partial K_{\bar{\epsilon}}} (x + \bar{\alpha}(x) \mathbb{B})$.

Next, we show that, for each $i \in 
\{1,2,...,N\}$, we can find $\sigma_i > 0$ such that 
 $$  \bigcup_{i \in \{1,2,...,N\}} R^b_{\Sigma^{-}_\epsilon}(t_i, \partial K_i + \sigma_i \mathbb{B})        \bigcap  
 \bigcup_{x \in \partial K_{\bar{\epsilon}}} (x + \bar{\alpha}(x) \mathbb{B})
= \emptyset. $$
Indeed, the latter follows from a direct combination of Lemmas \ref{lem1-} and \ref{lemprepre}
 in the Appendix and \eqref{eqhelp2}. 

As a result, the function $\sigma : \partial K_{\bar{\epsilon}} \rightarrow \mathbb{R}_{>0}$ given by 
$$ \sigma(x) := \min \{\sigma_i : 
i \in \{1,2,...,N\} ~ \text{s.t.} ~ x \in \partial K_i \} $$
 is lower semicontinuous and allows, using   \cite[Theorem 1]{katvetov1951real},  to conclude the existence of a continuous function $\bar{\sigma} : \partial K_{\bar{\epsilon}} \rightarrow \mathbb{R}_{>0}$ 
such that
$$  \bar{\sigma}(x) \leq \sigma(x) \qquad  
\forall x \in \partial K_{\bar{\epsilon}}, $$ 
and, at the same time, the solutions  to $\Sigma^-_{\epsilon}$ starting from the set
$ \bigcup_{x \in \partial K_{\bar{\epsilon}}} (x + \bar{\sigma}(x) \mathbb{B}) \subset $
leave the set 
$ K_{\bar{\epsilon}} \cup \bigcup_{x \in \partial K_{\bar{\epsilon}}} (x + \bar{\alpha}(x) \mathbb{B}). $
Finally, after continuously extending 
$\bar{\sigma}$ to $\mathbb{R}^n$, we take
$$ \delta(x) := 
\left\{
\begin{matrix}
\min \{ \bar{\sigma}(x), |x|_{\partial (K_{\bar{\epsilon}} \cup K_1)} \} & \text{if}~ x \in \text{int} (K_{\bar{\epsilon}})
\\
\bar{\sigma}(x) & \text{otherwise},
\end{matrix}
\right. \qquad K_1 := \bigcup_{x \in \partial K_{\bar{\epsilon}}} (x + \bar{\sigma}(x) \mathbb{B}). $$
Using \ref{item:R2}, we conclude that the solutions to $\Sigma^{-}_{\epsilon}$ starting from $ K_{\bar{\epsilon}, \delta} \cap U_1$ 
leave the set $K_{\bar{\epsilon}, \delta}$. Hence, \ref{item:I3} is verified. 
 
To prove Claim \ref{clm1}, we use Lemma \ref{lem2} to conclude that, for each $i \in \{1,2,...,N\}$ and for each $t_i > 0$, we have 
$$ R^b_{\Sigma^{-}_\epsilon}(t_i, \partial K_i) \cap \cl(K_{\bar{\epsilon}}) = \emptyset. $$
Furthermore, using Lemma \ref{lemprepre} in the Appendix, we conclude that, for each $i \in \{1,2,...,N\}$, we can make  $t_i > 0$ small enough such that
$R^b_{\Sigma^{-}_\epsilon}(t_i, \partial K_i)$ is compact. Furthermore, in view of \eqref{eqhelp} and Lemma \ref{lem1-}, we can make each $t_i$ even smaller such that
\begin{align} \label{eqRconst}
(\partial K_j + \beta_j \mathbb{B}) \cap  R_{\Sigma^{-}_\epsilon}(t_i, \partial K_i) = \emptyset \qquad \forall j \notin \mathcal{N}_i.
\end{align}
As a result, for each $i \in \{1,2,...,N\}$, we can find 
$\alpha_i \in (0,\beta_i)$ such that
\begin{align} \label{eqRconstbis} 
R^b_{\Sigma^{-}_\epsilon}(t_i, \partial K_i)  \bigcap  \bigcup_{j \in 
\{ i,\mathcal{N}_i \}} \left( \partial K_j + \alpha_j \mathbb{B} \right) = \emptyset. 
\end{align}
The claim is proved by combining \eqref{eqRconst} and \eqref{eqRconstbis}. }

\subsection{ Proof of Lemma \ref{lem6+}}

The proof of \eqref{eqpropadd1-} is obvious for any  $\rho_o \in \mathcal{C}_+$. 

 \blue{
Furthermore,  to find  $\rho_o : \mathbb{R}^n \rightarrow \mathbb{R}_{>0}$ such that     
 \eqref{eqpropadd1}-\eqref{eqpropadd2} hold,  we propose to grid the set $\text{cl}(K_{\bar{\epsilon}})$  using a sequence of nonempty compact subsets 
$\{K_i\}^{N}_{i=1}$, where $N \in \{1,2,...,\infty\}$.  
That is,  we assume that 
$$ \bigcup^{N}_{i = 1} K_i = \text{cl}(K_{\bar{\epsilon}}) $$
and,  for each $i \in \{1,2,...,N\}$,  there exists a finite set 
$\mathcal{N}_i \subset \{1,2,...,N \}$ 
such that  
$$ K_i \cap K_j = \emptyset  \quad  \forall j \notin \mathcal{N}_i  \qquad \text{and}  \qquad  
\mbox{int}_{\text{cl}(K_{\bar{\epsilon}})} (K_i \cap K_j) = \emptyset  \quad  \forall j \in \mathcal{N}_i,  $$
where $\mbox{int}_{K}(\cdot)$ is the interior of $(\cdot)$ relative to $K$. 

A key step consists in proving the following claim. 

\begin{claim} \label{clm2}
For each $i \in \{1,2,... , N\}$,   there exists $\rho_i > 0$ such that 
$$ \bigcup_{t \geq 0} \bigcup_{i \in \mathbb{N}} R_{\Sigma_{\epsilon_1}} (t, K_i + \rho_i \mathbb{B})  \subset 
\text{int}(K_{\bar{\epsilon},\delta}).
   $$ 
\end{claim}
Next, we introduce the function $\rho : \text{cl}(K_{\bar{\epsilon}})  \rightarrow 
\mathbb{R}_{>0}$ given by 
$$ \rho(x) := \min \{\rho_i : 
i \in \{1,2,...,N\} ~ \text{s.t.} ~ x \in  K_i \} $$
which is lower semicontinuous and allows us to conclude that 
$$ \bigcup_{t \geq 0} \bigcup_{x \in \text{cl}(K_{\bar{\epsilon}})} R_{\Sigma_{\epsilon_1}} \left(t,  x + \rho(x) \mathbb{B} \right)  \subset \text{int} (K_{\bar{\epsilon},\delta}).  $$
Furthermore,   we use  \cite[Theorem 1]{katvetov1951real} to conclude the existence of 
$\rho_1 : K_{\bar{\epsilon}} \rightarrow \mathbb{R}_{>0}$ continuous and satisfying  $\rho_1(x) \leq \rho(x)$  for all $x \in K_{\bar{\epsilon}}$.    
Hence, 
\begin{align*} 
 K_{\bar{\epsilon},  \rho_1,  \epsilon_1}   = \bigcup_{t \geq 0} \bigcup_{x \in \text{cl}(K_{\bar{\epsilon}})}  R_{\Sigma_{\epsilon_1}} (t, x + \rho_1(x) \mathbb{B}) \subset \text{int}(K_{\bar{\epsilon},\delta}). 
\end{align*}

On the other hand,  using the continuous function $\rho_2 : \mathbb{R}^n \backslash 
\text{int}(K_{\bar{\epsilon}, \delta})  \rightarrow \mathbb{R}_{>0}$ given by  $\rho_2(x) := \frac{1}{2} |x|_{K_{\bar{\epsilon},\rho_1,\epsilon_1}}$,  we  conclude that  $$ x + \rho_2(x) \mathbb{B} \subset \mathbb{R}^n \backslash 
\cl(K_{\bar{\epsilon}, \rho_1, \epsilon_1}) \qquad   \forall  
x \in \mathbb{R}^n \backslash \text{int}(K_{\bar{\epsilon},\delta}).  $$    

Finally,  after extending   $\rho_2$ to $\mathbb{R}^n$ and taking  $ \rho_o (x) := \min \{ \rho_1(x), \rho_2(x) \}$ for all $x \in \mathbb{R}^n$,  \eqref{eqpropadd1}-\eqref{eqpropadd2}  follow.

To prove Claim \ref{clm2},  we start using the fact that the set 
$\text{cl}(K_{\bar{\epsilon}})$ is forward contractive for $\Sigma_{\epsilon_1}$ and the solutions to $\Sigma_{\epsilon_1}$ starting from $\text{int}(K_{\bar{\epsilon}})$ never reach $\partial K_{\bar{\epsilon}}$ for positive times,  to conclude that 
$$ R_{\Sigma_{\epsilon_1}} (t, K_i) \backslash  K_i \subset \text{int} (K_{\bar{\epsilon}}) \qquad \forall i \in \{1,2,...,N\}, \qquad \forall t > 0.  $$

Next,  using Lemma \ref{lemprepre} in the Appendix,   we conclude that,  for each $i \in \{1,2,...,N\}$,  there exist
 $\bar{\rho}_i>0$  and $T_i >0$ such that $R_{\Sigma_{\epsilon_1}} (T_i ,  K_i + \bar{\rho}_i \mathbb{B})$ is bounded.    
 Furthermore,  using Lemma \ref{lem1-} in the Appendix,  we conclude that,  for each $i \in \{1,2,...,N\}$,  there exists $\eta_i > 0$ and $\rho_i \in (0, \bar{\rho}_i)$ such that 
$$|y|_{\partial K_{\bar{\epsilon},\delta}}  \geq \eta_i 
\qquad \forall y \in  R_{\Sigma_{\epsilon_1}} (T_i, K_i),  
\qquad 
|y|_{\partial K_{\bar{\epsilon}}} \geq \eta_i \qquad 
\forall  y \in R^b_{\Sigma_{\epsilon_1}} (T_i, K_i),  $$ 
and at the same time 
\begin{align*}
|y|_{\partial K_{\bar{\epsilon},\delta}} 
 \leq \eta_i/2 \quad \forall y \in  R_{\Sigma_{\epsilon_1}} 
(T_i, K_i + \rho_i \mathbb{B})  \qquad \text{and} \qquad 
|y|_{\partial K_{\bar{\epsilon}_1}}  \leq \eta_i/2 \quad \forall y \in R^b_{\Sigma_{\epsilon_1}} (T_i, K_i + \rho_i \mathbb{B}).
\end{align*} 

In particular, we conclude that 
$$ R^b_{\Sigma_{\epsilon_1}} (T_i, K_i + \rho_i \mathbb{B}) \subset \text{int}(K_{\bar{\epsilon}}) \qquad \forall i \in \{1,2,...,N \},  $$
which implies that
\begin{align*}
|y|_{\partial K_{\bar{\epsilon},\delta}} 
 \leq \eta_i/2 \qquad \forall y \in  R_{\Sigma_{\epsilon_1}} 
(t, K_i + \rho_i \mathbb{B})   \qquad \forall t \geq 0. 
\end{align*} 
The latter is enough to prove the claim. }

To prove \ref{item:I8+},  we use \eqref{eqpropadd} to
 conclude that  $\cl(K_{\bar{\epsilon}, \delta}) \backslash \text{int} (K_{\bar{\epsilon}}) \subset \text{int} (U_1)$.   
 Hence, we can always find a closed set $\hat{U}_1 \subset \text{int}(U_1)$ such
$ \cl(K_{\bar{\epsilon}, \delta}) \backslash \text{int} (K_{\bar{\epsilon}}) \subset \text{int} (\hat{U}_1), $
which implies that \eqref{eqlem102} holds. 
Next,   we let the function $\rho_3 : \hat{U}_1 \rightarrow \mathbb{R}_{>0}$ given by 
$ \rho_3 (x) := |x|_{\mathbb{R}^n \backslash U_1}$.  After extending $\rho_3$ to $\mathbb{R}^n$ and taking $\rho_o(x) := \min \{\rho_i(x) : i \in \{1,2,3 \} \}$ for all $x \in \mathbb{R}^n$,  both \eqref{eqpropadd1}-\eqref{eqpropadd2} and 
\ref{item:I8+} follow. 

To prove \ref{item:I2+a}), we use the fact $K_{\bar{\epsilon}} \subset K_{\bar{\epsilon}, \rho_o, \epsilon_1}$ and \ref{item:R1}.  Similarly, to prove \ref{item:I2+b}), we combine \eqref{eqpropadd1}-\eqref{eqpropadd2} and \ref{item:I3}.  To prove \ref{item:I2+c}) and  \ref{item:I2+d}),  we use Lemma \ref{lem2} while replacing the sets 
$(X_o, K_{\bar{\epsilon}})$ therein by the sets  $(K_{\bar{\epsilon},\rho_o}, K_{\bar{\epsilon},\rho_o,\epsilon_1})$.

\section{Conclusion}
In this paper, we establish the equivalence between robust safety and the existence of a smooth barrier certificate,  in the context of continuous-time systems modeled by differential inclusions.   Our result requires only continuity of the set-valued dynamics and empty intersection between the closures of the initial and unsafe sets. We relax most of the assumptions used in existing literature such as boundedness of the safety region, smoothness of the system's dynamics,  and uniqueness of solutions. 
 Future works pertain to address the considered problem in the more general context of hybrid dynamical systems.

\appendix

\section*{Appendix}

\blue{ In the following lemma,  we recall three different consequences of Assumption \ref{ass1} on the regularity of the reachability set-valued maps $R_{\Sigma}$ and $R^b_{\Sigma}$.  A similar result can be found in   \cite[Theorem 1, Page 103]{aubin2012differential} under a slightly different set of assumptions.   }

\begin{lemma} \label{lem1-}
Consider system $\Sigma$ such that  Assumption \ref{ass1} holds and let $T$ and  $U \subset   \mathbb{R}^n$ such that  $R_\Sigma(T, U)$ is bounded.  Then,   the maps $R_{\Sigma}$ and $R^b_{\Sigma}$ are both upper and outer semicontinuous on $[0,T] \times U$.
\end{lemma}

\begin{proof}
In view of Remark \ref{remplus},   when the maps $R_{\Sigma}$ and $R^b_{\Sigma}$ are locally bounded,  then
showing outer semicontinuity is enough to conclude 
upper semicontinuity.  Thus, we propose to only show outer semicontinuity of $R^b_{\Sigma}$. 
Indeed,  let $(t,x) \in [0,T] \times U$ and let two sequences $\left\{ (t_{i},x_{i}) \right\}^{\infty}_{i=0} \subset [0,T] \times U$ and $\left\{ y_i \right\}^{\infty}_{i=0} \subset \mathbb{R}^n$ such that $\lim_{i \rightarrow \infty} (t_{i}, x_{i}) = (t,x)$, $y_i \in R^b_{\Sigma}(t_{i},x_{i})$, and $\lim_{i \rightarrow \infty} y_i = y \in \mathbb{R}^n$.   
Outer semicontinuity of $R^b_{\Sigma}$ at $(t,x)$ follows if we show that $y \in R^b_{\Sigma}(t,x)$.   To this end,  we consider a sequence of  solutions $\left\{\phi_i \right\}^{\infty}_{i=0}$  to $\Sigma$ such that 
$$\dom \phi_i = [0,t_i] \quad  \text{and}  \quad y_i = \phi_i(t_{i}) \qquad  \forall i \in \{0,1,...\}.  $$  
\blue{Now,  since the sequence $\left\{\phi_i \right\}^{\infty}_{i=0}$ is uniformly bounded,   by passing to an adequate subsequence,  we conclude the existence of a continuous function $\phi : [0,T] \rightarrow \mathbb{R}^n$ constituting the graphical limit of  the sequence $\left\{\phi_i \right\}^{\infty}_{i=0}$ such that 
$$ \lim_{i \rightarrow \infty} \phi_i(t) = \phi(t) \quad  
\forall t \in [0,T],   \qquad \quad T = \lim_{i \rightarrow \infty} t_i.  $$ Hence,  $\phi(0) = x$ and $y = \lim_{i \rightarrow \infty} \phi_i(t_{i}) = \phi(t)$.  Finally,  using \cite[Theorem 5.29]{goebel2012hybrid},  we conclude that $\phi$ is solution to $\Sigma$ and thus $y \in R^b_{\Sigma}(t,x)$.  }

Now,  to show outer semicontinuity of $R$, we consider two sequences $\left\{(t_{i},x_{i})\right\}^{\infty}_{i=0} \subset [0,T] \times U$ and 
$\left\{ y_i \right\}^{\infty}_{i=0} \subset \mathbb{R}^n$ such that 
$\lim_{i \rightarrow \infty} (t_{i}, x_{i}) = (t,x)$, 
$y_i \in R(t_{i},x_{i})$,  and  $\lim_{i \rightarrow \infty} y_i = y \in \mathbb{R}^n$. 
Outer semicontinuity of $R$ at $(t,x)$ follows if we show that $y \in R(t,x)$.  
Having $y_i \in R(t_{i},x_{i})$, for each $i \in \mathbb{N}$, implies the existence of $t'_i \in [0,t_{i}]$ such that $y_i \in R^b(t'_{i},x_{i})$, for each $i \in \mathbb{N}$. By passing to an adequate subsequence, we conclude the existence of $t' \in [0,t]$ such that $t' = \lim_{i \rightarrow \infty} t'_i$. Hence,  since $R^b$ is outer semicontinuous, we conclude that 
$y \in R^b(t',x) \subset R(t,x)$. 
\end{proof}

\blue{We next present a  useful consequence of having $F$ locally bounded.   \ifitsdraft \else A proof can be found in \cite{maghenem2022converse}.  \fi }

\begin{lemma}  \label{lemprepre}
Consider the differential inclusion $\Sigma$ and assume that the map $F$ therein is locally bounded.  
Then,  the following properties are true:

\blue{ 
\begin{itemize}
\item  For each compact set $K \subset \mathbb{R}^n$,   there exist $b>0$ and $\bar{T}>0$ such that the sets
$R_{\Sigma} (\bar{T},  (K + b \mathbb{B}))$
and $R_{\Sigma^-} (\bar{T},  (K + b \mathbb{B}))$
are bounded.   
\item Given a sequence of solutions $\left\{\phi_i \right\}^{\infty}_{i=0}$ to $\Sigma$ starting from a compact set $K$ such that,  for each $i \in \{1,2,...\}$,  $\dom \phi_i = [0,t_i]$, $t_i > 0$,  and $\lim_{i \rightarrow \infty} t_i = t_o \in (0,+\infty]$.   Then,  there exists $\bar{T} \in (0,  t_o]$ the largest time such that the sequence $\{ \phi_i \}^{\infty}_{i =1}$ is uniformly bounded on any $[0, T] \subset [0,\bar{T})$.    
\end{itemize}
}
\end{lemma}

 \ifitsdraft 

\textcolor{brown}{
\begin{proof}
Since $F$ is locally bounded,  we conclude the existence of $L$ and $M >0$ such that $F( K + L \mathbb{B}) \subset M \mathbb{B}$.  Now,  by taking $b := L/2$ and $\bar{T} := L/(2M)$,  it follows that 
$ F( K + (b + \bar{T} M) \mathbb{B}) \subset M \mathbb{B}.   $
This implies that the solutions to $\Sigma$ starting from $K + b \mathbb{B}$, over the window of time $[0,\bar{T}]$ or $[-\bar{T},0]$,  cannot leave the set $ K + (b + \bar{T} M) \mathbb{B}$.  As a result, we obtain 
$$ R_{\Sigma} (\bar{T},  (K + b \mathbb{B}) ) \subset  K + (b + \bar{T} M) \mathbb{B} \quad \text{and} \quad 
 R_{\Sigma^-} (\bar{T},  (K + b \mathbb{B}) ) \subset  K + (b + \bar{T} M) \mathbb{B}.  $$
Hence,  $R_{\Sigma} (\bar{T},  (K + b \mathbb{B}))$
and
$R_{\Sigma^-} (\bar{T},  (K + b \mathbb{B}))$
are bounded.      
The proof of the second item follows using the first item. 
\end{proof}}

\fi

\blue{The following Lemma follows from the combination of    \cite[Theorem 5.2.1 and Lemma 5.1.2]{Aubin:1991:VT:120830}. 
\ifitsdraft
\else
 A detailed proof can be found in \cite{maghenem2022converse}.  
 \fi }

\begin{lemma} \label{lemA14}
Consider system $\Sigma$ in \eqref{eq.1} such that Assumption \ref{ass1} holds.  A closed set $K \subset \mathbb{R}^n$ is forward invariant for $\Sigma$ if 
\begin{align} \label{eqtengApp} 
F(x) \subset C_K(y) \qquad \forall  x \in \mathbb{R}^n \backslash K, \quad \forall y \in \text{Proj}_{K}(x).
\end{align}
\end{lemma}

\ifitsdraft

\begin{proof}
\textcolor{brown}{We start using \cite[Theorem 5.2.1]{Aubin:1991:VT:120830} to conclude that the set
$K$ is forward invariant if
\begin{equation}  
\label{eq.econe}
\begin{aligned} 
F(x) \subset E_K(x) & \qquad \forall x \in \mathbb{R}^n \backslash K,
\end{aligned}
\end{equation}
where 
$$E_K(x) := \left\{ v\in \mathbb{R}^n :\liminf_{h \rightarrow 0^+} \frac{|x+hv|_K - |x|_K}{h} \leq 0  \right\}.  $$
Next,  to complete the proof,  we use \cite[Lemma 5.1.2]{Aubin:1991:VT:120830} to conclude that
$$  C_K(y) \subset E_K(x) \qquad \forall y \in \text{Proj}_K(x), \quad \forall x \in \mathbb{R}^n \backslash K.   $$
Hence,   \eqref{eq.econe} is verified under  \eqref{eqtengApp}.}
\end{proof} 

\fi

\textcolor{blue}{ The following lemma can be deduced from \cite[Theorem 6.3]{clarke2008nonsmooth}.  
Although formulated in the nonsmooth setting and for locally-Lipschitz dynamics,  the same proof applies to our case.  \ifitsdraft
\else A detailed proof is in  \cite{maghenem2022converse}.  \fi }

\begin{lemma} \label{lemA9bis+}
Consider system $\Sigma$ in \eqref{eq.1} such that Assumption  \ref{ass3} holds.  Consider an open set  $O \subset \mathbb{R}^n$  and a continuously-differentiable function $B : \mathbb{R}^n \rightarrow \mathbb{R}$ such that, 
along each solution $\phi$ satisfying $\phi(\dom \phi) \subset O$,  the map $t \mapsto B(\phi(t))$ is nonincreasing.  Then, 
$$ \langle \nabla B(x) , \eta  \rangle \leq 0 \qquad  \forall \eta \in F(x),  \quad \forall x \in O.  $$ 
\end{lemma}

\ifitsdraft

\textcolor{brown}{
\begin{proof} 
Let $x_o \in O$ and $v_o \in F(x_o)$.  Since $F$ is continuous and has closed and convex images,  using Michael's selection theorem \cite{michael1956continuous},  we conclude the existence of a continuous selection $v : U(x_o) \rightarrow \mathbb{R}^n $ such that $v(x) \in F(x)$ for all $x \in U(x_o)$ with $v(x_o) = v_o$.  
Next,  using \cite[Proposition 3.4.2]{Aubin:1991:VT:120830},  we conclude the existence of a nontrivial continuously differentiable solution $\phi$ starting from $x_o$ solution to the system $\dot{x} = v(x)$; thus,  $\phi$ is also solution to $\Sigma$.  Furthermore, we consider a sequence  $\left\{ t_i \right\}^{\infty}_{i=0} \subset \dom \phi$ such that $\lim_{i \rightarrow \infty} t_i = 0$.  Note that
$
\frac{d}{dt} (B(\phi(t)))|_{t=0}  = \langle \nabla B(x_o),v(x_o) \rangle 
 = \lim_{t_i \rightarrow 0} \frac{B(\phi(t_i)) - B(\phi(0)) }{t_i}  \leq 0. 
$
\end{proof}}

\fi

\blue{ The remaining lemmas are deduced from  \cite{goebel2012hybrid},  where they are  formulated for the general context of hybrid inclusions.   
\ifitsdraft
\else
Simpler proofs than those in \cite{goebel2012hybrid} can be found in \cite{maghenem2022converse}.  \fi }

\begin{lemma} [Lemma 7.37.  \cite{goebel2012hybrid}] \label{lem8}
Consider system $\Sigma$ in \eqref{eq.1} such that Assumption  \ref{ass3} holds.  Then, for each  $\rho_1 \in \mathcal{C}_+$,  there exists a smooth function $\rho_2 \in \mathcal{C}_+$ satisfying  
\begin{align} \label{eqrhorho0}
\rho_2(x) \leq \rho_1(x) \qquad \forall x \in \mathbb{R}^n
\end{align}
such that the following property holds. 
\begin{enumerate} [label={P\ref{lem8}\arabic*)},leftmargin=*]
\item \label{item:star} 
For each $\phi \in \mathcal{S}_{\Sigma}(x_o)$,  for each $v \in \mathbb{B}$,  and for each $t \in \dom \phi \cap \mathbb{R}_{\geq 0}$,  the function  $\psi: [0,t] \rightarrow  \mathbb{R}^n$ given by $\psi(s) := \phi(s) + \rho_2(\phi(s)) v$  is solution to $\Sigma_{\rho_1}$. 
\end{enumerate}
\end{lemma}
 
\ifitsdraft

\textcolor{brown}{
\begin{proof}
Let $\phi \in \mathcal{S}_{\Sigma}(x_o)$,  $t \in \dom \phi \cap \mathbb{R}_{\geq 0}$,    and  $v \in \mathbb{B}$.  Given $\rho_2 : \mathbb{R}^n \rightarrow \mathbb{R}_{>0}$ continuous,  the function 
$\psi$ satisfies  
\begin{align*}
\dot{\psi}(s)  = \dot{\phi}(s) + \nabla \rho_2(\phi(s)) \dot{\phi}(s) v  \subset  F(\phi(s)) + |\nabla \rho_2(\phi(s))| f(\phi(s)) \mathbb{B} \quad  \text{for a. a.} ~ s \in [0,t],
\end{align*}
where $f(x) :=  \sup \{ |\zeta| : \zeta \in F(x) \}$,  which is  continuous in our case.   Now,   given 
$\rho : \mathbb{R}^n \rightarrow \mathbb{R}_{>0}$ continuous,  we show the existence of $\rho_2 : \mathbb{R}^n \rightarrow \mathbb{R}_{>0}$ continuous
such that 
\begin{align} \label{eqTCK}
 \rho_2(x) \leq \rho(x) ~ \text{and} ~  |\nabla \rho_2(x)|  \leq  \frac{\rho(x)}{  f(x) +1}  \quad \forall x \in \mathbb{R}^n.  
\end{align}
Under \eqref{eqTCK}, we conclude that 
\begin{align} \label{eqIntermuse}
\dot{\psi}(s)   \subset  F(\phi(s)) + \rho(\phi(s)) \mathbb{B} \quad \text{for a. a.} ~ s \in [0,t].  
\end{align}
Now,  to verify \eqref{eqTCK},  we partition $\mathbb{R}^n$ using a locally finite cover $\{K_i\}^\infty_{i=1}$ with $\cl(K_i)$ compact, and subordinate to this
cover a smooth partition of unity $\{ \psi_i \}^{\infty}_{i=1}$.   Finally,   we let 
$$ \rho_2(x) :=  \sum^\infty_{i=1} \frac{2^{1-i} a_i}{\max_{z \in K_i} \max \{ \psi_i(z),  |\nabla \psi_i(z)|  \} }  \psi_i(x), $$
where $a_i \in (0,1)$ such that $a_i \leq \rho(x)$ for all $x \in K_i$ and
$ a_i \sup_{z \in K_i }  f(z)  \leq 1$.
\\
To complete the proof,  we use Lemma \ref{lemA333} in the Appendix to conclude that given
  $\rho_1 : \mathbb{R}^n \rightarrow \mathbb{R}_{>0}$ continuous,  we can find $\rho \in \mathcal{C}_+$  such that 
\begin{equation}
\label{eqinterPr}
\begin{aligned}
\rho(x) & \leq \rho_1(x) \qquad \forall x \in \mathbb{R}^n,
\\
F(x) + \rho(x) \mathbb{B} &  \subset F(y) + \rho_1(y) \mathbb{B}  \qquad \forall y \in x + \rho(x) \mathbb{B}. 
\end{aligned}
\end{equation}
As a result,  since 
$$ \psi(s) \in \phi(s) + \rho_2(\phi(s)) \mathbb{B} \subset \phi(s) + \rho(\phi(s)) \mathbb{B} \quad  \forall s \in [0,t], 
$$    
 applying \eqref{eqinterPr},  we conclude that 
$$ F(\phi(s)) + \rho(\phi(s)) \mathbb{B}   \subset F(\psi(s)) + \rho_1(\psi(s)) \mathbb{B} \quad \forall s \in [0,t],  $$
and using \eqref{eqIntermuse},  we obtain 
$$ \dot{\psi}(s)  \subset  F(\psi(s)) + \rho_1(\psi(s)) \mathbb{B} \quad \text{for almost all} ~ s \in [0,t]. $$
The latter implies that $\psi$ is solution to $\Sigma_{\rho_1}$. 
\end{proof}}
  
  \fi

\begin{lemma} [Section 7.6.0.3.  in \cite{goebel2012hybrid}] \label{lemsmoothing}
Let $B_o : \mathbb{R}^n \rightarrow \mathbb{R}$ be  locally bounded.  Assume that,  at every $x \in \mathbb{R}^n$,  $B_o$ is either upper or lower semicontinuous.   Then,  the function 
$B : \mathbb{R}^n \rightarrow \mathbb{R}$ 
given by
\begin{align} \label{eqBinteg} 
B(x) := \int_{\mathbb{R}^n} B_{o}(x + \rho_o(x) v) \Psi(v) dv,  
\end{align}
where $\rho_o : \mathbb{R}^n \rightarrow \mathbb{R}_{>0}$ and $\Psi : \mathbb{R}^n \rightarrow [0,1]$ are smooth functions such that
\begin{equation}
\label{eqPsi} 
\begin{aligned} 
\Psi(v) = 0  \quad \forall v \in \mathbb{R}^n \backslash \mathbb{B} \quad \text{and} \quad
\int_{\mathbb{R}^n}  \Psi(v) dv = 1, 
\end{aligned}
\end{equation}
 is continuously differentiable. 
\end{lemma}

\ifitsdraft

\textcolor{brown}{
\begin{proof}
 Consider the change of coordinate 
$ w := x + \rho_o(x) v. $
Furthermore, we introduce the closed set  
$ \mathbb{B}_{\rho_o}(x) := x + \rho_o(x) \mathbb{B}. $ 
Hence, we obtain 
$$B(x)= \int_{\mathbb{R}^n} B_{o}(w) \Psi\left(\frac{w-x}{\rho_o(x)} \right) \frac{dw}{\rho_o(x)}.$$
Next, we let 
$$ g(x):= \int_{\mathbb{R}^n} f(x,w) dw  :=  \int_{\mathbb{R}^n} B_{o}(w) \Psi\left(\frac{w-x}{\rho_o(x)} \right) dw. $$
\\
We will show that $g$ is continuously differentiable,  which would imply that $B$ is continuously differentiable,  since $\rho_o$ is smooth. 
\begin{itemize}
\item[a.] Note that 
$ w \mapsto f(x,w)$ 
is $\mathcal{L}^1$ since 
$ w \mapsto \Psi\left(\frac{w-x}{\rho_o(x)} \right) $
is null outside the bounded set  $\mathbb{B}_{\rho_o}(x)$. 
\item[b.] Note that the map $ x \mapsto \Psi \left( \frac{w-x}{\rho_o(x)} \right) $ is smooth and null outside $\mathbb{B}_{\rho_o}(x)$.  Hence,  $ w \mapsto  \nabla_x f(x,w) \in \mathcal{L}^1$. 
 \item[c.] Since   the map $ x \mapsto  \nabla_x \left( \Psi\left(\frac{w-x}{\rho_o(x)} \right) \right)$ is smooth and null outside the bounded set $\mathbb{B}_{\rho_o}(x)$, we conclude the existence of a positive constant $k$ such that
 $$ \sup_{y \in \mathbb{B}_{\rho_o}(x)} \displaystyle \left\lvert \nabla_x \left( \Psi\left(\frac{w-y}{\rho_o(y)} \right) \right) \displaystyle \right\rvert \leq k. $$
 Hence, 
 $ |\nabla_x f(x,w)| \leq  b(w), $
  where $b : \mathbb{R}^n \rightarrow \mathbb{R}_{\geq 0}$ is the upper semicontinuous function given by
 $$ b(w) := 
 \left\{ 
 \begin{matrix}
 k |B_{o}(w)| & \text{if} ~ w \in \mathbb{B}_{\rho_o}(x) 
 \\ 
 0 & \text{otherwise}.
 \end{matrix}
 \right. $$
Being null outside the  set  $\mathbb{B}_{\rho_o}(x)$,    $w \mapsto b(w)$ is $\mathcal{L}_1$. 
\end{itemize} 
As a result,   using  \cite[Lemma 7.38]{goebel2012hybrid},  we conclude that $g$ is differentiable and 
$$ \nabla_x g(x) = \int_{\mathbb{R}^n} B_{o}(w) \nabla_x \left(\Psi\left(\frac{w-x}{\rho_o(x)} \right) \right)^\top dw,  $$
which is continuous.  
\end{proof} }

\begin{lemma}  \label{lemA333}
Let $F: \mathbb{R}^n \rightrightarrows \mathbb{R}^n$ such that Assumption \ref{ass3} holds.   Then,  for each   $\rho_1 \in \mathcal{C}_+$,  there exists  $\rho \in \mathcal{C}_+$ such that
\begin{align} \label{eqcontic}
F(x) + \rho(x) \mathbb{B}  \subset F(y) +  \rho_1 (y)\mathbb{B} \quad 
\forall y \in x + \rho(x) \mathbb{B}. 
\end{align}
\end{lemma}

\textcolor{brown}{
\begin{proof}
We grid $\mathbb{R}^n$ using a sequence of nonempty compact subsets $\{I_i\}_{i=1}^{N} \subset \mathbb{R}^n$,  where $N \in \{1,2,3, \dots ,\infty \}$,  such that,  for each $i \in \{1,2,...,N\}$,  there exists 
$\mathcal{N}_i \subset \{1,2,...,N \}$  finite such that 
$ I_i \cap I_j = \emptyset$  for all  $j \notin \mathcal{N}_i$,   and $ \mbox{int}(I_i \cap I_j) = \emptyset$ for all $j \in \mathcal{N}_i\backslash \{i\}$.
Since $F$ is continuous; thus,  uniformly continuous  on each $I_i + \mathbb{B}$,  
$i \in \{1,2,...,N\}$ \footnote{  
The map $F: K \rightrightarrows \mathbb{R}^n$ is \textit{uniformly continuous} if, for each $\epsilon > 0$, there exists 
$\delta > 0$ such that, for each $x \in K$, we have  
$$ 
|F(x_1) - F(x_2) |_H \leq \epsilon \qquad  \forall x_1, x_2 \in x + \delta \mathbb{B}.  
 $$  
When $K$ is compact,  using the same arguments as in the  single-valued case,  continuity and uniform continuity become equivalent. },  
  we conclude that,  for each $\varepsilon_i > 0$,  there exists 
$\delta_i \in (0,1]$ such that 
$ F(x) + \delta_i  \mathbb{B}  \subset F(y)+ \varepsilon_i \mathbb{B}$  for all $x, y \in I_i + \mathbb{B}$  such that $ |y - x| \leq  \delta_i$.    In particular,  since $\delta_i \leq 1$,  we conclude that
$ F(x) + \delta_i  \mathbb{B}  \subset F(y)+ \varepsilon_i \mathbb{B}$ for all $x \in I_i$,   for all $y \in x +  \delta_i \mathbb{B}$.     
Now,  if we let $\varepsilon_i := \min_{y \in I_i} \rho_1(y)$,  we conclude the existence of $\delta_i > 0$ such that 
$ F(x) + \delta_i \mathbb{B}  \subset F(y)+ \rho_1(y) \mathbb{B}$ for all $x \in I_i$,  for all $y \in x +  \delta_i \mathbb{B}$.   Next, we introduce the function $\delta :\mathbb{R}^n \to \mathbb{R}_{>0}$ given by
$\delta(x) :=\text{min}\{\delta_i: i\in \{1,2,\dots , N\} \;\ \text{such that} \;\ x\in I_i\} $.  By definition,  $\delta$ is lower semi-continuous; hence,  $-\delta$ is upper semicontinuous.  
Using   \cite[Theorem 1]{katvetov1951real},  we conclude the existence of $\rho \in \mathcal{C}_+$  and satisfying 
$- \rho(x) \geq -\delta(x)$ for all $x \in \mathbb{R}^n$.
Hence, \eqref{eqcontic} follows. 
\end{proof}}

\fi

\balance

\bibliographystyle{unsrt}      
\bibliography{biblio.bib}

\end{document} \documentclass{cocv}
\usepackage{amsfonts,amssymb,amsmath}
\usepackage{graphicx,graphics,epsfig,color}
\usepackage{multicol}
\usepackage[colorinlistoftodos, textwidth=35mm, shadow]{todonotes}
\usepackage{rgsMacros}
\usepackage{refcheck}

\def\startmodifnew{\color{red}}
\def\stopmodifnew{\color{black}\normalcolor}

\let\labelindent\relax
\usepackage{enumitem}
\usepackage{balance}
\usepackage{enumitem}

\newcommand\red[1]{{\color{red}#1}}
    \definecolor{gray}{rgb}{0.33,0.4,0.47}\def\gray#1{{\color{gray}#1}}
    \definecolor{steelblue}{rgb}{0,.42,.7}\def\steelblue#1{{\color{steelblue}#1}}
    \definecolor{britishgreen}{rgb}{0,0.26,0.15}\def\britishgreen#1{{\color{britishgreen}#1}}
    \definecolor{navyblue}{rgb}{0,0,.8}\def\navyblue#1{{\color{navyblue}#1}}
    \definecolor{olivegreen}{rgb}{0.14,0.29,0}\def\olivegreen#1{{\color{olivegreen}#1}}
    \definecolor{myred}{rgb}{0.86,0.1,0.16}\def\myred#1{{\color{myred}#1}}

\newif\ifitsdraft
\def\itsdraft{\global\itsdrafttrue}

\newtheorem{exe}{Example}
\newtheorem{corol}{Corollary}
\newtheorem{ass}{Assumption}
\newtheorem{defin}{Definition}
\newtheorem{cla}{Claim}
\newtheorem{rem}{Remark}
\newtheorem{lem}{Lemma}
\newtheorem{prop}{Proposition}
\newtheorem{thm}{Theorem}
\newtheorem{fct}{Fact}
\newtheorem{prob}{Problem}
\newenvironment{lemma}{\begin{lem}}{\hfill $\square$ \end{lem}}
\newenvironment{proposition}{\begin{prop}}{\hfill $\square$ \end{prop}}
\newenvironment{corollary}{\begin{corol}}{\hfill $\square$ \end{corol}}
\newenvironment{example}{\begin{exe}\rm }{\hfill $\square$ \end{exe}}
\newenvironment{remark}{\begin{rem}\rm }{\hfill $\bullet$ \end{rem}}
\newenvironment{assumption}{\begin{ass}}{\hfill $\bullet$ \end{ass}}
\newenvironment{theorem}{\begin{thm}}{\hfill $\square$ \end{thm}}
\newenvironment{definition}{\begin{defin}}{\hfill $\bullet$ \end{defin}}
\newenvironment{claim}{\begin{cla}}{\hfill $\bullet$ \end{cla}}
\newenvironment{fact}{\begin{fct}}{\hfill $\bullet$ \end{fct}}
\newenvironment{problem}{\begin{prob}}{\hfill $\bullet$ \end{prob}}

\begin{document}

\title{\LARGE \bf A converse robust-safety theorem for differential inclusions}

\author{Mohamed Maghenem} \address{CNRS, Gipsa-lab,  Grenoble INP,  Universit{\'e} Grenoble Alpes,  Grenoble, France. Email: mohamed.maghenem@gipsa-lab.fr.}

\author{Masoumeh Ghanbarpour} \address{Department of Electrical and Computer Engineering, University of Colorado Bulder,  Colorado,  USA.  
\\
 Email:  masoumeh.ghanbarpour@gmail.com }

\subjclass{93A10-26B05}

\begin{abstract}
This paper establishes the equivalence between robust safety and the existence of a barrier function certificate for differential inclusions.   More precisely,  for a robustly-safe differential inclusion, a barrier function is constructed as the 
time-to-impact function with respect to a specifically-constructed reachable set.  Using techniques from 
set-valued and nonsmooth analysis, we show that such a function,  although being possibly discontinuous,  certifies robust safety by verifying a condition involving the system's solutions.  Furthermore, we refine this construction,  using smoothing techniques from the literature of converse Lyapunov theory, to provide a smooth barrier certificate that certifies robust safety by verifying a condition involving only the barrier function and the system's dynamics. In comparison with existing converse robust-safety theorems, our results are more general as they allow the safety region to be unbounded, the dynamics to be a general continuous set-valued map,  and the solutions to be non-unique.  
\end{abstract}   
            
\keywords{ Robust safety;
Differential inclusions; 
Barrier functions; 
Converse theorem. }

\maketitle

\section{Introduction}  

Safety for a dynamical system requires the solutions  starting from a given set of initial conditions to never reach a given unsafe set  \cite{prajna2007framework}.   Depending on the application,  reaching the unsafe set may correspond to 
non-applicability of a predefined feedback law,  due to saturation or a change in the dynamics or,  simply,  due to collisions with physical obstacles. Ensuring safety is in fact key in many engineering applications including traffic regulation \cite{ersal2020connected},  aerospace \cite{9656550}, and human-robot interactions \cite{9788028}. 

\subsection{Motivation}

This notion of safety is not robust in nature, as it is possible to construct safe differential equations that become unsafe when arbitrarily small perturbations are added to their right-hand side \cite[Example 1]{9683684}.   As a result, we say, roughly speaking, that a dynamical system is robustly safe if it remains safe in the presence of a perturbation term added to its dynamics.  This 
robust-safety notion was first introduced in \cite{wisniewski2016converse} for systems defined on compact manifolds.  A similar notion is studied in \cite{ratschan2018converse,  9444774} for continuous-time systems modeled by differential equations.  The same notion is considered in \cite{9683684} in the context of differential inclusions,  which generalize differential equations by allowing the right-hand side to be a general set-valued map \cite{aubin2012differential},  and thus the solutions to be non-unique.  

As the analytical expression of the solutions of a dynamical system are usually impossible to obtain, and since their precise approximation can be computationally expensive,    barrier functions are widely used to study safety and robust safety without
computing or approximating the  solutions. This is analogous to Lyapunov theory for stability.  We recall that a \textit{barrier function candidate} is a scalar function with opposite signs on the initial and the unsafe subsets. Furthermore,  it certifies safety, or robust safety, by satisfying an inequality constraint involving the barrier function candidate itself and the system's dynamics.  In which case,  the barrier function candidate becomes a \textit{barrier certificate} \cite{9705088, 9683684}.   
 Such conditions are well documented in the literature of safety under different smoothness properties of the barrier function candidate  \cite{ames2014control,
konda2019characterizing, 10.1007/978-3-642-39799-8_17} and the system's dynamics \cite{draftautomatica}. 
Furthermore, in the context of robust safety,  when an upper bound on the perturbation is known,   safety conditions involving  the (worst-case) perturbed dynamics are used in \cite{liu2020converse, JANKOVIC2018359}.  In \cite{seiler2021control},  specific classes of perturbations, solution to  some dynamical models and verifying a certain integral constraint are considered.  Perturbation-free conditions ensuring robust safety are proposed in  \cite{9683684, RubSafPI},  provided that mild regularity assumptions on the dynamics and the barrier function candidate hold.  Showing the necessity of such perturbation-free conditions is the main subject of the current paper. 

 \subsection{Background}

Converse safety and robust-safety problems pertain to show the existence of a barrier certificate for safety and robust safety, provided that the system is safe and robustly safe, respectively.

\begin{itemize}
\item In the context of safety, it is shown in \cite{9705088} that, in general, the existence of a continuous barrier certificate is not necessary for safety,   unless special cases are considered \cite{prajna2005necessity}.   Alternatively, time-varying barrier certificates are introduced and their existence is shown in \cite{9705088} to be necessary as well as sufficient,  under some assumptions on the system.    
\item In the context of robust safety,   \cite{wisniewski2016converse} solved the converse robust-safety problem by constructing a smooth barrier certificate for systems defined on smooth and compact manifolds,  provided that \blue{the system's dynamics are represented by a smooth single-valued map},  the initial and unsafe sets are compact and disjoint,  and a \textit{Meyer} function exists.     

\item In \cite{ratschan2018converse},  when the  system's dynamics are represented by a smooth 
single-valued map,  the complement of the unsafe set is bounded, and the closures of the initial and unsafe sets are disjoint,  robust safety is shown to be equivalent to the existence of a smooth barrier certificate. In particular,  to prove the converse robust-safety theorem  in \cite{ratschan2018converse}, the reachable set, denoted by $K_{\bar{\epsilon}}$, along the solutions to a perturbed version of the system starting from the initial set is introduced, where the subscript  $\bar{\epsilon}$ stands for the perturbation term added to the 
original-system's dynamics. 
After that, \blue{a barrier function candidate is defined, at any point in the state space, as the first to impact the boundary of the set $K_{\bar{\epsilon}}$ by a solution starting from that point}.  Such a function is shown to be a valid barrier function candidate. Although being only continuous, \blue{it is shown to be strictly decreasing along the solutions to the original system, when $\bar{\epsilon}$ is a robustness margin}.    
After that, a smooth barrier certificate is deduced using boundedness of the safe set \blue{and the density property of the class of smooth functions} in the space of continuous functions.   

\item A very similar converse robust-safety theorem is established in \cite{9444774} using converse Lyapunov theorems for asymptotic stability.   Indeed,  when either the reachable set $K_{\bar{\epsilon}}$ is bounded or the system's dynamics are represented by a globally Lipschitz map, and a uniform separation exists between the unsafe region and the set $K_{\bar{\epsilon}}$,  the latter set is shown to be uniformly asymptotically stable for the original system.  As a result,  existing converse Lyapunov theorems are used to show that a smooth barrier certificate exists.  
\end{itemize}
 
 \subsection{Contribution} 
 
In this paper,  we prove two  converse robust-safety theorems under mild regularity assumptions on the system's dynamics. \blue{Indeed, the latter is allowed to be represented by a set-valued map. Moreover,  we do not restrict the safety region to be bounded}. As in \cite{ratschan2018converse}, we show that blue{a specifically defined} 
time-to-impact function with respect to  the boundary of a reachable set $K_{\bar{\epsilon}}$,  when $\bar{\epsilon}$ is a robustness margin, is strictly decreasing when evaluated along the solutions to the original system lying on a neighborhood of $K_{\bar{\epsilon}}$.  However,  \blue{since the solutions are not necessarily unique},  this function is not necessarily continuous. Nonetheless, for an appropriate choice of the robustness margin $\bar{\epsilon}$,  the constructed function is shown to be a (non-smooth) barrier certificate. That is, it satisfies a sufficient condition for robust safety that is \textit{non-infinitesimal}; namely,  a condition involving the system's solutions.  
To construct a smooth barrier certificate, inspired by \cite{SUBBARAMAN201654}, we propose to smoothen the constructed nonsmooth one.  However, for the resulting smooth function to be a barrier function candidate; namely,  to have opposite signs on the initial and unsafe sets,  we need to carefully choose the set with respect to which the time-to-impact function is defined, as well as the different parameters involved in its construction.  As a consequence,  we show the existence of smooth barrier certificate provided that the system's dynamics are represented by a set-valued map that is continuous and the closures of the initial and unsafe sets are disjoint.   Finally,  we show the utility of our converse result in the context of safety for self-triggered control systems. 

 The rest of the paper is organized as follows.  Preliminaries on set-valued maps,  differential inclusions,  and  invariance and attractivity notions are in Section \ref{Sec.2}.  The problem formulation,  the main results,  and a motivational example  are in Section \ref{Sec.3}.   Preparatory materials towards the proofs of the main results are in Section \ref{Sec.4}.  
The proofs  of the main results are in Sections \ref{Sec.5} and \ref{Sec.6}. \blue{Finally, intermediate technical results are reported in the Appendix.} 

\blue{A preliminary version of this work is in \cite{9682926}, where the solutions are assumed to be forward complete. This  is not the case here. Furthermore,  proofs, detailed explanations, and the application example are not present in latter reference. }

\textbf{Notation.} 
For $x$,  $y \in \mathbb{R}^n$,  $x^{\top}$ denotes the transpose of $x$,  $|x|$ the Euclidean norm of $x$ and $\langle x, y \rangle:= x^\top y $ the inner product between $x$ and $y$.  For a set $K \subset \mathbb{R}^n$, 
we use $\mbox{int}(K)$ to denote its interior,  $\partial K$ to denote its boundary,  $U(K)$ to denote any open neighborhood of the set $K$,  and $|x|_K$ to denote the distance between $x$ and the set $K$.   Furthermore, we use $C_K(x)$ to denote the \textit{contingent} cone of $K$ at $x$, which is given by 
$$ C_K(x) := \left\{ v \in \mathbb{R}^n: \liminf_{h \rightarrow 0^+} |x + h v|_K/h = 0 \right\}.   $$ 
 For $O \subset \mathbb{R}^n$,  $K \backslash O$ denotes the subset of elements of $K$ that are not in $O$.    \textcolor{blue}{For a function $\phi : \dom \phi \rightarrow \mathbb{R}^m$,  $\dom \phi \subset \mathbb{R}^n$ denotes the domain of definition of $\phi$}.   By $F : \mathbb{R}^m \rightrightarrows \mathbb{R}^n $,  we denote a set-valued map associating each element $x \in \mathbb{R}^m$ \blue{with} a subset $F(x) \subset \mathbb{R}^n$.  In particular,  $\text{Proj}_{K} : \mathbb{R}^n \rightrightarrows K$ represents the projection set-valued map on $K$; namely, 
$ \text{Proj}_{K}(x) := \{y \in K : |x-y| = |x|_{K} \}$.  For a set $D \subset \mathbb{R}^m$,  $F(D) := \{ \eta \in F(x) : x \in D \}$.    For a differentiable map $B : \mathbb{R}^n \rightarrow \mathbb{R}$,  \textcolor{blue}{$\nabla_{x_i} B$ denotes  the derivative of $B$ with respect to $x_i$},  $i \in \{1,2,...,n\}$,  and $\nabla B$ denotes the gradient of $B$ with respect to $x$.   \ifitsdraft We say that $B \in \mathcal{L}^1$ if $|B|_1 := \int_{\mathbb{R}^n} |B(x)| dx$ is finite.  \fi
Finally,  $\mathbb{B}$ denotes the closed unit ball centered at the origin. 

\section{Preliminaries} \label{Sec.2}

\subsection{Set-valued vs single-valued maps} 

Consider a set-valued map $F: K \rightrightarrows \mathbb{R}^n$, where $K \subset \mathbb{R}^m$. 

\begin{itemize}
\item $F$ is \textit{outer semicontinuous} at $x \in K$ if,  for every sequence $\left\{x_i\right\}^{\infty}_{i=0} \subset K$ and for every sequence  $\left\{ y_i \right\}^{\infty}_{i=0} \subset \mathbb{R}^n$ with $\lim_{i \rightarrow \infty} x_i = x$, $\lim_{i \rightarrow \infty} y_i = y \in \mathbb{R}^n$, and $y_i \in F(x_i)$ for all $i \in \mathbb{N}$, we have $y \in F(x)$;  see \cite{rockafellar2009variational}. 

\item  $F$ is \textit{upper semicontinuous} at $x \in K$ if,  for each $\varepsilon > 0$,  there exists a neighborhood of $x$, 
denoted by $U(x)$,  such that for each $y \in U(x) \cap K$, $F(y) \subset F(x) + \varepsilon \mathbb{B}$; see \cite[Definition 1.4.1]{aubin2009set}.

\item  $F$ is  \textit{continuous} at $x \in K$  if,  for each $\epsilon > 0$,  there exists $\delta > 0$ such that 
\begin{align} \label{eqContin}
|F(x_1) - F(x_2) |_H \leq \epsilon \qquad  \forall x_1, x_2 \in x + \delta \mathbb{B},  
\end{align}  
where $|F(x) - F(y) |_H$ stands for the Hausdorff distance between the sets $F(x)$ and $F(y)$.   
 
\item  $F$ is \textit{locally bounded} at $x \in K$ if there exists a neighborhood of $x$, denoted by $U(x)$,  and 
$\beta > 0$ such that  $|\zeta| \leq \beta$ for all $\zeta \in F(y)$ and for all $y \in U(x) \cap K$.  
\end{itemize}

Furthermore,  the map $F$ is upper,  outer semicontinuous,  continuous, or locally bounded if, respectively,  so it is for all $x \in K$.

Consider a single-valued map $B: K \rightarrow \mathbb{R}$,  where $K \subset \mathbb{R}^m$. 

\begin{itemize}
\item  $B$ is \textit{lower semicontinuous} at $x \in K$ if, for every sequence $\left\{ x_i \right\}_{i=0}^{\infty} \subset K$ such that $\lim_{i \rightarrow \infty} x_i = x$, we have $\liminf_{i \rightarrow \infty} B(x_i) \geq B(x)$. 
\item  $B$ is \textit{upper semicontinuous} at $x \in K$ if, for every sequence $\left\{ x_i \right\}_{i=0}^{\infty} \subset K$ such that $\lim_{i \rightarrow \infty} x_i = x$, we have $\limsup_{i \rightarrow \infty} B(x_i) \leq B(x)$.  
\item  $B$ is \textit{continuous} at $x \in K$ if it is both upper and lower semicontinuous at $x$. 
\end{itemize}

Furthermore,  $B$ is upper, lower semicontinuous, or continuous if, respectively,  so it is for all $x \in K$. 

\subsection{Differential inclusions} 
\textcolor{blue}{We recall the notion of a 
Carath{\'e}odory solution to a differential inclusion of the form}
\begin{align} \label{eq.1}
\Sigma : \quad  \dot x \in F(x) \qquad x \in \mathbb{R}^n.
\end{align}
\begin{definition}
A function $\phi : \dom \phi  \rightarrow \mathbb{R}^n$,  \blue{with $\dom \phi \subset \mathbb{R}$ an interval containing $\{0\}$},  is a solution to $\Sigma$ if it is locally absolutely continuous and  $\dot{\phi}(t) \in F(\phi(t))$ for almost all $t \in \dom \phi$.
\end{definition}

A solution $\phi$ to $\Sigma$ is said to start from $x$ if $\phi(0) = x$.   A solution $\phi$ to $\Sigma$ is maximal if there is no solution $\psi$ to  $\Sigma$  such that $\psi(t) = \phi(t)$ for all $t \in \dom \phi$ and $\dom \phi$ \blue{is} strictly included in $\dom \psi$.  Furthermore,  given $T \geq 0$,  we use $\mathcal{S}^T_{\Sigma}(x)$ to denote the set of solutions $\phi$ to $\Sigma$ starting from $x$ with $\dom \phi = [0,T]$.  Finally,   we use $\mathcal{S}_{\Sigma}(x)$ to denote the set of maximal solutions $\phi$ to $\Sigma$ starting from $x$. 

\textcolor{blue}{ We now propose to view the sets of points reached by the solutions to $\Sigma$, starting from a given initial condition and over 
a given window of time,  as set-valued maps.  Indeed,  as in \cite[Section 4.2.]{refId0} and \cite[Page 104]{aubin2012differential},  we,  respectively,    recall the set-valued maps $R_{\Sigma} : \mathbb{R} \times \mathbb{R}^n \rightrightarrows \mathbb{R}^n$ and $R^b_{\Sigma} : \mathbb{R} \times \mathbb{R}^n \rightrightarrows \mathbb{R}^n \cup \emptyset $ given by 
\begin{align*}
 R_{\Sigma}(t,x)  := \{ \phi(s): \phi \in \mathcal{S}_{\Sigma}(x), ~ 
s \in \dom \phi \cap I_t \},  \quad
 R^b_{\Sigma}(t, x) & := \left\{ \phi(t) : \phi \in \mathcal{S}_{\Sigma}(x),  ~t \in \dom \phi  \right\},
\end{align*}
where $I_t := [\min\{0,t\}, \max\{0,t\}]$.  }
 In simple words,  the set $R_{\Sigma}(t,x)$ includes all the elements reached by the solutions to $\Sigma$ starting from $x$ over the interval $I_t$.  Furthermore,  the set $R^b_{\Sigma}(t,x)$ includes the value of the solutions starting from $x$ at $t$,  when $t$ is part of their domain.  
 
Finally,  we introduce the following assumption on $F$. 
\begin{assumption} \label{ass1} 
The map $F$ is upper semicontinuous and $F(x)$ is nonempty, compact,  and convex for all $x \in \mathbb{R}^n$.
\end{assumption} 

\textcolor{blue}{Assumption \ref{ass1} guarantees  the existence of a non-trivial solution from any $x \in \mathbb{R}^n$ as well as useful structural properties for the set of solutions to $\Sigma$;  see \cite{aubin2012differential,  refId0,  filippov2013differential}.  }
 
\begin{remark} \label{remplus}
\blue{Given a set-valued map $F : K \rightrightarrows \mathbb{R}^n$, where $K \subset \mathbb{R}^m$,  
we recall,  based on \cite[Theorem 5.19]{rockafellar2009variational} and \cite[Lemma 5.15]{goebel2012hybrid},  that  the following two properties are equivalent.
\begin{itemize}
\item  $F$ is upper semicontinuous and  $F(x)$ is compact for all $x \in K$. 
\item $F$ is outer semicontinuous and locally bounded.
\end{itemize}}
\end{remark}

\section{Problem formulation and results}  \label{Sec.3}

Given a set of initial conditions $X_o \subset \mathbb{R}^n$ and an unsafe set $X_u \subset \mathbb{R}^n$  such that $X_o \cap X_u = \emptyset$, 
we recall that $\Sigma$ is safe with respect to $(X_o,X_u)$ if,  for each solution $\phi$ with $\phi(0) \in X_o$,  we have 
$\phi(\dom \phi  \cap \mathbb{R}_{\geq 0} ) \subset \mathbb{R}^n \backslash X_u$. 
Note that safety with respect to  $(X_o, X_u)$ is verified if and only if there exists a set $K \subset \mathbb{R}^n$, with $X_o \subset K$ and $K \cap X_u = \emptyset$,  that is forward invariant.   
In turn, a set $K \subset \mathbb{R}^n$ is forward invariant if,  for each solution $\phi$ to $\Sigma$ with $ \phi(0) \in K$, $\phi(\dom \phi \cap \mathbb{R}_{\geq 0}) \subset K$.  

Next, we consider the perturbed version of $\Sigma$, denoted by $\Sigma_\epsilon$, and given by 
\begin{align} \label{eq.2}
\Sigma_\epsilon : \quad \dot{x} \in F(x) + \epsilon(x) \mathbb{B} \qquad  x \in \mathbb{R}^n.
\end{align} 
Following \cite{9683684}, we introduce the robust-safety notion studied in this paper. 
\begin{definition}[Robust safety]
$\Sigma$ is robustly safe with respect to $(X_o, X_u)$ if there exists $\epsilon \in \mathcal{C}_+$  such that $\Sigma_{\epsilon}$ in \eqref{eq.2} is safe with respect to $(X_o,X_u)$. 
The function $\epsilon$ is in this case named
 \textit{robust-safety margin}.
\end{definition}
We next recall the notion of a barrier function candidate.  
\begin{definition}
A scalar function $B : \mathbb{R}^n \rightarrow \mathbb{R}$ is a barrier function candidate with respect to $(X_o,X_u)$ if
\begin{align*} 
\begin{matrix} 
B(x) > 0 & \forall x \in X_u \quad \text{and} \quad B(x) \leq 0 & \forall x \in X_o. 
\end{matrix} 
\end{align*} 
\end{definition}
Note that a barrier function candidate 
$B$ defines the set 
\begin{align} \label{eq.4} 
K := \left\{ x \in \mathbb{R}^n : B(x) \leq 0 \right\}, 
\end{align} 
satisfying  $X_o \subset K$ and $K \cap X_u = \emptyset$.

\subsection{Sufficient conditions for robust safety} 

Since robust safety for $\Sigma$ coincides with safety for  
$\Sigma_\epsilon$,  when $\epsilon$ is a robustness margin,   then robust safety follows if 
\begin{enumerate} [label={C0)},leftmargin=*]
\item \label{item:C0} There exists a barrier function candidate $B : \mathbb{R}^n \rightarrow \mathbb{R}$ such that the set $K$ in \eqref{eq.4} is forward invariant for  $\Sigma_\epsilon$,  for some   $\epsilon :\mathbb{R}^n \rightarrow \mathbb{R}_{>0}$ continuous.  
\end{enumerate}

\blue{ The latter allows us to recall the following solution-dependent sufficient condition for robust safety \cite{9705088}. }

\begin{proposition}
$\Sigma$ is robustly safe with respect to $(X_o,X_u) \in \mathbb{R}^n \times \mathbb{R}^n$ provided that   
\begin{enumerate} [label={C1)},leftmargin=*]
\item \label{item:C2bis}  
\blue{There exists a barrier function candidate $B$ such that the set $K$ in \eqref{eq.4} is closed, there exists $U(\partial K)$ an open neighborhood of $\partial K$,  and  there exists 
$\epsilon \in \mathcal{C}_+$ such that,  along every solution (not necessarily maximal)  
$\phi$  to $\Sigma_{\epsilon}$ satisfying $\phi(\dom \phi) \subset U(\partial K)$,    the map 
$t \mapsto B(\phi(t))$ is non increasing. }
\end{enumerate} 
\end{proposition}

Note that \ref{item:C2bis} does not require from $B$ to be smooth neither $F$ to satisfy Assumption \ref{ass1}. However, it involves the solutions to $\Sigma_\epsilon$ (hence, it also involves the perturbation term $\epsilon$).   

\blue{We now recall from \cite{9683684,RubSafPI} a sufficient condition for robust safety that uses only the barrier function candidate $B$ and the nominal dynamics $F$.}

\begin{proposition}
Consider system $\Sigma$ such that  Assumption \ref{ass1} holds. $\Sigma$ is robustly safe with respect to $(X_o,X_u) \in \mathbb{R}^n \times \mathbb{R}^n$ provided that 
\begin{enumerate} [label={C2)},leftmargin=*]
\item \label{item:C3} There exists 
a continuously differentiable barrier function candidate $B$ such that
\begin{align} 
\langle \nabla B (x), \eta \rangle < 0 \qquad \forall \eta \in F(x), \quad \forall x \in \partial K.   \label{eq.2cbis}
\end{align}
\end{enumerate} 
\end{proposition}
\blue{Note that \eqref{eq.2cbis} involves only the barrier function candidate $B$,  which is now required to be continuously differentiable,  and the nominal dynamics $F$, which is now required to verify Assumption \ref{ass1}.} 

\subsection{Converse robust-safety theorems}

We start introducing the following two assumptions. 

\begin{assumption} \label{ass3} 
The map $F$ is continuous and $F(x)$ is nonempty, compact, and convex for all $x \in \mathbb{R}^n$.
\end{assumption}

\begin{assumption} \label{ass4-} 
$\cl(X_o) \cap X_u = \emptyset$. 
\end{assumption}

\blue{Using the latter two assumptions, we can formulate our first converse robust-safety theorem.}

\begin{theorem} \label{thm3}
Consider system $\Sigma$ that is robustly safe with respect to $(X_o,X_u)$ and such that Assumptions  \ref{ass3} and \ref{ass4-} hold.  Then,   \ref{item:C2bis} holds. 
\end{theorem}

The next converse robust-safety theorem uses the following assumption instead of Assumption \ref{ass4-}. 

\begin{assumption} \label{ass4} 
$\cl(X_o) \cap \cl (X_u) = \emptyset$. 
\end{assumption}   

\begin{theorem} \label{thm4}
Consider system $\Sigma$ that is robustly safe with respect to $(X_o,X_u)$ and such that Assumptions \ref{ass3}  and \ref{ass4} hold.  Then,   \ref{item:C3} holds.
\end{theorem}

\begin{remark} \label{remthm2}
Our proof of Theorem \ref{thm4} actually allows us to conclude that, for any $k \in \{1,2,...\}$, their exists a continuously-differentiable barrier certificate $B \in \mathcal{C}^k$ such that
$$ \langle \nabla B (x), \eta \rangle < -1 \qquad \forall \eta \in F(x),  \quad \forall x \in \partial K.  $$   
\end{remark}

\begin{remark}
\blue{ We believe that relaxing Assumption 2, and using Assumption 1 instead, would require a totally different approach to prove Theorems \ref{thm3} and \ref{thm4}.  Indeed,  a key property allowing our barrier-function construction is established in Lemma \ref{lem1} below.   This property does not hold when only Assumption \ref{ass1} is verified; 
see Example \ref{exp1}.   Similarly,  relaxing Assumptions \ref{ass4-} and \ref{ass4} would require using a completely different approach to prove Theorems \ref{thm3} and \ref{thm4},  respectively.  Indeed,  such separations between the two sets $X_o$ and $X_u$ are necessary to squeeze the zero-level set of the constructed barrier function between the two sets; see the forthcoming Section \ref{Sec.Separ} for more details. }
\end{remark}

\subsection{Application: Safety for self-triggered control systems}

Consider the control system $\Sigma_u$ given by
$$ \Sigma_u : \dot{x} = f(x,u)  \qquad  (x,u) \in  \mathbb{R}^n \times  \mathbb{R}^{m},  $$  
 where $f : \mathbb{R}^n \times  \mathbb{R}^{m} \rightarrow \mathbb{R}^n$ is a continuous function. Furthermore, we let   $\kappa: \mathbb{R}^n \rightarrow \mathbb{R}^{m}$ be a continuous feedback law such that the resulting closed-loop system 
 $$ \Sigma : \dot{x} =  F(x) := f(x,\kappa(x))  \quad x \in \mathbb{R}^n $$
 is safe with respect to $(X_o,X_u) \subset \mathbb{R}^n \times \mathbb{R}^n$.  
 
  In a self-triggered (ST) control framework \cite{di_benedetto_digital_2013},  we construct a monotonically increasing sequence  
  $\{t_i\}^{\infty}_{i=0} \subset \mathbb{R}_{\geq 0}$ and we force the controller to   remain constant between each two time samples $t_i$ and $t_{i+1}$,  i.e.,  
  $$ u(t) = \kappa(x(t_i)) \qquad  \forall t \in [t_i, t_{i+1}), \qquad i \in \mathbb{N}. $$  
  Hence,  a solution 
$\phi$ to the ST closed-loop system satisfies
\begin{align*} 
\dot{\phi}(t) =  F(\phi(t))  + \Gamma(\phi(t),\phi(t_i))   \quad \forall t \in [t_i, t_{i+1}), \quad \Gamma(\phi(t),\phi(t_i)) :=  f(\phi(t),\kappa(\phi(t_i))) - f(\phi(t),\kappa(\phi(t))).
\end{align*}

Our goal here is to address the following problem.  

\begin{problem} \label{probap}
\blue{Knowing that $\Sigma$ is robustly safe, prove the existence of a sequence $\{t_i\}^{\infty}_{i=0} \subset \mathbb{R}_{\geq 0}$ and $T > 0$ such that $t_{i+1} - t_i  > T$ for all $i \in \mathbb{N}$ and the ST closed-loop system is safe with respect to $(X_o,X_u)$.}  
\end{problem}

Inspired by \cite{9483096},  we will show that Theorem \ref{thm4} is key to solve Problem \ref{probap}. Indeed, Theorem \ref{thm4} allows us to formulate the following intermediate 
result. 

\begin{corollary} \label{coroo}
Consider the control system $\Sigma_u$ and the feedback law $\kappa$ such that the resulting closed-loop system $\Sigma$ is robustly safe with respect to $(X_o,X_u)$. 

Then, there exist a continuously-differentiable barrier certificate $B$ and
continuously-differentiable functions $\alpha : \mathbb{R}^n \rightarrow \mathbb{R}$ and $\gamma : \mathbb{R}^n \times \mathbb{R}^n \rightarrow \mathbb{R}$ such that
\begin{align} \label{eqConC-} 
\alpha(x) \geq 3/4 \quad \forall x \in \partial K,  \qquad 
\gamma (x,x) \leq 1/8 \quad \forall x \in \mathbb{R}^n, 
\end{align}
and
\begin{align} \label{eqConC}
\hspace{-0.6cm}
\langle \nabla B(x),  f(x,\kappa(y)) \rangle   \leq - \alpha(x) + \gamma(x,y) \qquad  \forall (x,y) \in \mathbb{R}^n \times \mathbb{R}^n.  
\end{align}
\end{corollary}

\begin{proof}
\blue{The existence of the 
continuously-differentiable barrier certificate $B$ is guaranteed by Theorem \ref{thm4}. Furthermore,  according to Remark \ref{remthm2}, we can choose the barrier certificate $B$ to satisfy
$$ \hat{\alpha}(x) := - \langle \nabla B(x) ,  F(x)  \rangle > 1 \quad \forall x \in \partial K \quad \text{and} \quad   \hat{\gamma}(x,x) := \langle \nabla B(x) ,   \Gamma(x,x)  \rangle  = 0 \quad \forall x \in \mathbb{R}^n. $$
As a result, using Whitney approximation theorem for continuous functions, we conclude that we can always find  continuously-differentiable functions $\gamma$ and $\alpha$ such that \eqref{eqConC-} holds and at the same time
 \begin{align} \label{eqprop31}
 \alpha(x) \leq  \hat{\alpha}(x) \quad  \forall x \in \mathbb{R}^n
 \quad \text{and}  \quad 
 \gamma(x,y)  \geq \hat{\gamma}(x,y) \quad \forall (x,y) \in \mathbb{R}^n \times \mathbb{R}^n.  
\end{align}
The choice in \eqref{eqprop31} allows us to verify 
\eqref{eqConC}. }
\end{proof}

To solve Problem \ref{probap},
we introduce the following assumption.

\begin{assumption} \label{assexp}
\blue{The set $\mathbb{R}^n \backslash X_u$ is bounded, and there exists $\tau>0$ such that the solutions to $\Sigma_y : \dot{x} = f(x,\kappa(y))$ starting from $y \in \mathbb{R}^n \backslash X_u$ cannot blow up on the time  interval $[0,\tau]$.}  
\end{assumption}

Under Assumption \ref{assexp}, we conclude that the zero sub-level set $K$ of the barrier certificate  
 is compact. Hence, there exist $\beta>0$ and $T_1>0$ such that the following properties hold. 
 
\blue{\begin{itemize}
\item[Pr1)] For each $y \in G:=  \{x \in K: |x|_{\partial K} \leq \beta\}$, we have $\alpha(y) - \gamma(y,y)  \geq \alpha(y) - 1/8  \geq 1/4$. 
This is true because $\alpha$ is continuously differentiable. 
\item[Pr2)] The solutions  to $\Sigma_y : \dot{x} = f(x,\kappa(y))$ starting from  $y \in L :=   \{x \in K: |x|_{\partial K} \geq \beta\}$
remain in $K$ on $[0,T_1]$.
This is because $f$ and $\kappa$ are continuous and $K$ is compact.
\end{itemize}}

The following result  addresses Problem \ref{probap}. 

\begin{proposition}
\blue{Consider the control system $\Sigma_u$ whose dynamics $f$ is continuous. Consider $(X_o,X_u) \subset \mathbb{R}^n \times \mathbb{R}^n$  and a continuous feedback law $\kappa$ such that Assumption \ref{assexp} holds.  Then,  a solution to Problem \ref{probap} is given by
$$ t_{i+1} := 
\left\{ 
\begin{matrix} 
t_i +  \max \{ T_1,  T_r(\phi(t_i)) \}  & \text{if} ~ \phi(t_i) \in L
\\
t_i +  T_r(\phi(t_i))    &  \text{if} ~ \phi(t_i) \in G,
\end{matrix} 
\right.
$$ 
where $L$ and $(T_1,G)$ are introduced in Pr1) and Pr2), respectively. Furthermore, the map $y \mapsto T_r(y)$ is defined, for all $y \in K$, as
\begin{align*}
T_r(y) & := 0 &  \text{if} ~ \alpha(y) - \gamma(y,y) \leq 0, 
\\
T_r(y) & := \left\{ 
\begin{matrix}
\tau &  \text{if} ~ M_r(\tau,y) \leq 0
\\ 
\min \left\{ \tau, \frac{2(\alpha(y) - \gamma(y,y))}{M_r(\tau,y)} \right\} & \mbox{otherwise}
\end{matrix}
\right.   &  \text{otherwise}, 
\end{align*}
where $\tau$ is introduced in Assumption \ref{assexp} and 
\begin{equation}
\label{eqMreq}
\begin{aligned}
M_r(\tau, y) & :=   \sup \{\langle \nabla_x \gamma(x, y), f(x, \kappa(y)) \rangle:  x \in R_{\Sigma_y}(\tau, y) \} 
+ \sup \{\langle - \nabla \alpha(x) , f(x, \kappa(y)) \rangle: x \in R_{\Sigma_y}(\tau, y)\},
\end{aligned}
\end{equation}
$R_{\Sigma_y}(\tau, y)$ is the set of points reached by the solutions to $\Sigma_y$ starting from $y$ over the window of time $[0,\tau]$, 
and the functions $\alpha$ and $\gamma$ are introduced in Corollary \ref{coroo}. }
\end{proposition}

\begin{proof}
Using Pr1), we conclude that $T_r(y)>0$ for all $y \in G$. Furthermore, using Pr2), we conclude that $t_{i+1} - t_i > 0$ for all $i \in \mathbb{N}$. 

Next, we show that any solution $\phi$ to $\dot{x} = f(x,\kappa(y))$,  starting from $y \in K$,  satisfies 
$ \phi([0,T_r(y)]) \subset K$. This is enough to conclude that the self-triggered closed-loop system is safety with respect $(X_o,X_u)$.
To this end, we note that $\phi$ is locally absolutely continuous,  and  since $\alpha$ and $\gamma$ are continuously differentiable, it follows that  $t \mapsto \alpha(\phi(t))$ and $t \mapsto \gamma(\phi(t),y)$ are also locally absolutely continuous.   Hence,  for each $t \in \dom \phi$,  there exists a sequence 
 $\{\tau_n\}_{n=0}^N\subset [0,t]$,  with  $N \in \mathbb{N}^* \cup \{\infty\}$, such that 
 $$ \lim_{n\rightarrow N} \tau_n = t, \qquad   \tau_n - \tau_{n-1} > 0, $$  
 and the maps $t \mapsto \alpha(\phi(t))$ and $t\mapsto \gamma(\phi(t), y)$ are  differentiable on each $(\tau_{n-1}, \tau_{n})$. 

 \textcolor{blue}{Consider the map} $\bar{\gamma}(\cdot) := \gamma(\cdot,y)$ 
 and note that
$$ \bar{\gamma}(\phi(t)) - \bar{\gamma}(y) = \sum_{n=1}^N  \left( \bar{\gamma}(\phi(\tau_n)) - \bar{\gamma}(\phi(\tau_{n-1})) \right), \quad \alpha(\phi(t)) - \alpha(y)  = \sum_{n=1}^N \left[\alpha(\phi(\tau_n)) - \alpha(\phi(\tau_{n-1}))\right].  $$
As a result, using the classical mean-value theorem,  we conclude that, for each $n\in\{1,2,...,N\}$, there exist $c_n, d_n \in (\tau_{n-1}, \tau_{n})$ such that
\begin{align*}
\bar{\gamma}(\phi(t)) - \bar{\gamma}(y)  = \sum_{n=1}^{N} \left( \frac{d}{dt} \bar{\gamma}(\phi(t)) \Big\vert_{t=c_n} (\tau_{n} - \tau_{n-1}) \right),
\quad 
\alpha(\phi(t)) - \alpha(y)  = \sum_{n=0}^{N}
\left(\frac{d}{dt} \alpha(\phi(t)) \Big\vert_{t=d_n}  (\tau_{n+1} - \tau_n)\right). 
\end{align*}
As a result,  when $t \in [0,\tau]$,  we conclude that
\begin{equation}
    \label{eqciteit}
    \begin{aligned} 
    \bar{\gamma}(\phi(t))  - \bar{\gamma}(y) & \leq  t \sup\{\langle \nabla\bar{\gamma}(x),f(x,\kappa(y)) \rangle:  x \in R_{\Sigma_y}(\tau,y) \},
\\
-\alpha(\phi(t)) + \alpha(y)    
   & \leq t \sup \{\langle - \nabla \alpha(x), f(x,\kappa(y)) \rangle :   x \in R_{\Sigma_y}(\tau,y) \}.
\end{aligned}
\end{equation}
Next, we note that 
 $$ \frac{d}{dt} B(\phi(t)) = \langle \nabla B(\phi(t)), \dot{\phi}(t) \rangle \qquad   \text{for almost all} ~ t \in \dom \phi.  $$  
Integrating the previous equality 
from $0$ to $t \leq \tau$, we obtain
\begin{align*} 
 B(\phi(t)) -B(y) & \leq \int_0^t [-\alpha(\phi(s)) + \bar{\gamma} (\phi(s))]  ds  \leq \int_0^t [-\alpha(y) + \bar{\gamma}(y) + s M_r(\tau, y)] ds    
\\ & 
\leq  - t (\alpha(y) - \bar{\gamma}(y))  + \frac{t^2}{2} M_r(\tau, y) \qquad \forall t \in [0,\tau].
\end{align*} 
To obtain the latter inequalities, we used \eqref{eqConC}, \eqref{eqciteit}, and \eqref{eqMreq}. 

Now, we note that, when $\alpha(y) - \bar{\gamma}(y) > 0$ and $M_r(\tau,y) > 0$, it follows that
$$ B(\phi(t)) - B(y) \leq 0 \qquad  \forall t \in \left[ 0, 2 \frac{ \left(\alpha(y) - \bar{\gamma}(y) \right)}{M_r(\tau,y)} \right] \cap [0,\tau]. $$ 
 Otherwise,  when $\alpha(y) - \bar{\gamma}(y) > 0$ and $M_r(\tau,y) \leq 0$, we conclude that
$$ B(\phi(t)) - B(y) \leq 0 \qquad  \forall t \in [0,\tau]. $$ 
The latter is enough to conclude that the proposed triggering sequence guarantees safety for the resulting self-triggered 
closed-loop system.

In the rest of the proof, 
we show the existence of $T > 0$ such that $ t_{i+1} - t_i  > T$ for all  $i \in \mathbb{N}$. To do so,  it is enough to show that the maps $y \mapsto T_r(y)$ is lower semicontinuous on $G$. Indeed, since it is already positive on $G$ and $G$ is compact,  it would follow using \cite[Theorem B.2]{puterman2014markov} that $T_r$ reaches its minimum on $G$. As a result, we can take $T := \min \{ \min \{T_r(y): y \in G \}, T_1 \}$. 

To show that $T_r$ is lower semicontinuous on $G$, we start noting that the set-valued map $y \mapsto R_{\Sigma_y}(\tau,y)$ is upper semicontinuous with compact images on $G$; see Lemma \ref{lem1-} in the Appendix. Next,  we use \cite[Theorem 1.4.16]{aubin2009set},  under smoothness properties of $\gamma$ and $\alpha$,  to conclude that the single-valued map $y \mapsto M_r(\tau,y)$ is upper semicontinuous on $G$. \blue{ It is also locally bounded on $G$ since so is the set-valued map $y \mapsto R_{\Sigma_y}(\tau,y)$.} 

Now,  since  $\alpha$ is positive and continuous on $G$,  we conclude that $y \mapsto \frac{\alpha(y) - \bar{\gamma}(y)}{M_r(\tau,y)}$ is lower semicontinuous and positive on $G$. To complete the proof, we consider a sequence $ \{x_i\}^{\infty}_{i=0} \subset G$ that converges to $x_o \in G$. Since $y \mapsto M_r(\tau,y)$ is upper semicontinuous, we conclude that $\limsup_{i \rightarrow \infty} M_r(\tau,x_i) \leq M_r(\tau,x_o)$.   \textcolor{blue}{Moreover, by selecting an appropriate subsequence, we can assume, without loss of generality, that} $ \liminf_{i \rightarrow \infty} T_r(x_i)  =    \lim_{i \rightarrow \infty} T_r(x_i)  = a > 0$. Also, since $M_r$ is locally bounded,  \textcolor{blue}{one can select another subsequence to conclude the existence of} $\beta \in \mathbb{R}$ such that 
$  \lim_{i \rightarrow \infty} M_r(\tau,x_i) = \beta \leq M_r(\tau,x_o)$. 
 \textcolor{blue}{To finalize the proof, we} distinguish the following three scenarios:  
\begin{enumerate}
\item  When $\beta < 0$, we conclude that 
$M_r(\tau,x_i) < 0$ for all $i \in \mathbb{N}$ sufficiently large. In this case,  $T_r(x_i) = \tau$ for all $i \in \mathbb{N}$ sufficiently large, thus, $\lim_{i \rightarrow \infty} T_r(x_i) =  
\tau \geq T_r(x_o)$.   
\item  When $\beta > 0$, we conclude that $M_r(\tau,x_i) > 0$ for all $i \in \mathbb{N}$ sufficiently large.   In this case, we conclude that 
$T_r(x_i) = \min \left\{ \tau,  \frac{\alpha(x_i) - \bar{\gamma}(x_i)}{M_r(\tau,x_i)}  \right\}$  for all $i \in \mathbb{N}$ large. Hence,   
$
\lim_{i \rightarrow \infty} T_r(x_i)  =  
\min \left\{\tau, \lim_{i \rightarrow \infty} \frac{ \alpha(x_i) - \bar{\gamma}(x_i)}{ M_r(\tau,x_i)} \right\} 
\geq \min \left\{ \tau,   \frac{\alpha(x_o) - \bar{\gamma}(x_o)}{M_r(\tau,x_o)} \right\} \geq T_r(x_o)$.
\item   When $\beta = 0$, we conclude that 
$| \frac{\alpha(x_i) - \bar{\gamma}(x_i)}{ M_r(\tau,x_i)} | \geq \tau$ for all $i \in \mathbb{N}$ sufficiently large. Hence, 
$T_r(x_i) = \tau$ for all $i \in \mathbb{N}$ sufficiently large.
\end{enumerate}
\end{proof}

\section{Preparatory material} \label{Sec.4}

In this section, we present key intermediate results that allow us to construct a smooth barrier certificate for robustly-safe systems. 

\subsection{Perturbed differential inclusions} 

Given $x \in \mathbb{R}^n$, 
we investigate conditions on $F$, under which,  for any  $\epsilon \in \mathcal{C}_+$,  there exists $T > 0$ 
such that  $R^b_{\Sigma}(T,x) \subset  \text{int}\left( R^b_{\Sigma_\epsilon}(T,x) \right)$. A similar result requiring $F$ to be locally Lipschitz can be found in \cite{puri1995varepsilon} and \cite{ratschan2018converse}.

\begin{lemma} \label{lem1}
Consider system $\Sigma$  such that Assumption  \ref{ass3} holds.  Then,  for each $x \in \mathbb{R}^n$ and for each  $\epsilon \in \mathcal{C}_+$,  
there exists $T>0$ such that,  for each 
$t \in (0,T]$,  there exists $\delta > 0$ such that 
\begin{align}
y + \delta \mathbb{B} & \subset R^b_{\Sigma_\epsilon}(t,x) \qquad \forall y \in R^b_{\Sigma}(t,x) \backslash \{x\}, 
\label{eqreachcover} \\
x + \delta \mathbb{B} & \subset 
 R^b_{\Sigma_\epsilon}(-t,y) = R^b_{\Sigma^-_\epsilon}(t,y) \qquad \forall y \in R^b_{\Sigma}(t,x) \backslash \{x\}, \label{eqreachcover1} 
\end{align}
where $\Sigma^- :   \dot{x} \in - F(x)$ with $x \in \mathbb{R}^n$.  
\end{lemma}

\begin{proof}
Given $x \in \mathbb{R}^n$ and  $\epsilon \in \mathcal{C}_+$,  we pick $\Delta \in (0,1)$ such that,  for each $x_1,x_2 \in x + \Delta \mathbb{B}$ and for each $f_1 \in F(x_1)$, there exists $f_2 \in F(x_2)$ such that
\begin{align} \label{eqcontin}
|f_1 - f_2| \leq \underline{\epsilon}/2, \quad \underline{\epsilon} := \min \{\epsilon(y) : y \in x + \mathbb{B} \}. 
\end{align}  
The latter is possible since the set-valued map $F$ is assumed to be continuous. 

Now,  using Lemma \ref{lemprepre} in the Appendix,  we conclude the existence of $b$ and $\bar{T}>0$ such that $R_{\Sigma} (\bar{T},  (x + b \mathbb{B}))$ and 
$R_{\Sigma^-} (\bar{T},  (x + b \mathbb{B}))$
are bounded. As a result,  we invoke Lemma \ref{lem1-} in the appendix to conclude that  $R_\Sigma$ and $R_{\Sigma^-}$ are outer semicontinuous and locally bounded  on 
$ [0,\bar{T}] \times (x + b \mathbb{B})$,  which implies,  according to Remark \ref{remplus},  that $R_\Sigma$ and $R_{\Sigma^-}$ are upper semicontinuous and have  compact images on 
$[0,\bar{T}] \times (x + b \mathbb{B})$. 
Hence,   we can find $T \in (0,\bar{T}]$ sufficiently small 
such that 
\begin{align} 
R_{\Sigma}(T,x) & \subset 
\left(x + \frac{\Delta}{2} \mathbb{B} \right), 
\label{eqcontin1a}
\\
R_{\Sigma^{-}}(T,R_{\Sigma}(T,x)) & \subset 
\left(x + \frac{\Delta}{2} \mathbb{B} \right).
\label{eqcontin1b}
\end{align}

To prove \eqref{eqreachcover}, we let $\delta \in (0, \Delta/2]$ and we let $y \in R^b_{\Sigma}(t,x) \backslash \{ x \}$ for some $t \in (0,T]$.  Hence,  there exists a solution $\phi \in \mathcal{S}_{\Sigma}(x)$ such that $\phi(0) = x$ and $\phi(t) = y$.

Now,  given $z \in y + \delta \mathbb{B}$, we consider the function $\eta : [0,t] \rightarrow \Delta \mathbb{B}$ given by
$$ \eta(s) := \phi(s) - \frac{s}{t} (y-z) \qquad \forall s \in [0,t]. $$
Note that
$$ \dot{\eta}(s) = \dot{\phi}(s) - \frac{1}{t} (y-z) \qquad \text{for almost all} \quad s \in [0,t]. $$
Next, for almost all $s \in [0,t]$, we let 
$f_\eta(s) \in F(\eta(s))$ be such that 
$ |f_\eta(s) - \dot{\phi}(s)| \leq \underline{\epsilon}/2. 
$
This is possible using \eqref{eqcontin}.  Hence, we conclude that 
\begin{align*} 
\dot{\eta}(s) & = f_\eta(s) + (\dot{\phi}(s) - f_\eta(s)) - \frac{1}{t} (y-z)  \in F(\eta(s)) + 
\left( \frac{\underline{\epsilon}}{2} + \frac{\delta}{t} \right) \mathbb{B} \qquad \quad \text{for almost all} ~~ s \in [0,t]. 
\end{align*} 
As a result, by taking  $\delta := \frac{\min \{ \underline{\epsilon}t, \Delta \}}{2}$,  we conclude that
$$
\dot{\eta}(s)  \in F(\eta(s)) + 
\underline{\epsilon} \mathbb{B} \subset F(\eta(s)) + 
\epsilon(\eta(s)) \mathbb{B} \qquad \text{for almost all} \quad s \in [0,t]. 
$$ 
Hence,  $\eta : [0,t] \rightarrow x + \Delta \mathbb{B}$ is a solution to 
$\Sigma_{\epsilon}$ verifying $\eta(0) = x$ and $\eta(t) = z$,  which proves \eqref{eqreachcover} since $z$ can be arbitrary in $y + \delta \mathbb{B}$. 

The proof of \eqref{eqreachcover1} follows using the same steps used to prove \eqref{eqreachcover},  while invoking \eqref{eqcontin1b} instead of \eqref{eqcontin1a}. 
\end{proof}

It is important to note that,  when only  
Assumption \ref{ass1} is verified;  namely,  when $F$ is not required to be continuous,  then we can find a system $\Sigma$,  $x \in \mathbb{R}^n$,  and  $\epsilon \in \mathcal{C}_+$  such that 
$$R^b_{\Sigma}(t,x) \backslash \{x\}   \not\subset  \text{int}(R^b_{\Sigma_\epsilon}(t,x)) \qquad \forall t > 0.  $$ 

\begin{example} \label{exp1}
Consider system $\Sigma$ with $n = 2$ and 
$$ F(x) := \left\{ 
\begin{matrix} 
[0 \quad 1]^\top  & \text{if}~ x_2 < 0
\\
\co \{[0 \quad 1]^\top, [-1 \quad 0]^\top\} & \text{otherwise}.   
\end{matrix} 
\right. $$
Let the constant function $\epsilon : \mathbb{R}^2 \rightarrow \{1/2\}$.

Note that $\Sigma$ admits a solution $\phi$ starting from $x_o = 0$ that is given by 
$ \phi(t) := [-t \quad 0]^\top$ for all $t \geq 0$.  
Note that,  for each $t>0$,  we have 
$ \phi(t) \in R^b_{\Sigma}(t,x_o) \backslash \{x_o\} $. 

Now, given $t \geq 0$,  we show that
\begin{align*} 
 \phi(t) - [0 \quad 1/i]^\top \notin R_{\Sigma_\epsilon}(t,x_o) \qquad  \forall i \in \{1,2,... \}.  
 \end{align*}
To do so, it is enough to show that the set 
$K := \{x \in \mathbb{R}^2 : x_2 \geq 0 \}$ is forward invariant for $\Sigma_{\epsilon := \frac{1}{2}}$. 
We show the latter by applying  Lemma \ref{lemA14} in the Appendix,  via verifying 
\begin{align} \label{eqteng} 
F(x) + \mathbb{B}/2 \subset C_K(y)  \qquad  \forall y \in \text{Proj}_{K}(x) \qquad \forall  x \in \mathbb{R}^2 \backslash K. 
\end{align}
To show \eqref{eqteng},  
we note that 
$\mathbb{R}^2 \backslash K = \{ x \in \mathbb{R}^2 : x_2 < 0 \}$ and that the projection of $x := [x_1 \quad x_2]^\top \in \mathbb{R}^2 \backslash K$ on $K$ is $y := [x_1 \quad 0]^\top$. Furthermore,   using 
\cite[Lemma 3]{draftautomatica},  we conclude that 
$$ C_K(y) = \{v \in \mathbb{R}^2 : v_2 \geq 0 \} \qquad \forall y \in \partial K.  $$  
Now,  for each $x \in \mathbb{R}^2 \backslash K$,  we have  
$$ F(x) + \mathbb{B}/2 = 
\{[ \alpha/2 \quad (\beta/2) + 1 ]^\top : \sqrt{\alpha^2 + \beta^2} \leq 1 \}. $$ 
As a result, for each $v := [v_1 \quad v_2]^\top \in F(x) + 
\mathbb{B}/2$, we conclude that $v_2 \geq 0$, which means that $v \in C_K(y)$ for all $y \in \text{Proj}_{K}(x)$.   Hence, \eqref{eqteng} follows.
\end{example}
 
\subsection{The time to impact contractive and recurrent sets}

We start recalling a definition of forward contractivity of a closed set for a differential inclusion $\Sigma$ \cite{Blanchini:1999:SPS:2235754.2236030}.  

\begin{definition}[Forward contractivity] 
A closed subset $K \subset \mathbb{R}^n$ is forward contractive for $\Sigma$ if it is forward invariant for $\Sigma$ and,  for every $x_o \in \partial K$ and for every solution $\phi \in \mathcal{S}_{\Sigma}(x_o)$,  there exists 
$T > 0$ such that $\phi(t) \in \mbox{int}(K)$ for all $t \in \dom \phi \cap (0,T]$. 
\end{definition}

\blue{Next, we recall a definition of global reccurence 
of a subset $K \subset \mathbb{R}^n$ for a differential inclusion $\Sigma$ }
\cite{SUBBARAMAN201654}. 
\begin{definition}[Global recurrence] 
\blue{A set $K \subset \mathbb{R}^n$ is globally recurrent for $\Sigma$ if,  for each $\phi \in \mathcal{S}_{\Sigma}(x)$ with  
$x \in \mathbb{R}^n$,  there exists $t \in \dom \phi \cap \mathbb{R}_{\geq 0}$ such that $\phi(t) \in \mbox{int}(K)$.} 
\end{definition}

\blue{Finally, we deduce a definition   of local recurrence,  which we will exploit in our work.}

\begin{definition}[Local recurrence]
A set $K \subset \mathbb{R}^n$ is locally recurrent for $\Sigma$ on a neighborhood of $\partial K$,  denoted by $U(\partial K)$,  
if, for each $\phi \in \mathcal{S}_{\Sigma}(x)$ with $x \in U(\partial K)$, there exists $t \in \dom \phi \cap \mathbb{R}_{\geq 0}$ such that 
$\phi(t) \in \mbox{int}(K)$. 
\end{definition}

\begin{remark}
We note that local (respectively, global) recurrence of the set $K \subset \mathbb{R}^n$ implies,  respectively,   
local (resp. global) recurrence of $\cl(K)$ and vice versa.
\end{remark}

\blue{We recall here a variant of   the time-to-impact function (also known as the hitting-time function \cite{Aubin:1991:VT:120830}), denoted by $B_K : \mathbb{R}^n \rightarrow \mathbb{R}$, with respect to a closed subset $K \subset \mathbb{R}^n$  along the solutions to $\Sigma$.   We will show that such a function is well defined when $K$ is locally recurrent for $\Sigma$. Furthermore, when $K$ is additionally contractive, then $B_K$ is shown to be strictly decreasing along the solutions that remain close to $\partial K$.}  

We start introducing the  set-valued map 
$\hat{B}_{K} : U_1 \rightrightarrows  \mathbb{R} \cup \{\pm \infty\}$ 
given by 
\begin{align} \label{eqhatB}
\hat{B}_K(x) := 
\left\{ 
\begin{matrix}
\inf \{ T_{K} (\phi) : \phi \in \mathcal{S}_{\Sigma}(x) \} 
& \text{if $x \in \text{int}(K)$}
\\
0 & \text{if $x \in \partial K$}
\\
\sup \{ T_K (\phi) : \phi \in \mathcal{S}_{\Sigma}(x) \} & \text{otherwise},
\end{matrix} \right.
\end{align}
where,  for each $\phi \in \mathcal{S}_{\Sigma}(x)$,  the functional $T_K : \mathcal{S}_{\Sigma}(U_1) \rightarrow \mathbb{R} \cup \{  \pm \infty \} $ is given by
$$ T_K (\phi) :=
\argmin \{ |t| : t \in \dom \phi,~ \phi(t) \in \partial K \}.  $$

\blue{Roughly speaking,  $T_K(\phi)$ associates to each $\phi \in \mathcal{S}_{\Sigma}(x)$ the time, with the smallest norm,  at which $\phi$ hits $\partial K$.   
Then, $\hat{B}_K(x)$ is defined to be the smallest 
(respectively,   the largest) value among such times,  over all the $\phi$s in $\mathcal{S}_{\Sigma}(x)$,  when $x \in \text{int}(K)$ (respectively,  when $x \in \mathbb{R}^n \backslash K$).   }

\blue{ In \cite{Aubin:1991:VT:120830},  functionals similar to $T_K$ are introduced under the name hitting- and the exit-time functionals, denoted by $\theta_K$ and 
$\tau_K$, respectively.
Note that,  under Assumption \ref{assnew} below,  the latter two functionals coincide.  Furthermore,  under the same assumption,  we conclude that 
$$ T_{K} (\phi) =  \theta_{\mathbb{R}^n \backslash K}(\phi) = \tau_{\mathbb{R}^n \backslash K}(\phi) \qquad \forall \phi \in \mathcal{S}_\Sigma(x) ~ \text{and} ~ x \in \mathbb{R}^n \backslash K.  $$
}

\begin{assumption} \label{assnew}
There exists a closed subset 
$U_{1} \subset \mathbb{R}^n$ such that $\partial K \subset \text{int}(U_1)$ and
\begin{enumerate} [label={A\ref{assnew}\arabic*)},leftmargin=*]
\item \label{item:71--} The set $K$ is locally recurrent for $\Sigma$ on $U_1$.
\item \label{item:72--} The set $\mathbb{R}^n\backslash K$ is locally recurrent for $\Sigma^-$ on $U_1$.
\item \label{item:73--} The set $K$ is forward contractive for $\Sigma$.
\item \label{item:74--} The set $\mathbb{R}^n \backslash \mbox{int}(K)$ is forward contractive for $\Sigma^-$.
\end{enumerate}
\end{assumption}

\blue{
Using \cite[Proposition 4.2.4]{Aubin:1991:VT:120830},  we can guarantee that 
$\hat{B}_K(x) = \sup\{T_K(\phi) : \phi \in \mathcal{S}_\Sigma(x)\}$
is upper semicontinuous on $\mathbb{R}^n \backslash K$ provided that $F$ is strict Marchaud.
In the following lemma, we establish the same conclusions, among others, using only Assumptions \ref{ass1} and \ref{assnew} ($F$ is not required to be strict Marchaud). }

\begin{lemma} \label{lem4-}
Consider system $\Sigma$ such that Assumption \ref{ass1} holds.  \blue{Furthermore,  consider a closed subset $K \subset \mathbb{R}^n$ for which Assumption \ref{assnew} holds.}  
Then,  the following properties are satisfied.
\begin{enumerate} [label={P\ref{lem4-}\arabic*)},leftmargin=*]
\item \label{item:71-}   $\hat{B}_{K}$ is well defined on $U_1$; i.e.,  for each $x \in U_1$,  $\hat{B}_{K}(x)$ exists and it is finite.
\item \label{item:72-}  $\hat{B}_K$ is upper 
semicontinuous on $U_1 \backslash \text{int}(K)$.  
\item \label{item:73-} $\hat{B}_K$ is lower semicontinuous on $K \cap U_1$.  
\end{enumerate}
\end{lemma}     

\begin{proof}
We first use \ref{item:71--} and \ref{item:72--} to conclude that each maximal solution 
$\phi \in \mathcal{S}_{\Sigma}(x)$ with $x \in U_1$ reaches $\partial K$ at some time, which can be either positive or negative.  Hence,  $T_{K}(\phi)$ exists and it is finite for all 
$\phi \in \mathcal{S}_{\Sigma}(U_1)$.   

Next,  we will show that the map $x \mapsto T_K(\mathcal{S}_\Sigma(x))$ is locally bounded on $U_1 \backslash K$.   
To find a contradiction,   we pick $x_o \in U_1 \backslash K$, a sequence 
$\{ x_i \}^{\infty}_{i=1} \subset  U_1 \backslash K$ such that  $ \lim_{i \rightarrow \infty} x_i = x_o$, and a sequence 
$\{ \phi_i \}^{\infty}_{i=1}$ of solutions (not necessarily maximal) to $\Sigma$ such that $ \phi_i(0) = x_i$,   $\dom \phi_i = [0,T_{K}(\phi_i)]$ for all $i \in \{1,2,...\}$, and $ \lim_{i \rightarrow \infty} T_{K} (\phi_i) = +\infty$.   Without loss of generality,  we also assume that $i \mapsto T_K(\phi_i)$ is strictly increasing. 

Since $F$ is locally bounded,  using \red{Lemma \ref{lemprepre} in the Appendix},  we conclude the existence of $\bar{T} \in \mathbb{R}_{>0} \cup \{+ \infty\}$ the largest time  such that,  on any interval $[0,T] \subset [0,\bar{T})$,   the sequence $\{ \phi_i \}^{\infty}_{i=0}$ is uniformly bounded.   As a result,  by passing to an appropriate subsequence and \blue{using \cite[Theorem 5.29.]{goebel2012hybrid}} recursively on each closed interval $[0,T] \subset [0,\bar{T})$,  we conclude the existence of $\phi : [0,\bar{T}) \rightarrow \mathbb{R}^n$   solution  to $\Sigma$ starting from $x_o$ such that 
\begin{align} \label{eqlimit-}
\text{lim}_{i \rightarrow \infty} \phi_i(t) = \phi(t) \qquad \forall t \in [0,\bar{T}). 
\end{align}
Now,  since $x_o \in U_1 \backslash K$,  we conclude the existence of $\alpha \in (0, \bar{T})$ such that 
$T_{K} (\phi) < \alpha$. \blue{Indeed, $\alpha$ must be smaller than  $\bar{T}$ because  $T_K(\phi)$ exists and is finite, and when $\bar{T}$ is finite, we necessarily have  $\lim_{t \rightarrow \alpha} \phi(t) = + \infty$.}
Next, by forward contractivity of 
the set $K$ for $\Sigma$,  we conclude that 
$\phi(\alpha) \in \mbox{int}(K)$; thus, there exists 
$\beta>0$ such that 
$ | \phi(\alpha) |_{\mathbb{R}^n \backslash K}  \geq \beta$.
However, there exists $i^* \in \mathbb{N}$ such that 
$ \phi_i (\alpha) \notin K$ for all $i \geq i^*$,  
which implies that $ |\phi(\alpha) - \phi_i(\alpha)| \geq \beta$ for all $i \geq i^*$.    The latter contradicts \eqref{eqlimit-}; hence, the map $x \mapsto T_K(\mathcal{S}_\Sigma(x))$ is locally bounded and thus \ref{item:71-} follows as a direct consequence. 

To  prove \ref{item:72-},  we start re-expressing $\hat{B}$ as
$$  \hat{B}(x) :=
 \sup \{ t_\phi : t_\phi \in T_{K} (\mathcal{S}_{\Sigma}(x)) \} \qquad \forall x \in U_1 \backslash \text{int}(K). 
$$
Furthermore, we propose to show that the set-valued map 
$x \mapsto T_{K} (\mathcal{S}_{\Sigma}(x))$ is outer 
semicontinuous on $U_1 \backslash \text{int}(K)$. For this, we consider a sequence $\{x_i\}^\infty_{i=1} \subset U_1 \backslash \text{int}(K)$ that converges to $x_o \in U_1 \backslash \text{int}(K)$ and a sequence $\{t_i\}^\infty_{i=1} \subset \mathbb{R}_{\geq 0}$ that converges to $t_o \in \mathbb{R}_{\geq 0}$ such that 
$ t_i \in T_{K} (\mathcal{S}_{\Sigma}(x_i))$  for all i $\in \{1,2, ...\}$,   and we show that $t_o \in T_{K} (\mathcal{S}_{\Sigma}(x_o))$.  
Indeed,   let a sequence of solutions $ \{ \phi_i \}^{\infty}_{i =1}$,  with $\dom \phi_i = [0,  t_i]$, $t_i = T_K(\phi_i)$,  and $\phi_i(0) = x_i$ for all $i \in \{1,2,...\}$.  Here,  we distinguish between two cases.   
\begin{itemize}
\item When the sequence  $\{ \phi_i \}^{\infty}_{i =1}$ is uniformly bounded,  using \cite[Theorem 5.29.]{goebel2012hybrid}, we conclude the existence of a solution $\phi : [0,t_o] \rightarrow \mathbb{R}^n$ that is the graphical limit of an appropriate subsequence of  $\{ \phi_i \}^{\infty}_{i =1}$.  Next,  we show that 
$t_o = T_K(\phi) $.  To find a contradiction,  we assume that  $t_o \neq  T_K(\phi)$. Hence,  using \ref{item:71--} and \ref{item:73--},  it follows that there exists $\delta > 0$ such that either $|\phi(t_o)|_K \geq \delta$ or 
$|\phi(t_o)|_{\mathbb{R}^n \backslash K} \geq \delta$. 
On the other hand,  \blue{by graphical convergence and continuity of solutions}, we conclude that
$ \lim_{i \rightarrow \infty} \phi_i(t_i) = \phi(t_o)$.  Hence,  there exists $i^\star \in \{1,2,...\}$ such that $|\phi_i(t_i) - \phi(t_o)| \leq \delta/2$ for all $i \geq i^\star$.  The latter implies that 
$|\phi_i(t_i)|_{\partial K} \geq \delta/2$  for all $i \geq i^\star$,   and the contradiction follows. 

\item If the sequence is not uniformly bounded, we   exploit local boundedness of $F$ and \red{Lemma \ref{lemprepre} in the Appendix}, to conclude the existence of $\bar{T} \in (0,  t_o]$ the largest time such that the sequence $\{ \phi_i \}^{\infty}_{i =1}$ is uniformly bounded on any $[0, T] \subset [0,\bar{T})$.  
\blue{Using \cite[Theorem 5.29.]{goebel2012hybrid} recursively on each interval $[0, T] \subset [0, \bar{T})$},  we conclude the existence of an unbounded  solution $\phi : [0,\bar{T}) \rightarrow \mathbb{R}^n$  such that,  after passing to an appropriate subsequence,  we obtain  
 \begin{align} \label{eqlimit}  
 \lim_{i \rightarrow \infty} \phi_i(t)  = \phi(t) \qquad \forall t \in [0,\bar{T}).   
 \end{align}
 In this case,  we conclude that $\phi$ must reach $\partial K$ at some $ \alpha := T_K(\phi) \in [0, \bar{T}) < t_o$.  Since $K$ is contractive,  then we can find 
$\beta \in (\alpha,  t_o)$ such that 
$|\phi(\beta)|_{\mathbb{R}^n \backslash K} = \gamma$ for some $\gamma > 0$.  Now,  using \eqref{eqlimit}, we conclude that $ \lim_{i \rightarrow \infty} |\phi_i(\beta)|_{\mathbb{R}^n \backslash K}  = \gamma  > 0$ with 
$\beta < t_o$; hence, $\beta < t_i$ for $i$ large enough. This yields to a contradiction implying that $\alpha = t_o$.   
\end{itemize}

Having the map $x \mapsto T_K(\mathcal{S}_\Sigma(x))$ locally bounded and  outer semicontinuous on $U_1 \backslash \text{int}(K)$ implies that the map $x \mapsto T_K(\mathcal{S}_\Sigma(x))$ is upper semicontinuous with compact images on $U_1 \backslash \text{int}(K)$; see Remark \ref{remplus}.   Therefore,  using  \cite[Theorem 1.4.16]{aubin2009set},  we conclude that $\hat{B}_{K}$ is upper semicontinuous on $ U_1 \backslash \text{int}(K)$.

Finally,  to prove \ref{item:73-},  it is enough to 
notice that   
\begin{align*}
\inf \{ T_{K} (\phi) : \phi \in \mathcal{S}_{\Sigma}(x) \}  & =  - \sup\{ - T_{K} (\phi) : \phi \in \mathcal{S}_{\Sigma}(x) \}  =  - \sup\{  T_{K} (\phi) : \phi \in \mathcal{S}_{\Sigma^-}(x) \}.
\end{align*}
Therefore,  the previous arguments apply,  under \ref{item:72--} and \ref{item:74--},   to conclude that $-\hat{B}_K$ is upper semicontinuous on $K \cap U_1$; thus, 
 $\hat{B}_K$ is lower semicontinuous on $K \cap U_1$.
\end{proof}

Next, we introduce the function $B_K : \mathbb{R}^n \rightarrow \mathbb{R}$, which extends $\hat{B}_K$ to the entire $\mathbb{R}^n$ and is given by
\begin{equation}
\label{eqbarcand}
\begin{aligned} 
\hspace{-0.3cm} B_K(x) := 
\left\{
\begin{matrix}
\inf \{ \hat{B}_K(y) : y \in \text{Proj}_{U_1}(x)  \}   & \hspace{-0.2cm} \text{if} ~ x \in \text{int}(K)
\\
0 & \hspace{-0.2cm} \text{if $x \in \partial K$}
\\
\sup \{ \hat{B}_K(y) : y \in \text{Proj}_{U_1}(x) \} & \hspace{-0.2cm} \text{otherwise}.
\end{matrix}
\right. 
\end{aligned}
\end{equation}

\begin{lemma} \label{lem5-}
Consider system $\Sigma$ such that Assumption \ref{ass1} holds.  \blue{Furthermore,  consider $(X_o,X_u) \subset \mathbb{R}^n \times \mathbb{R}^n$ and a closed subset $K \subset \mathbb{R}^n$  such that $X_o \subset K$, $X_u \cap K = \emptyset$, and 
Assumption \ref{assnew} holds}.   

Then, the following properties are satisfied.  
\begin{enumerate} [label={P\ref{lem5-}\arabic*)},leftmargin=*]
\item \label{item:B1--}  $B_K$ is a barrier function candidate with respect to $(X_o,X_u)$ and $K = \{ x \in \mathbb{R}^n : B_K(x) \leq 0\}$. 
\item \label{item:B3--} $B_K$ is upper semicontinuous on $\mathbb{R}^n \backslash \text{int}(K)$.
\item \label{item:B4--} $B_K$ is lower semicontinuous on $K$.
\item \label{item:B5--} Along each solution $\phi$ to $\Sigma$ (not necessarily maximal) such that $\phi(\dom \phi) \subset U_1$,  the map $t \mapsto B_{K}(\phi(t))$ is  decreasing. 
In particular, we have 
$$ B_{K}(\phi(t')) - B_{K}(\phi(t)) \leq  t-t'  \quad \forall t, t' \in \dom \phi  ~~ \text{with} ~~ t \leq t'.  $$
\end{enumerate}
\end{lemma}

\begin{proof}
To prove \ref{item:B1--},  we use 
the fact  that the set $K$ is forward contractive for $\Sigma$. Hence, the solutions to $\Sigma$ starting outside the set $K$ cannot reach $K$ at negative times. This implies that $ B_K(x) > 0$ for all $x \in \mathbb{R}^n \backslash K$.   Similarly,  the set $\mathbb{R}^n \backslash 
\mbox{int}(K)$ is forward contractive for 
$\Sigma^-$. Hence, the solutions to $\Sigma$ starting from $K$ cannot reach $\mathbb{R}^n \backslash \mbox{int}(K)$ at  positive times. This implies that $B_K(x) \leq 0$ for all $x \in K$. 

To prove \ref{item:B3--}, we start noting that 
\begin{align} \label{eqBBhat}
B_{K}(x) = \hat{B}_{K}(x)  \qquad \forall x \in U_1, 
\end{align}

which implies, using \ref{item:72--}, that $B_K$ is upper semicontinuous on the set $U_1 \backslash K$. Furthermore,  using \cite[Example 5.23]{rockafellar2009variational},  we conclude that $\text{Proj}_{U_1}$ is outer semicontinuous and locally bounded; hence, admits compact images. 
Thus,   using \cite[Theorem 1.4.16]{aubin2009set},  \ref{item:B3--} follows.   

Similarly,  to prove \ref{item:B4--}, we use \eqref{eqBBhat} and 
\ref{item:73-} to conclude that $B_K$ is lower semicontinuous on $U_1 \cap K$.    Furthermore, for each $x \in K \backslash U_1$,  we have 
$ B_K(x) := - \sup \{ -\hat{B}_K(y) : y \in \text{Proj}_{U_1}(x) \}.   $
Hence,  using \cite[Example 5.23]{rockafellar2009variational},  we know that $\text{Proj}_{U_1}$ is outer semicontinuous and locally bounded; thus,  admits compact images.  Thus,  using \cite[Theorem 1.4.16]{aubin2009set},  we conclude that $-B_K$ is upper semicontinuous; thus, $B_K$ is lower semicontinuous on $K$.

Finally, to prove \ref{item:B5--}, we let $\phi$ be a solution to $\Sigma$ (not necessarily maximal) starting from $x_o \in U_1 $ and such that $\phi(\dom \phi) \subset U_1$. Furthermore, we let
$t \in \dom \phi$ and $\sigma>0$ such that $[t,t+ \sigma] \subset \dom \phi$.  Hence,
$$ B_{K}(\phi(s)) =  \hat{B}_{K}(\phi(s)) \qquad \forall s \in [t,t+\sigma]. $$  
Next, we distinguish three situations. 
\begin{enumerate} [label={\arabic*)},leftmargin=*]
\item  When $\phi([t,t+\sigma]) \subset 
U_1 \backslash K$, there exists a solution $\psi \in \mathcal{S}_{\Sigma}(\phi(t+
\sigma))$ such that
$ B_K(\phi(t+\sigma))=T_K(\psi)$. Additionally, we introduce the solution 
$\hat{\psi} \in \mathcal{S}_{\Sigma}(\phi(t))$ satisfying
\begin{align*} 
\hat{\psi}(s) = 
\left\{
\begin{matrix}
\psi (s-\sigma) & \forall s \geq \sigma, \\
\phi(t+s) &  \forall s \in [0, \sigma]. \end{matrix} \right.
\end{align*}
Then, $T_{K}(\hat{\psi}) = \sigma + T_{K}(\psi)$. 
Hence, \eqref{eqhatB} implies that
$$ B_{K}(\phi(t)) \geq T_{K}(\hat{\psi}) = \sigma + B_{K}(\phi(t+\sigma)). $$

\item  When $\phi([t,t+\sigma]) \subset U_1 \cap \mbox{int}(K)$,
 there exists a solution 
 $\psi \in \mathcal{S}_{\Sigma}(\phi(t))$ 
such that $B_{K}(\phi(t))=T_{K}(\psi)$. 
Additionally, we introduce the solution $\hat{\psi} \in \mathcal{S}_{\Sigma}(\phi(t+
\sigma))$ satisfying
\begin{align*} 
\hat{\psi}(s) = 
\left\{ 
\begin{matrix}
\psi(s+\sigma) & \forall s \leq - \sigma, \\
 \phi(t+\sigma+s) & \forall s \in [- \sigma, 0].
\end{matrix} \right.
\end{align*}
As a result, we have $T_{K}(\hat{\psi}) = -\sigma + T_{K}(\psi)$ and \eqref{eqhatB} implies that
$$ B_{K}(\phi(t+\sigma)) 
\leq T_{K}(\hat{\psi}) = - \sigma + 
B_{K}(\phi(t)). $$

\item Finally, we consider a solution 
$\phi$ satisfying 
$ \phi(t) \in U_1  \backslash K$,  $\phi(t+\sigma) \in  U_1 \cap \mbox{int}(K)$,  
and such that there exists a unique $t_1 \in (t, t+\sigma)$ for which  
$ \phi(t_1) \in \partial K$.   Indeed, the aforementioned scenario complement the scenarios in the previous two items since, due to \ref{item:73--} and \ref{item:74--}, the solutions cannot slide on $\partial K$.  

Now, as in the previous steps, we conclude that 
$$ B_{K}(\phi(t+t'_1)) - B_{K}(\phi(t)) \leq -  t'_1 \qquad \forall t'_1 \in [0,t_1) $$
and
$$
B_{K}(\phi(t+\sigma)) - B_{K}(\phi(t + t'_1)) 
\leq  - (\sigma -  t'_1) \qquad \forall t'_1 \in (t_1,\sigma]. 
$$
Next, we note that
$$ B_{K}(\phi(t+t'_1)) > 0 = B_{K}(\phi(t+t_1))  \qquad \forall t'_1 \in [0,t_1) $$
and 
$$
 B_{K}(\phi(t+t'_1)) < 0 = 
 B_{K}(\phi(t+t_1)) \qquad \forall t'_1 \in (t_1,\sigma]. $$
Hence,  
$$ B_K(\phi(t+t_1)) - B_K(\phi(t)) \leq - t_1 \quad \text{and} 
\quad 
B_K(\phi(t+\sigma)) - B_K(\phi(t + t_1))  
\leq  - (\sigma -  t_1),
$$
 implying that  
$ B_K(\phi(t+\sigma)) - B_K(\phi(t))  
\leq  - \sigma. $  
\end{enumerate}
\end{proof}

\section{Proof of Theorem \ref{thm3}} 
\label{Sec.5}

\subsection{Proof steps} 

We propose to prove Theorem \ref{thm3}
by following 2 steps. 

\begin{itemize}
\item[\underline{Step 1}] We show the existence of a robustness margin $\bar{\epsilon}$ such that,  the set
\begin{align} \label{eqinfreachbis}
K_{\bar{\epsilon}} := \bigcup_{t \geq 0} \bigcup_{x \in X_o} R_{\Sigma_{\bar{\epsilon}}} (t, x)
\end{align}
satisfies $X_o \subset \cl(K_{\bar{\epsilon}})$  and  $\cl(K_{\bar{\epsilon}}) \cap X_u = \emptyset$. 
\item[\underline{Step 2}] We show the existence of a robustness margin $\epsilon_1$, which is smaller than $\bar{\epsilon}$,  such that
 the set $\cl(K_{\bar{\epsilon}})$ verifies Assumption \ref{assnew}  after replacing 
$(\Sigma, K)$ therein by 
$(\Sigma_{\epsilon_1}, \cl(K_{\bar{\epsilon}}))$.  
\end{itemize}
Finally,  using Lemma \ref{lem5-},  we conclude that \ref{item:C2bis} is verified for $\epsilon = \epsilon_1$,  $B = B_{\cl(K_{\bar{\epsilon}})}$,  and $B_{\cl(K_{\bar{\epsilon}})}$ as introduced in \eqref{eqhatB} and \eqref{eqbarcand} while replacing 
$(\Sigma, K)$ therein by 
$(\Sigma_{\epsilon_1}, \cl(K_{\bar{\epsilon}}))$.  

\subsection{Proof of Step 1 }  \label{Sec.Separ}

We will show that when $\Sigma$ is robustly safe with respect to $(X_o,X_u)$ and Assumption \ref{ass4-} holds,  then
\begin{align} \label{eqdisconnect-}
X_u \cap \cl(K_{\bar{\epsilon}}) = \emptyset, 
\end{align}
for appropriately chosen robustness margin $\bar{\epsilon}$.  Furthermore,  when additionally Assumption \ref{ass4} holds, we show that 
\begin{align} \label{eqinter}
\cl(X_u) \cap \cl(K_{\bar{\epsilon}}) = \emptyset.
\end{align}

\begin{lemma} [Separation between  $X_u$ and $\cl(K_{\bar{\epsilon}})$] \label{lem2bis}
Consider system $\Sigma$ such that Assumption \ref{ass3} holds. Furthermore, consider $(X_o,X_u) \subset \mathbb{R}^n \times \mathbb{R}^n$ such that $\Sigma$ is robustly safe with respect to $(X_o,X_u)$. Then, given a robustness margin $\bar{\epsilon}_o$, we conclude that, for each $\bar{\epsilon} \in \mathcal{C}_+$ satisfying 
\begin{align} \label{eqlosttrack}
\bar{\epsilon}(x) < \bar{\epsilon}_o(x) \qquad \forall x \in \mathbb{R}^n,
\end{align}
the set $K_{\bar{\epsilon}}$ introduced in \eqref{eqinfreachbis} satisfies the following properties. 
\begin{enumerate} [label={\arabic*)},leftmargin=*]
\item \label{item:P52b} If Assumption \ref{ass4} holds, then \eqref{eqinter} holds. 
\item \label{item:P52a} If Assumption \ref{ass4-} holds, then \eqref{eqdisconnect-} holds.
\end{enumerate} 
\end{lemma}

\begin{proof}
We prove \ref{item:P52b} using contradiction. That is, when \eqref{eqinter} is not satisfied, we conclude that there exists $x \in \partial X_u \cap \partial K_{\bar{\epsilon}}$ and such that $x \in \partial K_{\bar{\epsilon}_o}$. The latter is the only possibility since $\bar{\epsilon}_o$ is a robustness margin and since 
$K_{\bar{\epsilon}} \subset K_{\bar{\epsilon}_o}$.  Furthermore, using Assumption \ref{ass4}, we conclude that $x \notin \cl(X_o)$. Next, we consider a sequence $\{x_i\}^{\infty}_{i = 0} \subset K_{\bar{\epsilon}}$ 
that converges to $x$.  For 
$ \underline{\epsilon} := 
\min \{ \bar{\epsilon}_o(y)-\bar{\epsilon}(y) : 
y \in x + \mathbb{B} \}, $
we pick $\Delta \in (0,1)$ such that 
$ (x + \Delta \mathbb{B}) \cap \cl (X_o) = \emptyset $
and, for each 
$(x_1,x_2) \in (x + \Delta \mathbb{B}) \times (x + \Delta \mathbb{B})$ and for each $\eta_1 \in F(x_1)$, there exists $\eta_2 \in F(x_2)$ such that
\begin{align} \label{eqAdded}
|\eta_1 - \eta_2| \leq \frac{\underline{\epsilon}}{2}. 
\end{align}  
The latter is possible since the set-valued map $F$ is assumed to be continuous.
Without loss of generality, we assume that 
\begin{align} \label{eqstep1}
x_i \in x + \frac{\Delta}{4} \mathbb{B} \qquad \forall i \in \mathbb{N}. 
\end{align}
Furthermore, we consider a sequence of solutions 
$\{\phi_i\}^{\infty}_{i = 0}$ to $\Sigma_{\bar{\epsilon}}$ such that each solution $\phi_i$ starts from $x_{oi} \in \cl(X_o)$. Moreover, for each $i \in \mathbb{N}$, there exist $t^1_i > t^o_i > 0$ such that
$ \phi_i([t^o_i,t^1_i]) \subset x + \frac{\Delta}{2} \mathbb{B}$,   $\phi_i([t^o_i,t^1_i]) \subset K_{\bar{\epsilon}}, $ 
$ 
|\phi_i(t^o_i) - x| = \frac{\Delta}{2}$,  and
$\phi_i(t^1_i) = x_i$.   
Next, using \eqref{eqstep1} and local boundedness of the map $F + \bar{\epsilon} \mathbb{B}$, we conclude the existence of $t > 0$ such that 
\begin{align} \label{eqDelta+i}
t^1_i - t^o_i \geq t \qquad \forall i \in \mathbb{N}. 
\end{align}
Now, we let $\tilde{\Sigma} := \Sigma_{\bar{\epsilon}}$; which implies that  
$\tilde{\Sigma}_{\epsilon := \bar{\epsilon}_o - \bar{\epsilon}} = \Sigma_{\bar{\epsilon}_o}$. 
Furthermore, using a similar approach as in  the proof of Lemma \ref{lem1},  we  will show the existence of $\delta>0$ such that
\begin{align} \label{eqreachreach}
 \phi_i(t^1_i) + \delta \mathbb{B} \subset R^b_{\tilde{\Sigma}_{\epsilon}}(t^1_i - t^o_i,\phi_i(t^o_i)) \quad \forall i \in \mathbb{N}. 
\end{align}    
Indeed, let $\delta \in (0, \Delta/4]$ and note that 
$ x_i \in 
R^b_{\tilde{\Sigma}}(t^1_i - t^o_i,y_i) 
\backslash \{ y_i \}$,   where $y_i := \phi_i(t^o_i)$. 

Now, given $z \in x_i + \delta \mathbb{B}$, we consider the function 
$$ \eta_{zi}(s) := \phi_i(s+t^o_{i}) - \frac{s}{t^1_i-t^o_i} (x_i-z) \quad \forall s \in [0,t^1_i-t^o_i]. $$
Note that 
$ \eta_{zi}(s) \in x + \Delta \mathbb{B}$ for all 
$s \in [0,t^1_i-t^o_i]$, and for all $i \in \mathbb{N}$.   Furthermore, for almost all $s \in [0,t^1_i-t^o_i]$, we have 
$$ \dot{\eta}_{zi}(s) := \dot{\phi}_i(s + t^o_i) - \frac{1}{t^1_i-t^o_i} (x_i-z). $$
Next, for almost all 
$s \in [0, t^1_i-t^o_i]$, we let 
$\eta_{1i}(s) \in F(\eta_{zi}(s))$ be such that 
$|\eta_{1i}(s) - \dot{\phi}_i(s + t^o_i)| \leq \underline{\epsilon}/2$.  The latter is possible using \eqref{eqAdded}. Hence, for almost all $s \in [0,t^1_i-t^o_i]$, we have  
\begin{align*} 
\dot{\eta}_{zi}(s) & := \eta_{1i}(s) + (\dot{\phi}_i(s + t^o_i) - \eta_{1i}(s)) - \frac{1}{t^1_i-t^o_i} (x_i-z) \in F(\eta_{zi}(s)) + 
\left( \frac{\underline{\epsilon}}{2} + \frac{\delta}{t^1_i-t^o_i} \right) \mathbb{B}. 
\end{align*} 
Hence, by taking  $\delta := \min \{ \underline{\epsilon} t, \Delta \}/4,$
where $t$ comes from \eqref{eqDelta+i}, we conclude that, for almost all 
$s \in [0, t^1_i-t^o_i]$,
\begin{align*} 
\dot{\eta}_{zi}(s) & \in F(\eta_{zi}(s)) + 
\underline{\epsilon} \mathbb{B} \subset F(\eta_{zi}(s)) + 
\epsilon(\eta_{zi}(s)) \mathbb{B}. 
\end{align*}
Hence, $\eta_{zi} : 
[0,t^1_i-t^o_i] \rightarrow x + \Delta \mathbb{B}$ is a solution to 
$\tilde{\Sigma}_{\epsilon}$  with 
$\eta_{zi}(0) = y_i$ and $\eta_{zi}(t^1_i-t^o_i) = z$, which proves \eqref{eqreachreach}. 

Next, since $\delta$ is uniform in $i$, we conclude the existence of $i \in \mathbb{N}$ sufficiently large such that 
$$ x \in \mbox{int} \left( R^b_{\tilde{\Sigma}_{\epsilon}}(t^1_i - t^o_i,y_i) \right). $$    Hence, since $x \in \partial K_{\bar{\epsilon}_o}$, we conclude that there exists 
$y \notin \cl(K_{\bar{\epsilon}_o})$ such that $ y \in  R^b_{\tilde{\Sigma}_{\epsilon}}(t^1_i - t^o_i,y_i)$.  
However, $y_i = \phi_i(t^o_i) \in K_{\bar{\epsilon}} \subset K_{\bar{\epsilon}_o}$ and the latter set is, by definition, forward invariant for $\tilde{\Sigma}_{\epsilon} = \Sigma_{\bar{\epsilon}_o}$, which yields to a contradiction.  

Finally, we prove \ref{item:P52a}) using contradiction. That is, when \eqref{eqdisconnect-} is not satisfied, we conclude that there exists 
$x \in \partial X_u \cap 
\partial K_{\bar{\epsilon}} \cap X_u$ such that $x \in \partial K_{\bar{\epsilon}_o}$. The latter is the only possibility since $\bar{\epsilon}_o$ is a robustness margin and $K_{\bar{\epsilon}} \subset K_{\bar{\epsilon}_o}$ by definition. 
Furthermore, we distinguish two situations:

\begin{itemize}
\item When $x \notin \cl(X_o)$, the contradiction follows using the same steps as in the proof of \ref{item:P52b}).

\item When $x \in \cl(X_o)$, the contradiction follows using Assumption \ref{ass4-}.
\end{itemize}
\end{proof}

Clearly,  using Lemma \ref{lem2bis},  we conclude the existence of a robustness margin  $\bar{\epsilon}$ such that \eqref{eqdisconnect-} holds. Hence, Step 1 is verified. 

\subsection{Proof of Step 2} 
 
It is shown in \cite[Lemma 2]{ratschan2018converse} that,  given $\bar{\epsilon}: \mathbb{R}^n \rightarrow \mathbb{R}_{>0}$ continuous,  when  
$\cl(K_{\bar{\epsilon}})$ is bounded and $F$ is single valued and smooth,  every maximal solution to $\Sigma$  starting from $R_{\Sigma}(t,  \partial K_{\bar{\epsilon}})$,   for some $t \in \mathbb{R}$,   must cross $\partial K_{\bar{\epsilon}}$ only one time. Motivated by this fact, we establish contractivity of  
$\cl(K_{\bar{\epsilon}})$ for $\Sigma_\epsilon$, for an appropriate choice of the perturbation $\epsilon$.

\begin{lemma} [Contarctivity of the set 
$\cl(K_{\bar{\epsilon}})$]  \label{lem2}
Consider system $\Sigma$  such that Assumption \ref{ass3} holds.   
Consider $X_o \subset\mathbb{R}^n$,  a  $\bar{\epsilon}\in \mathcal{C}_+$,  and the set $K_{\bar{\epsilon}}$ introduced in \eqref{eqinfreachbis}.  Then,  for each  $\epsilon \in \mathcal{C}_+$ satisfying  
\begin{align} \label{eqepsineq}
\epsilon(x) < \bar{\epsilon}(x) \qquad \forall x \in \mathbb{R}^n,
\end{align}
the following properties hold:

\begin{enumerate} 
[label={P\ref{lem2}\arabic*)},leftmargin=*]
\item \label{item:P3a} The set 
$\cl(K_{\bar{\epsilon}})$ is forward contractive for $\Sigma_{\epsilon}$.
\item \label{item:P3b} The set $\mathbb{R}^n \backslash \mbox{int}(K_{\bar{\epsilon}})$ is forward contractive for $\Sigma^-_{\epsilon}$.
\item \label{item:P3c}  
The solutions 
to $\Sigma_\epsilon$ starting from $\text{int}(K_{\bar{\epsilon}})$ never reach 
$\partial K_{\bar{\epsilon}}$ for positive times.
\item \label{item:P3d} 
The solutions to $\Sigma_\epsilon$ starting from $\mathbb{R}^n \backslash \cl(K_{\bar{\epsilon}})$ never reach 
$\partial K_{\bar{\epsilon}}$ for negative times.
\end{enumerate}
\end{lemma}

\begin{proof}
We prove \ref{item:P3a} by first proving that the set $\cl(K_{\bar{\epsilon}})$ is forward invariant for $\Sigma_{\epsilon}$. To find a contradiction, we assume the existence of $x \in \partial K_{\bar{\epsilon}}$ and a solution $\phi \in \mathcal{S}_{\Sigma_{\epsilon}}(x)$ such that, for some $T_1 >0$,
$ \phi(s) \in \mathbb{R}^n \backslash \cl(K_{\bar{\epsilon}})$  for all 
$s \in (0,T_1]$.   Next, we take $t \in (0, \min\{T_1,T\}]$, where $T>0$ comes from Lemma \ref{lem1}. Hence, using Lemma \ref{lem1}, we conclude the existence of 
$\delta>0$ such that
$ x + \delta \mathbb{B} \subset 
 R^b_{\Sigma_{\bar{\epsilon}}}(-t, \phi(t)). $
Hence, there exists 
$y \in \mbox{int} (K_{\bar{\epsilon}})$ and a solution $\psi$ to $\Sigma_{\bar{\epsilon}}$ such that $\psi(0) = y$ and  $\psi(t) = \phi(t)$. 
The latter implies that the set 
$K_{\bar{\epsilon}}$ is forward invariant for 
$\Sigma_{\bar{\epsilon}}$, which yields to a contradiction since,  by its definition, the set $K_{\bar{\epsilon}}$ is forward invariant for 
$\Sigma_{\bar{\epsilon}}$. 

To complete the proof, we show that the solutions to $\Sigma_{\epsilon}$ cannot slide on $\partial K_{\bar{\epsilon}}$; namely, for each $x \in \partial K_{\bar{\epsilon}}$ and for each 
$\phi \in \mathcal{S}_{\Sigma_{\epsilon}}(x)$, there exists $T>0$ such that 
\begin{align} \label{eqcontractive}
\phi(t) \in \mbox{int} (K_{\bar{\epsilon}}) \qquad \forall t \in (0,T].  
\end{align}
To find a contradiction,  we assume the existence of a solution $\phi$ to $\Sigma_{\epsilon}$ starting from $x_o \in \mathbb{R}^n$ and an interval $(t_1, t_2) \subset \dom \phi$,  with $t_2 > t_1$,  such that 
$ \phi(s) \in \partial K_{\bar{\epsilon}}$ for all $s \in (t_1, t_2)$. 
Next, we pick $t_3 \in (t_1, t_2)$ and we let
$ x := \phi(t_3) \in \partial K_{\bar{\epsilon}}$. 
We take $t \in (0,\min\{(t_2 - t_3),T\})$, where $T>0$ comes from Lemma \ref{lem1}.  Furthermore, we let $y  := \phi( t + t_3) \in \partial K_{\bar{\epsilon}}$.    Hence,  using Lemma \ref{lem1},  we conclude that
\begin{align} \label{eqball+}
y + \delta \mathbb{B} \subset R^b_{\Sigma_{\bar{\epsilon}}}(t,\phi(t_3)). 
\end{align} 
Similarly, we take $t \in (0,\min\{(t_3 - t_1),T\})$, where $T>0$ comes from Lemma \ref{lem1} while using $(\Sigma^{-}_{\epsilon}, \Sigma^-_{\bar{\epsilon}})$ instead of $(\Sigma, \Sigma_{\bar{\epsilon}})$. Furthermore, we let $ y  := \phi( -t + t_3 ) \in \partial K_{\bar{\epsilon}} $.     Hence, using Lemma \ref{lem1}, we conclude that
\begin{align} \label{eqball-}
y + \delta \mathbb{B} \subset 
R^b_{\Sigma^-_{\bar{\epsilon}}}(t,\phi(t_3)). 
\end{align} 
 Finally, combining \eqref{eqball+} and \eqref{eqball-}, we conclude the existence of a solution to $\Sigma_{\bar{\epsilon}}$ starting from $\mbox{int} (K_{\bar{\epsilon}})$ and that leaves the set $K_{\bar{\epsilon}}$. This latter fact yields to a contradiction since, by definition, the set $K_{\bar{\epsilon}}$ is forward invariant for $\Sigma_{\bar{\epsilon}}$.

To prove \ref{item:P3b}, we recall that the solutions to $\Sigma_{\epsilon}$ cannot slide on $\partial K_{\bar{\epsilon}}$ along positive time intervals. Hence, the same property must hold for the solutions to $\Sigma^{-}_{\epsilon}$. As a result, if \ref{item:P3b} is not satisfied, then there exists a maximal solution $\phi \in \mathcal{S}_{\Sigma^{-}_{\epsilon}}(x_o)$, for some $x_o \in \partial K_{\bar{\epsilon}}$, and $T_o>0$ such that
$\phi((0,T_o]) \subset \mbox{int} (K_{\bar{\epsilon}})$.
Using Lemma \ref{lem1},  we conclude the existence of $T \in (0, T_o]$ and $\delta > 0$ such that,  for $y := \phi(T) \in K_{\bar{\epsilon}}$,  we have 
$$ x_o + \delta \mathbb{B}  \subset 
 R^b_{\Sigma^-_{\bar{\epsilon}}}(-T,y) = R^b_{\Sigma_{\bar{\epsilon}}}(T,y).  $$   
 However, the latter contradicts forward invariance of the set $K_{\bar{\epsilon}}$ for $\Sigma_{\bar{\epsilon}}$.

To prove \ref{item:P3c} using contradiction,   we assume the existence of a solution $\phi$  to $\Sigma_\epsilon$ starting from 
$x \in \text{int}(K_{\bar{\epsilon}})$  such that, for some $t_1>0$, we have 
 $ \phi(t_1) \in 
\partial K_{\bar{\epsilon}}$ and $\phi([0,t_1)) \in 
\text{int} (K_{\bar{\epsilon}}). $
Next, we take 
$ y := \phi(t_1) \in R^b_{\Sigma_{\epsilon}}(t_1) \backslash 
\{x\} $
and 
$$ z := \phi(t_1-t) \in R^b_{\Sigma_{\epsilon}}
(-t,y) \backslash \{y\} = R^b_{\Sigma^{-}_{\epsilon}}(t,y) \backslash \{y\}, $$
for some $t \in (0,t_1]$ to be determined. 
Hence, $z \in \text{int}(K_{\bar{\epsilon}})$. Next, using Lemma \ref{lem1} while 
replacing $(x,y,\Sigma,\Sigma_{\epsilon})$ by 
$(y,z,\Sigma^{-}_{\epsilon}, 
\Sigma^{-}_{(\epsilon + \bar{\epsilon})/2})$, we conclude the existence of $T>0$ such that, 
for $t := \min\{T,t_1\}$, there exists 
$\delta>0$ such that 
$  y + \delta \mathbb{B} \subset 
R^b_{\Sigma_{(\bar{\epsilon} + \epsilon)/2}}(t,z). $
However, according to \ref{item:P3a}, the 
set $\cl (K_{\bar{\epsilon}})$ must be forward invariant for $\Sigma_{(\bar{\epsilon} + \epsilon)/2}$, since 
\begin{align} \label{eqeps}
0 < (\bar{\epsilon}(x) + \epsilon(x))/2 < \bar{\epsilon}(x) \qquad \forall x \in \mathbb{R}^n, 
\end{align}
which yields to a contradiction.

Similarly, to prove \ref{item:P3d} using contradiction, we assume the existence of a solution $\phi$ to $\Sigma_\epsilon$ starting from $x \in \mathbb{R}^n \backslash \cl(K_{\bar{\epsilon}})$ 
such that, for some $t_1>0$, we have 
 $ \phi(-t_1) \in \partial K_{\bar{\epsilon}}$ and $\phi([-t_1,0)) \in \text{int} (K_{\bar{\epsilon}})$.   
Next, we take 
$$ y := \phi(-t_1) \in R^b_{\Sigma_{\epsilon}}(-t_1,x) \backslash 
\{x\} $$
and 
$$ z := \phi(t-t_1) \in R^b_{\Sigma_{\epsilon}}
(t,y) \backslash \{y\} = R^b_{\Sigma^{-}_{\epsilon}}(-t,y) \backslash \{y\}, $$
for some $t \in (0,t_1]$ to be 
determined. Hence, $z \in \mathbb{R}^n \backslash \cl(K_{\bar{\epsilon}})$. Next, using Lemma \ref{lem1} while 
replacing $(x,y,\Sigma,\Sigma_{\epsilon})$ by 
$(y,z,\Sigma_{\epsilon}, 
\Sigma_{(\epsilon + \bar{\epsilon})/2})$, we conclude the existence of $T>0$ such that, 
for $t := \min\{T,t_1\}$, there exists 
$\delta>0$ such that 
$$  y + \delta \mathbb{B} \subset 
R^b_{\Sigma^{-}_{(\bar{\epsilon} + \epsilon)/2}}(t,z) = R^b_{\Sigma_{(\bar{\epsilon} + \epsilon)/2}}(-t,z). $$
However, according to \ref{item:P3b}, the 
set $\cl (K_{\bar{\epsilon}})$ must be forward invariant for $\Sigma_{(\bar{\epsilon} + \epsilon)/2}$,  since  \eqref{eqeps} holds, 
which yields to a contradiction.
\end{proof}

It is shown in  \cite[Lemma 3]{ratschan2018converse} that,  when the set $\cl(K_{\bar{\epsilon}})$ is bounded and when $F$ is single-valued and smooth,  there exists $U \subset \mathbb{R}^n$ 
a neighborhood of $\partial K_{\bar{\epsilon}}$ such that every maximal solution to $\Sigma$ starting from $U$ reaches $\partial K_{\bar{\epsilon}}$ in finite (positive or negative) time.   
Inspired by this observation and using Lemma \ref{lem2},  we establish recurrence properties of the set $K_{\bar{\epsilon}}$ on a neighborhood of $\partial K_{\bar{\epsilon}}$.

\begin{lemma} [Recurrence of the set 
$\cl(K_{\bar{\epsilon}})$] \label{lem3}
Consider system $\Sigma$ such that Assumption  \ref{ass3} holds.  
Consider $X_o \subset \mathbb{R}^n$ and  $\epsilon_1, \bar{\epsilon} \in \mathcal{C}_+$ such that 
\begin{align} \label{eqepsilon1}
\epsilon_1(x) < \bar{\epsilon}(x) \qquad \forall x \in \mathbb{R}^n.
\end{align}
Then, there exists a closed subset $U_{1} \subset \mathbb{R}^n$ such that 
\begin{align} \label{eqinter0}
\partial K_{\bar{\epsilon}} \subset \text{int} (U_1),
\end{align}
where the set $K_{\bar{\epsilon}}$ is introduced in \eqref{eqinfreachbis}.  Moreover,   for each 
$\epsilon \in \mathcal{C}_+$ satisfying  
\begin{align} \label{eqepsilon2}
\epsilon(x) \leq \epsilon_1(x) \qquad \forall x \in \mathbb{R}^n, 
\end{align}
the following properties hold:

\begin{enumerate} [label={P\ref{lem3}\arabic*)},leftmargin=*]
\item \label{item:R1} The set $K_{\bar{\epsilon}}$ is locally recurrent for $\Sigma_{\epsilon}$ on $U_1$.
\item \label{item:R2} The set 
$\mathbb{R}^n\backslash K_{\bar{\epsilon}}$ is locally recurrent for $\Sigma^-_{\epsilon}$ on $U_1$.
\end{enumerate}
\end{lemma}

\begin{proof}
To prove \ref{item:R1},  we will show that,  for each $x_o \in \partial K_{\bar{\epsilon}}$,  there exists $\delta > 0$ such that,  for each $x \in (x_o + \delta \mathbb{B}) \backslash \cl(K_{\bar{\epsilon}})$ and for each solution $\phi \in \mathcal{S}_{\Sigma_{\epsilon}}(x)$, there exists $T_\phi > 0$ such that $ \phi(T_\phi) \in \mbox{int} (K_{\bar{\epsilon}}) $.  First,  using Lemma \ref{lem2} and  according to \ref{item:P3a} and \ref{item:P3c},  
we conclude that $ R^b_{\Sigma_\epsilon}(T,x_o) \subset \mbox{int}(K_{\bar{\epsilon}})$ for all $T>0$.  Furthermore,  since $F$ is locally bounded,  using Lemmas \ref{lemprepre} and \ref{lem1-},  we conclude that,  for each $T>0$ small,  the set $R^b_{\Sigma_\epsilon}(T,U)$ is compact,  for $U$ a sufficiently small neighborhood of $x_o$.  

Hence, there exists $\alpha>0$ such that 
$ \min \{|y|_{\partial K_{\bar{\epsilon}}} : y \in  R^b_{\Sigma_\epsilon}(T,x_o)\}  \geq \alpha.  $   
Now,  using the second item in Lemma \ref{lem1-} in the Appendix,  we conclude that there exists $\delta>0$ such that, for each $x \in (x_o + \delta \mathbb{B}) \backslash \cl(K_{\bar{\epsilon}})$, we have $R^b_{\Sigma_\epsilon}(T, x) \subset R^b_{\Sigma_\epsilon}(T, x_o) + \alpha/2 \mathbb{B}$. 
Hence, $ R^b_{\Sigma_\epsilon}(T, x) \subset \mbox{int}(K_{\bar{\epsilon}})$, which proves \ref{item:R1}. Finally, \ref{item:R2} can be proved following the exact steps, while using \ref{item:P3b} and \ref{item:P3d} instead of \ref{item:P3a} and \ref{item:P3c}.  
\end{proof}

At this point,  to prove  Step 2,  we let $U_1$ and $\epsilon_1$ be,  respectively,  the set and the robustness margin introduced in Lemma \ref{lem3}.  Furthermore,  we show that Assumption \ref{assnew} is verified  after replacing 
$(\Sigma, K)$ therein by 
$(\Sigma_{\epsilon_1}, \cl(K_{\bar{\epsilon}}))$.  Indeed,  we  use \ref{item:R1} to conclude that \ref{item:71--} holds.  Furthermore,  we use \ref{item:R2} to conclude that \ref{item:72--} holds.  We use \ref{item:P3a} to conclude that \ref{item:73--} holds.  Finally,  we use \ref{item:P3b} to conclude that \ref{item:74--} holds.

\section{Proof of Theorem \ref{thm4}} \label{Sec.6}

\subsection{Proof steps}
 
Given  $\delta, ~\bar{\epsilon}  \in \mathcal{C}_+$, we introduce the inflation of the set $K_{\bar{\epsilon}}$ given by
\begin{align} \label{Kinter1}
K_{\bar{\epsilon},\delta} := \bigcup_{x \in K_{\bar{\epsilon}}} (x + \delta(x) \mathbb{B}). 
\end{align}

Since we fail to show that the set $\cl(K_{\bar{\epsilon}, \delta})$ is contractive,  even when the inflation $\delta$ is arbitrarily small,  we next propose an alternative set, denoted by $K_{\bar{\epsilon},\rho_o,\epsilon_1}$, and given by
\begin{equation} 
\label{Kinter1+}
\begin{aligned} 
K_{\bar{\epsilon},\rho_o,\epsilon_1}   := \bigcup_{t \geq 0} \bigcup_{x \in K_{\bar{\epsilon},\rho_o}} R_{\Sigma_{\epsilon_1}} (t, x),
\quad 
K_{\bar{\epsilon},\rho_o}  := \bigcup_{x \in K_{\bar{\epsilon}}} (x + \rho_o(x) \mathbb{B}), 
\end{aligned}
\end{equation}
where 
$\rho_o, \epsilon_1 : \mathbb{R}^n \rightarrow \mathbb{R}_{>0}$
are continuous functions.

The proof of Theorem \ref{thm4} follows in five steps.

\begin{itemize}
\item [\underline{Step 1}] We show the existence of  $\delta \in \mathcal{C}_+$, a subset $U_1 \subset \mathbb{R}^n$,  
and  robustness margins $\bar{\epsilon}$ and $\epsilon_1$  such that
\eqref{eqepsilon2} holds,
\begin{align} \label{eqpropadd}
\cl(K_{\bar{\epsilon}, \delta}) \backslash \text{int} (K_{\bar{\epsilon}}) \subset \text{int} (U_1),  \qquad  \cl(K_{\bar{\epsilon}, \delta}) \cap \cl(X_u) = \emptyset,  
\end{align}
and, at the same time,  $K_{\bar{\epsilon},\delta}$ enjoys the same recurrence properties established for $K_{\bar{\epsilon}}$ in Lemma \ref{lem3}.   

\item[\underline{Step 2}] We show the existence of $\rho_o \in \mathcal{C}_+$ such that $\cl(K_{\bar{\epsilon}, \rho_o, \epsilon_1})$ enjoys the same contractivity and recurrence properties established for the set $\cl(K_{\bar{\epsilon}})$ in Lemmas \ref{lem2} and \ref{lem3},  respectively.   Moreover, we carefully choose $\rho_o$ such that 
\begin{align} 
  \cl (K_{\bar{\epsilon}, \rho_o, \epsilon_1}) \subset \text{int}( K_{\bar{\epsilon}, \delta}),  \qquad  K_{\bar{\epsilon}} \subset  \text{int}(K_{\bar{\epsilon}, \rho_o, \epsilon_1}),  \label{eqpropadd1}
\\
\bigcup_{x \in \mathbb{R}^n \backslash \text{int}(K_{\bar{\epsilon},\delta})}  \hspace{-0.4cm} (x + \rho_o(x) \mathbb{B})  \cap \cl(K_{\bar{\epsilon},\rho_o,\epsilon_1}) = \emptyset,  \label{eqpropadd2}
\end{align}
and such that there exists $\hat{U}_1 \subset \mathbb{R}^n$ closed and verifying
\begin{align*}     
\cl(K_{\bar{\epsilon}, \delta}) \backslash \text{int} (K_{\bar{\epsilon}}) & \subset \text{int} (\hat{U}_1), \label{eqlem102}  \\
 \bigcup_{x \in \hat{U}_1}  (x + \rho_o(x) \mathbb{B}) & \subset U_1. \label{eqlem101} 
\end{align*}

\item[\underline{Step 3}] We show the existence of a robustness margin  $\epsilon_2$ such that \ref{item:C2bis} holds for $\epsilon = \epsilon_2$ and $B := B_{\cl(K_{\bar{\epsilon},\rho_o,\epsilon_1})}$,  where  $B_{\cl(K_{\bar{\epsilon},\rho_o,\epsilon_1})}$ is defined as in \eqref{eqbarcand} while replacing  $(K, \Sigma)$ therein by $(\cl (K_{\bar{\epsilon},\rho_o,\epsilon_1}),\Sigma_{\epsilon_2})$. We prove the latter  using Lemma \ref{lem5-},  \eqref{eqpropadd},  and \eqref{eqpropadd1},  under the recurrence and contractivity properties of the set $\cl(K_{\bar{\epsilon},\rho_o,\epsilon_1})$.  

\item[\underline{Step 4}] We show the existence of a smooth function $\rho_2 : \mathbb{R}^n \rightarrow \mathbb{R}_{>0}$ such that,  for each  $v \in \mathbb{B}$,   \ref{item:C2bis} holds for $\epsilon = 0$ and $B = B^v_{\cl(K_{\bar{\epsilon},\rho_o,\epsilon_1})}  : \mathbb{R}^n \rightarrow \mathbb{R}$,  where 
\begin{align} \label{eq.hatB}
B^{v}_{\cl(K_{\bar{\epsilon},
\rho_o,\epsilon_1})}(x) :=  B_{\cl(K_{\bar{\epsilon},
\rho_o,\epsilon_1})}(x + \rho_2(x) v).
\end{align}
We prove the latter using Lemma \ref{lem8} in the Appendix and \eqref{eqpropadd2}-\eqref{eqlem101}.

\item[\underline{Step 5}] Finally, thanks to Lemma \ref{lemsmoothing} in the Appendix,  we show that \ref{item:C3} holds for
\begin{equation}
\label{eqBcertificate} 
\begin{aligned} 
B(x) & := 
\int_{\mathbb{R}^n}  B_{\cl(K_{\bar{\epsilon},
\rho_o,\epsilon_1})}(x + \rho_2(x) v) \Psi(v) dv,  
\end{aligned}
\end{equation}
where $\Psi : \mathbb{R}^n \rightarrow [0,1]$ is smooth and satisfies 
\begin{equation}
\label{eqPsiprr} 
\begin{aligned} 
\Psi(v) = 0  \quad \forall v \in \mathbb{R}^n \backslash \mathbb{B} \qquad \text{and} \qquad
\int_{\mathbb{R}^n}  \Psi(v) dv = 1, 
\end{aligned}
\end{equation}
\end{itemize}

\subsection{ Proof of Step 1 }

The next statement is similar to \cite[Theorem 2]{SUBBARAMAN201654}. Although formulated for general hybrid systems in the aforementioned reference,  the set to render recurrent was assumed to be bounded.    

\begin{lemma}[Recurrence of the set $K_{\bar{\epsilon},\delta}$] \label{lem6}
Consider system $\Sigma$ such that Assumption \ref{ass3} holds. Consider two subsets $(X_o,X_u) \in \mathbb{R}^n \times \mathbb{R}^n$ such that $\cl(X_o) \cap \cl(X_u) = \emptyset$ and 
$\bar{\epsilon}, \epsilon_1 : \in \mathcal{C}_+$  such that 
\begin{align} \label{eqepsilon1n}
\epsilon_1(x) < \bar{\epsilon}(x) \qquad \forall x \in \mathbb{R}^n.
\end{align}
Consider the set $K_{\bar{\epsilon}}$ introduced in \eqref{eqinfreachbis} and a closed set $U_1 \subset \mathbb{R}^n$ such that the conclusions of Lemmas \ref{lem2} and \ref{lem3} hold. 
Then, there exists   
$\delta \in \mathcal{C}_+$  such that the set $K_{\bar{\epsilon},\delta}$ in \eqref{Kinter1} satisfies \eqref{eqpropadd} and,  for each  
$\epsilon \in \mathcal{C}_+$  satisfying 
\begin{align} \label{eqepsilon2n}
\epsilon(x) \leq \epsilon_1(x) \qquad \forall x \in \mathbb{R}^n, 
\end{align}
the following properties hold. 
\begin{enumerate} [label={P\ref{lem6}\arabic*)},leftmargin=*]
\item \label{item:I2} The set $K_{\bar{\epsilon},\delta}$ is locally recurrent for $\Sigma_{\epsilon}$ on $U_1$.
\item \label{item:I3} The set $\mathbb{R}^n \backslash K_{\bar{\epsilon},\delta}$ is locally recurrent for 
$\Sigma^{-}_{\epsilon}$ on $U_1$.
\end{enumerate}
\end{lemma}

\begin{proof}
\blue{We start using 
Lemma \ref{lem2bis} to conclude that $\cl(X_u) \cap \cl(K_{\bar{\epsilon}}) = \emptyset$. 
We next use 
Lemma \ref{lem3} to conclude the existence of a closed subset $U_1 \subset \mathbb{R}^n$ such that \eqref{eqinter0} holds and \ref{item:R1}-\ref{item:R2} hold 
for each  
$\epsilon \in \mathcal{C}_+$ 
 satisfying \eqref{eqepsilon2n}.
As a result, to prove the existence of  $\delta \in \mathcal{C}_+$ for which \eqref{eqpropadd} holds,  we combine \eqref{eqinter0} and the fact that $\cl(X_u) \cap \cl(K_{\bar{\epsilon}}) = \emptyset$ to conclude that \eqref{eqpropadd}  is verified for
$$ \delta(x) := \frac{1}{2} \min \left\{ |x|_{\mathbb{R}^n \backslash U_1},  |x|_{\text{cl}(X_u)} \right\}.  $$

To prove \ref{item:I2}, we use \ref{item:R1}, which shows that the set $K_{\bar{\epsilon}}$ is locally recurrent for $\Sigma_{\epsilon}$ on $U_1$. Hence, \ref{item:I2} follows since, by definition, 
$K_{\bar{\epsilon}} \subset K_{\bar{\epsilon},\delta}$. 

To prove \ref{item:I3}, we grid the set 
$\partial K_{\bar{\epsilon}}$ using a sequence of nonempty compact subsets 
$\{\partial K_i\}^{N}_{i=1}$, where $N \in \{1,2,...,\infty\}$. That is, we assume that 
 $\bigcup^{\infty}_{i = 1} \partial K_i = \partial K_{\bar{\epsilon}}$.   Furthermore, for each $i \in \{1,2,...,N\}$, there exists a finite set  $\mathcal{N}_i \subset \{1,2,...,N \}$ and $\beta_i>0$ such that  
 \begin{align} \label{eqhelp}
 (\partial K_i + \beta_i \mathbb{B}) \cap \partial K_j = \emptyset \quad  
 \forall j \notin \mathcal{N}_i  \quad  \text{and} \quad  
\mbox{int}_{\partial K_{\bar{\epsilon}}} (\partial K_i \cap \partial K_j) = \emptyset \quad  
\forall j \in \mathcal{N}_i, 
\end{align}
where $\mbox{int}_{\partial K_{\bar{\epsilon}}}(\cdot)$ denotes the  interior of $(\cdot)$ relative $\partial K_{\bar{\epsilon}}$.   Such a decomposition always exists according to the Whitney Covering Lemma \cite{10.2307/1989708}.

Next, we make the following claim that we prove later. 
\begin{claim} \label{clm1}
For each $i \in 
\{1,2,...,N\}$,  we can find $t_i>0$ and $\alpha_i>0$ such that
\begin{align} \label{eqtoshow}
 \bigcup_{i \in 
 \{1,2,...,N\}} R^b_{\Sigma^{-}_\epsilon}(t_i, \partial K_i)        
\bigcap  
\bigcup_{i \in \{1,2,...,N\}} \bigcup_{x \in \partial K_i} (x + \alpha_i \mathbb{B})
= \emptyset. 
\end{align}
\end{claim}
Under Claim \ref{clm1}, we introduce the function $\alpha : \partial K_{\bar{\epsilon}} \rightarrow \mathbb{R}_{>0}$ given by 
$$ \alpha(x) := \min \{\alpha_i/2 : 
i \in \{1,2,...,N\} ~ \text{s.t.} ~ x \in \partial K_i \}. $$
The function $\alpha$ is lower semicontinuous and allows, in view of \eqref{eqtoshow}, to conclude that
$$ \bigcup_{i \in \{1,2,...,N\}} R^b_{\Sigma^{-}_\epsilon}(t_i, \partial K_i)        
\bigcap  
 \bigcup_{x \in \partial K_{\bar{\epsilon}}} (x + \alpha(x) \mathbb{B})
= \emptyset. $$
Hence,  using   \cite[Theorem 1]{katvetov1951real},  we conclude the existence of a continuous function $\bar{\alpha} : \partial K_{\bar{\epsilon}} \rightarrow \mathbb{R}_{>0}$ 
such that
$$  \bar{\alpha}(x) \leq \alpha(x) \qquad  
\forall x \in \partial K_{\bar{\epsilon}}. $$  Hence,  
\begin{align} \label{eqhelp2}  
\bigcup_{i \in \{1,2,...,N\}} R^b_{\Sigma^{-}_\epsilon}(t_i, \partial K_i)        
\bigcap  
 \bigcup_{x \in \partial K_{\bar{\epsilon}}} (x + \bar{\alpha}(x) \mathbb{B})
= \emptyset. 
\end{align}
That is, every solution to $\Sigma^{-}_{\epsilon}$ starting from $\partial K_{\bar{\epsilon}}$ leaves the set
 $K_{\bar{\epsilon}} \cup \bigcup_{x \in \partial K_{\bar{\epsilon}}} (x + \bar{\alpha}(x) \mathbb{B})$.

Next, we show that, for each $i \in 
\{1,2,...,N\}$, we can find $\sigma_i > 0$ such that 
 $$  \bigcup_{i \in \{1,2,...,N\}} R^b_{\Sigma^{-}_\epsilon}(t_i, \partial K_i + \sigma_i \mathbb{B})        \bigcap  
 \bigcup_{x \in \partial K_{\bar{\epsilon}}} (x + \bar{\alpha}(x) \mathbb{B})
= \emptyset. $$
Indeed, the latter follows from a direct combination of Lemmas \ref{lem1-} and \ref{lemprepre}
 in the Appendix and \eqref{eqhelp2}. 

As a result, the function $\sigma : \partial K_{\bar{\epsilon}} \rightarrow \mathbb{R}_{>0}$ given by 
$$ \sigma(x) := \min \{\sigma_i : 
i \in \{1,2,...,N\} ~ \text{s.t.} ~ x \in \partial K_i \} $$
 is lower semicontinuous and allows, using   \cite[Theorem 1]{katvetov1951real},  to conclude the existence of a continuous function $\bar{\sigma} : \partial K_{\bar{\epsilon}} \rightarrow \mathbb{R}_{>0}$ 
such that
$$  \bar{\sigma}(x) \leq \sigma(x) \qquad  
\forall x \in \partial K_{\bar{\epsilon}}, $$ 
and, at the same time, the solutions  to $\Sigma^-_{\epsilon}$ starting from the set
$ \bigcup_{x \in \partial K_{\bar{\epsilon}}} (x + \bar{\sigma}(x) \mathbb{B}) \subset $
leave the set 
$ K_{\bar{\epsilon}} \cup \bigcup_{x \in \partial K_{\bar{\epsilon}}} (x + \bar{\alpha}(x) \mathbb{B}). $
Finally, after continuously extending 
$\bar{\sigma}$ to $\mathbb{R}^n$, we take
$$ \delta(x) := 
\left\{
\begin{matrix}
\min \{ \bar{\sigma}(x), |x|_{\partial (K_{\bar{\epsilon}} \cup K_1)} \} & \text{if}~ x \in \text{int} (K_{\bar{\epsilon}})
\\
\bar{\sigma}(x) & \text{otherwise},
\end{matrix}
\right. \qquad K_1 := \bigcup_{x \in \partial K_{\bar{\epsilon}}} (x + \bar{\sigma}(x) \mathbb{B}). $$
Using \ref{item:R2}, we conclude that the solutions to $\Sigma^{-}_{\epsilon}$ starting from $ K_{\bar{\epsilon}, \delta} \cap U_1$ 
leave the set $K_{\bar{\epsilon}, \delta}$. Hence, \ref{item:I3} is verified. 
 
To prove Claim \ref{clm1}, we use Lemma \ref{lem2} to conclude that, for each $i \in \{1,2,...,N\}$ and for each $t_i > 0$, we have 
$$ R^b_{\Sigma^{-}_\epsilon}(t_i, \partial K_i) \cap \cl(K_{\bar{\epsilon}}) = \emptyset. $$
Furthermore, using Lemma \ref{lemprepre} in the Appendix, we conclude that, for each $i \in \{1,2,...,N\}$, we can make  $t_i > 0$ small enough such that
$R^b_{\Sigma^{-}_\epsilon}(t_i, \partial K_i)$ is compact. Furthermore, in view of \eqref{eqhelp} and Lemma \ref{lem1-}, we can make each $t_i$ even smaller such that
\begin{align} \label{eqRconst}
(\partial K_j + \beta_j \mathbb{B}) \cap  R_{\Sigma^{-}_\epsilon}(t_i, \partial K_i) = \emptyset \qquad \forall j \notin \mathcal{N}_i.
\end{align}
As a result, for each $i \in \{1,2,...,N\}$, we can find 
$\alpha_i \in (0,\beta_i)$ such that
\begin{align} \label{eqRconstbis} 
R^b_{\Sigma^{-}_\epsilon}(t_i, \partial K_i)  \bigcap  \bigcup_{j \in 
\{ i,\mathcal{N}_i \}} \left( \partial K_j + \alpha_j \mathbb{B} \right) = \emptyset. 
\end{align}
The claim is proved by combining \eqref{eqRconst} and \eqref{eqRconstbis}. } 
\end{proof}

\subsection{ Proof of Step 2} 

We show that the set  introduced in \eqref{Kinter1+} preserves both recurrence and contractivity.   

\begin{lemma} [Recurrence and contractivity of the set $\cl(K_{\bar{\epsilon},\rho_o,\epsilon_1})$] \label{lem6+}
Consider system $\Sigma$ such that Assumption  \ref{ass3} holds.   Consider two subsets $(X_o,X_u) \subset \mathbb{R}^n \times \mathbb{R}^n$ and 
$\bar{\epsilon}, \epsilon_1, \epsilon_2\in \mathcal{C}_+$ such that 
\begin{align} \label{eqepsilon3}
\epsilon_2 (x) < \epsilon_1(x) < \bar{\epsilon}(x) \qquad \quad 
\forall x \in \mathbb{R}^n.
\end{align}
Consider the set $K_{\bar{\epsilon}}$ introduced in \eqref{eqinfreachbis},
 a closed subset $U_1 \subset \mathbb{R}^n$, and  $\delta \in \mathcal{C}_+$ such that the conclusions in Lemmas \ref{lem2}, \ref{lem3} and \ref{lem6} hold.
 
Then, there exists   
$\rho_o \in \mathcal{C}_+$
 such that the set $K_{\bar{\epsilon},\rho_o,\epsilon_1}$ introduced in \eqref{Kinter1+} satisfies 
\begin{align}     
\cl(K_{\bar{\epsilon}, \delta}) \backslash \text{int} (K_{\bar{\epsilon}}) & \subset \text{int} (\hat{U}_1), \label{eqlem102}  \\
 \bigcup_{x \in \hat{U}_1}  (x + \rho_o(x) \mathbb{B}) & \subset U_1, \label{eqlem101} 
\end{align}
and the following properties hold.  

\begin{enumerate} [label={P\ref{lem6+}\arabic*)},leftmargin=*]
\item \label{item:I8+} There exists a closed subset $\hat{U}_1 \subset \text{int}(U_1)$ such that \eqref{eqlem101} and \eqref{eqlem102} hold.

\item  For each continuous function
$\epsilon : \mathbb{R}^n \rightarrow \mathbb{R}_{\geq 0}$ satisfying
\begin{align} \label{eqepsilon4}
\epsilon (x) \leq  \epsilon_2(x) \qquad \forall x \in \mathbb{R}^n,
\end{align}
\begin{enumerate}
\item \label{item:I2+a}   $K_{\bar{\epsilon},\rho_o,\epsilon_1}$ is locally recurrent for $\Sigma_{\epsilon}$ on $U_1$.
\item \label{item:I2+b} $\mathbb{R}^n \backslash K_{\bar{\epsilon},\rho_o,\epsilon_1}$ is locally recurrent for 
$\Sigma^{-}_{\epsilon}$ on $U_1$.
\item \label{item:I2+c} 
$\cl(K_{\bar{\epsilon},\rho_o,\epsilon_1})$ is forward contractive for $\Sigma_{\epsilon}$.

\item \label{item:I2+d}  $\mathbb{R}^n \backslash \mbox{int}(K_{\bar{\epsilon},\rho_o,\epsilon_1})$ is forward contractive for $\Sigma^-_{\epsilon}$.
\end{enumerate}
\end{enumerate}
\end{lemma}

\begin{proof}
To prove \eqref{eqpropadd1}-\eqref{eqpropadd2},  we first grid the set $\text{cl}(K_{\bar{\epsilon}})$  using a sequence of nonempty compact subsets 
$\{K_i\}^{N}_{i=1}$, where $N \in \{1,2,...,\infty\}$. That is, we assume that $\bigcup^{N}_{i = 1} K_i = \text{cl}(K_{\bar{\epsilon}}). $
and, for each $i \in \{1,2,...,N\}$, there exists a finite set 
$\mathcal{N}_i \subset \{1,2,...,N \}$ 
such that  $ K_i \cap K_j = \emptyset$  for all 
$j \notin \mathcal{N}_i$,  and 
$\mbox{int}_{\text{cl}(K_{\bar{\epsilon}})} (K_i \cap K_j) = \emptyset$ for all $j \in \mathcal{N}_i,$
where $\mbox{int}_{K}(\cdot)$ is the interior of $(\cdot)$ relative to $K$. 

Next, since the set $\text{cl}(K_{\bar{\epsilon}})$ is forward contractive for $\Sigma_{\epsilon_1}$ and the solutions to $\Sigma_{\epsilon_1}$ starting from $\text{int}(K_{\bar{\epsilon}})$ never reach $\partial K_{\bar{\epsilon}}$ for positive times, we conclude,  for each $i \in \{1,2,...,N\}$ and for each $T_i > 0$,  that $R_{\Sigma_{\epsilon_1}} (T_i, K_i) \backslash  K_i \subset \text{int} (K_{\bar{\epsilon}})$.
Furthermore,   using \red{Lemmas \ref{lemprepre} and \ref{lem1-}},   we conclude,  for each $i \in \{1,2,...,N\}$, the existence of $\bar{\rho}_i>0$ such that $R_{\Sigma_{\epsilon_1}} (T_i ,  K_i + \bar{\rho}_i \mathbb{B})$ is compact for $T_i > 0$ sufficiently small.  As a result,  \red{using Lemma \ref{lem1-} in the Appendix},  we conclude, for each $i \in \{1,2,...,N\}$,  the existence 
of $\delta_i > 0$ such that 
$|y|_{\partial K_{\bar{\epsilon},\delta}}  \geq \delta_i $ for all $y \in  R_{\Sigma_{\epsilon_1}} (T_i, K_i)$ and 
$|y|_{\partial K_{\bar{\epsilon}}} \geq \delta_i$ for all 
$y \in R^b_{\Sigma_{\epsilon_1}} (T_i, K_i)$.
Next,  \red{using Lemma \ref{lem1-}},  we conclude upper semicontinuity of 
$R_{\Sigma_{\epsilon_1}}$ and $R^b_{\Sigma_{\epsilon_1}}$ on each $[0,T_i] \times (K_i + \bar{\rho}_i \mathbb{B})$.  
In particular,  for each $i \in \{1,2,...,N\}$,  there exists 
$\rho_i \in (0,  \bar{\rho}_i]$ such that 
\begin{align*}
|y|_{R_{\Sigma_{\epsilon_1}} (T_i, K_i)} 
& \leq \delta_i/2 \qquad \forall y \in  R_{\Sigma_{\epsilon_1}} 
(T_i, K_i + \rho_i \mathbb{B}), \\
|y|_{R^b_{\Sigma_{\epsilon_1}} 
(T_i, K_i)} & \leq \delta_i/2 \qquad \forall y \in R^b_{\Sigma_{\epsilon_1}} (T_i, K_i + \rho_i \mathbb{B}).
\end{align*} 
Hence, for each $i \in \{1,2,...,N\}$, 
$|y|_{\partial K_{\bar{\epsilon},\delta}}  \geq \delta_i/2$ for all $y \in  R_{\Sigma_{\epsilon_1}} (T_i, K_i + \rho_i \mathbb{B}), 
$
and
$R^b_{\Sigma_{\epsilon_1}} (T_i, K_i + \rho_i \mathbb{B}) \subset \text{int}(K_{\bar{\epsilon}})$.
As a result, we conclude that 
$$ \bigcup_{t \geq 0} \bigcup_{i \in \mathbb{N}} R_{\Sigma_{\epsilon_1}} (t, K_i + \rho_i \mathbb{B}) =   \bigcup_{i \in \mathbb{N}}  R_{\Sigma_{\epsilon_1}} (T_i , K_i + \rho_i \mathbb{B}), $$ 
which is a closed set,  and
$ \bigcup_{t \geq 0} \bigcup_{i \in \mathbb{N}}  R_{\Sigma_{\epsilon_1}} (t,  K_i + \rho_i \mathbb{B}) \subset 
\text{int}(K_{\bar{\epsilon},\delta}). $
Next, we introduce the function $\rho : \text{cl}(K_{\bar{\epsilon}})  \rightarrow 
\mathbb{R}_{>0}$ given by 
$ \rho(x) := \min \{\rho_i : 
i \in \{1,2,...,N\} ~ \text{s.t.} ~ x \in  K_i \}. $

The function $\rho$ is lower semicontinuous and allows us to conclude that the solutions  to $\Sigma_{\epsilon_1}$ 
starting from 
$\bigcup_{x \in \text{cl}(K_{\bar{\epsilon}})} (x + \rho(x) \mathbb{B})$ 
remain in $\text{int} (K_{\bar{\epsilon},\delta})$ for any positive times.  

As a last step,  we use   \cite[Theorem 1]{katvetov1951real} to conclude the existence of 
$\rho_1 : K_{\bar{\epsilon}} \rightarrow \mathbb{R}_{>0}$ continuous and satisfying  $\rho_1(x) \leq \rho(x)$  for all $x \in K_{\bar{\epsilon}}$.    
Hence, 
\begin{align}  \label{eqextend+++}
 K_{\bar{\epsilon},  \rho_1,  \epsilon_1}   \subset \bigcup_{t \geq 0} \bigcup_{x \in \text{cl}(K_{\bar{\epsilon}})}  R_{\Sigma_{\epsilon_1}} (t, x + \rho_1(x) \mathbb{B}) \subset \text{int}(K_{\bar{\epsilon},\delta}). 
\end{align}

Let  $\rho_2 : \mathbb{R}^n \backslash 
\text{int}(K_{\bar{\epsilon}, \delta})  \rightarrow \mathbb{R}_{>0}$ as $\rho_2(x) = \frac{1}{2} |x|_{K_{\bar{\epsilon},\rho_1,\epsilon_1}}$  to conclude that  $ x + \rho_2(x) \mathbb{B} \subset \mathbb{R}^n \backslash 
\cl(K_{\bar{\epsilon}, \rho_1, \epsilon_1})$  for all  
$x \in \mathbb{R}^n \backslash \text{int}(K_{\bar{\epsilon},\delta})$.    Finally,  after extending   $\rho_2$ to $\mathbb{R}^n$ and taking  $ \rho_o (x) := \min \{ \rho_1(x), \rho_2(x) \}$ for all $x \in \mathbb{R}^n$,  \eqref{eqpropadd1}-\eqref{eqpropadd2}  follow.

To prove \ref{item:I8+},  we use \eqref{eqpropadd} to
 conclude that  $\cl(K_{\bar{\epsilon}, \delta}) \backslash \text{int} (K_{\bar{\epsilon}}) \subset \text{int} (U_1)$.   
 Hence, we can always find a closed set $\hat{U}_1 \subset \text{int}(U_1)$ such
$ \cl(K_{\bar{\epsilon}, \delta}) \backslash \text{int} (K_{\bar{\epsilon}}) \subset \text{int} (\hat{U}_1), $
which implies that \eqref{eqlem102} holds. 
Next,   we let the function $\rho_3 : \hat{U}_1 \rightarrow \mathbb{R}_{>0}$ given by 
$ \rho_3 (x) := |x|_{\mathbb{R}^n \backslash U_1}$.  After extending $\rho_3$ to $\mathbb{R}^n$ and taking $\rho_o(x) := \min \{\rho_i(x) : i \in \{1,2,3 \} \}$ for all $x \in \mathbb{R}^n$,  both \eqref{eqpropadd1}-\eqref{eqpropadd2} and 
\ref{item:I8+} follow. 

To prove \ref{item:I2+a}), we use the fact $K_{\bar{\epsilon}} \subset K_{\bar{\epsilon}, \rho_o, \epsilon_1}$ and \ref{item:R1}.  Similarly, to prove \ref{item:I2+b}), we combine \eqref{eqpropadd1}-\eqref{eqpropadd2} and \ref{item:I3}.  To prove \ref{item:I2+c}) and  \ref{item:I2+d}),  we use Lemma \ref{lem2} while replacing the sets 
$(X_o, K_{\bar{\epsilon}})$ therein by the sets  $(K_{\bar{\epsilon},\rho_o}, K_{\bar{\epsilon},\rho_o,\epsilon_1})$.
\end{proof}

\newpage 

\subsection{Proof of Step 3}

 Let $\bar{\epsilon}_o$ be a robustness margin.  Using Lemma \ref{lem2bis}, we conclude the existence of a robustness margin $\bar{\epsilon}$ satisfying \eqref{eqlosttrack} and such that \eqref{eqinter} holds.   Furthermore,    let $\epsilon_1$ and $\epsilon_2$ be robustness margins such that  \eqref{eqepsilon1} and \eqref{eqepsilon3} hold.  Let $K_{\bar{\epsilon}}$ be the reachable set defined in \eqref{eqinfreachbis} such that \eqref{eqinter} holds. Let $U_1$ be the closed set used in Lemma \ref{lem3}, and $\delta$ be the continuous function used in Lemma \ref{lem6}. Let $K_{\bar{\epsilon}, \rho_o,\epsilon_1}$ be the reachable set defined in \eqref{Kinter1+} and $\rho_o$ and $\hat{U}_1$ 
be, respectively,  the function and the closed set used in Lemma \ref{lem6+}.  Let the map  
$B_{\cl(K_{\bar{\epsilon},\rho_o,\epsilon_1})} : \mathbb{R}^n \rightarrow \mathbb{R} \cup 
\{\pm \infty \}$ defined as in \eqref{eqbarcand} while replacing $(K, \Sigma)$ therein by $(\cl (K_{\bar{\epsilon},\rho_o,
\epsilon_1}),\Sigma_{\epsilon_2})$.  

We show that 
\begin{enumerate} [label={S3\arabic*)},leftmargin=*]
\item \label{item:B1--+} 
$\cl(K_{\bar{\epsilon},\rho_o,\epsilon_1}) = \{ x \in \mathbb{R}^n : B_{\cl(K_{\bar{\epsilon},\rho_o,\epsilon_1})}(x) \leq 0 \}$ and  $B_{\cl(K_{\bar{\epsilon},\rho_o,\epsilon_1})}$ is a barrier candidate with respect to $(X_o,X_u)$. 
\item \label{item:B3--+}  $B_{\cl(K_{\bar{\epsilon},\rho_o,\epsilon_1})}$ is upper semicontinuous on $\mathbb{R}^n \backslash \text{int} (K_{\bar{\epsilon},\rho_o,\epsilon_1})$.

\item \label{item:B4--+}  $B_{\cl(K_{\bar{\epsilon},\rho_o,\epsilon_1})}$ is lower semicontinuous on $\cl(K_{\bar{\epsilon},\rho_o,\epsilon_1})$.

\item \label{item:B5--+} 
Along each solution $\phi$ to 
$\Sigma_{\epsilon_2}$ with
$\phi(\dom \phi) \subset U_1$,  the map 
$t \mapsto B_{\cl(K_{\bar{\epsilon},\rho_o,\epsilon_1})}(\phi(t))$ is strictly decreasing,  and, 
for each $(t,t') \in \dom \phi \times \dom \phi$ with $t' \geq t$,  we have 
$$ B_{\cl(K_{\bar{\epsilon},\rho_o,\epsilon_1})}
(\phi(t')) - B_{\cl(K_{\bar{\epsilon},\rho_o,\epsilon_1})}(\phi(t)) \leq - (t'-t). $$
\end{enumerate}

Indeed,  we establish the proof by verifying the conditions in  Lemma \ref{lem5-} while replacing $(K, \Sigma)$ therein by $(\cl(K_{\bar{\epsilon}, \rho_o, \epsilon_1}), \Sigma_{\epsilon_2})$. To do so, we start using \eqref{eqpropadd},  \eqref{eqpropadd1}-\eqref{eqpropadd2}, and the definition of the set $K_{\bar{\epsilon}, \rho_o, \epsilon_1}$, to conclude that  $\cl(K_{\bar{\epsilon}, \rho_o, \epsilon_1}) \cap \cl(X_u) = \emptyset$,  $ \partial K_{\bar{\epsilon}, \rho_o, \epsilon_1} \subset \text{int} (U_1)$,  and  $X_o \subset \text{int}(K_{\bar{\epsilon}, \rho_o, \epsilon_1})$.  
To complete the proof, it remains to show that conditions \ref{item:71--}-\ref{item:74--}, while replacing $(K, \Sigma)$ therein by $(\cl(K_{\bar{\epsilon}, \rho_o, \epsilon_1}), \Sigma_{\epsilon_2})$, are satisfied. Indeed, using \ref{item:I2+a}, we conclude that \ref{item:71--} holds. Using \ref{item:I2+b}, we conclude that \ref{item:72--} holds. Using \ref{item:I2+c}, we conclude that \ref{item:73--} holds. Finally, using \ref{item:I2+d}, we conclude that \ref{item:74--} holds.

\subsection{Proof of Step 4}

Given a smooth function $\rho_2 : \mathbb{R}^n \rightarrow \mathbb{R}_{>0}$,  and $v \in \mathbb{B}$, we let the function
$B^v_{\cl(K_{\bar{\epsilon},
\rho_o,\epsilon_1})} : \mathbb{R}^n \rightarrow \mathbb{R}$ introduced in \eqref{eq.hatB}.  
We show that,  for each  $\rho_2 \in \mathcal{C}_+$ satisfying  
\begin{align} \label{eqrho}
 \rho_2(x) \leq \rho_1(x) := \min \{ \epsilon_2(x), \rho_o(x) \} \qquad \forall x \in \mathbb{R}^n, 
\end{align}
the function $B^{v}_{\cl(K_{\bar{\epsilon},\rho_o,\epsilon_1})}$ satisfies the following properties:

\begin{enumerate} [label={S4\arabic*)},leftmargin=*]
\item \label{item:Bhat1}  $B^{v}_{\cl(K_{\bar{\epsilon},\rho_o,\epsilon_1})}$ is a barrier function candidate with respect to 
$(X_o,X_u)$.
\item \label{item:Bhat2}  $B^{v}_{\cl(K_{\bar{\epsilon},\rho_o,\epsilon_1})}$ is upper semicontinuous on the set $K_{vo} := \{x \in \mathbb{R}^n : B^{v}_{\cl(K_{\bar{\epsilon},\rho_o,\epsilon_1})}(x) \geq 0 \}$.

\item \label{item:Bhat3} $B^{v}_{\cl(K_{\bar{\epsilon},\rho_o,\epsilon_1})}$ is lower semicontinuous on the set 
$K_{vi} := \{x \in \mathbb{R}^n : B^{v}_{\cl(K_{\bar{\epsilon},\rho_o,\epsilon_1})}(x) \leq 0 \}$.
\end{enumerate}
Furthermore,  we show  the existence of a smooth function $\rho_2 : \mathbb{R}^n \rightarrow \mathbb{R}_{>0}$ satisfying \eqref{eqrho} such that, for each $v \in \mathbb{B}$, the function $B^{v}_{\cl(K_{\bar{\epsilon},\rho_o,\epsilon_1})}$ in \eqref{eq.hatB} satisfies the following property:
\begin{enumerate} 
[label={S45)},leftmargin=*]
\item \label{item:Bhhat5} Along each solution $\phi$ to $\Sigma$ with $\phi(\dom \phi) \subset \hat{U}_1$,  the map $t \mapsto B^{v}_{\cl(K_{\bar{\epsilon},\rho_o,\epsilon_1})}(\phi(t))$ is strictly decreasing,  and,  for each 
$(t',t) \in \dom \phi \times \dom \phi$,  with $t' \geq t$,  we have 
$$ B^{v}_{\cl(K_{\bar{\epsilon},\rho_o,\epsilon_1})}(\phi(t')) - B^{v}_{\cl(K_{\bar{\epsilon},\rho_o,\epsilon_1})}(\phi(t)) \leq - (t'-t).  $$
\end{enumerate}
~~ \\ ~~ \\

Indeed,  we showed  that $B_{\cl(K_{\bar{\epsilon}, \rho_o, \epsilon_1})}$ is well defined on $\mathbb{R}^n$. Therefore, $B^v_{\cl(K_{\bar{\epsilon}, \rho_o, \epsilon_1})}$  is also well defined on $\mathbb{R}^n$.  
Next,  by definition,  we know that 
$$ B_{\cl(K_{\bar{\epsilon}, \rho_o, \epsilon_1})}(y) \leq 0 \qquad  \forall y \in \cl(K_{\bar{\epsilon}, \rho_o, \epsilon_1}). $$
Then, given $v \in \mathbb{B}$, we conclude that
$$ B^{v}_{\cl(K_{\bar{\epsilon}, \rho_o, \epsilon_1})}(x) \leq 0 \qquad  \forall (x+\rho_2(x)v) \in \cl(K_{\bar{\epsilon}, \rho_o, \epsilon_1}). $$  

To show that 
\begin{align} \label{eqtoprove}
B^{v}_{\cl(K_{\bar{\epsilon}, \rho_o, \epsilon_1})}(x) \leq 0 \qquad \forall x \in X_o,
\end{align}
we use \eqref{eqrho} and \eqref{Kinter1+},  to conclude that
$\rho_2(x) \leq \rho_o(x)$ for all $x \in \mathbb{R}^n$  and 
$ x + \rho_o(x) v \in K_{\bar{\epsilon}, \rho_o, \epsilon_1}$ for all $x \in K_{\bar{\epsilon}}$.  
Next,  since $X_o \subset K_{\bar{\epsilon}}$, we conclude that
$ x + \rho_2(x) v \in K_{\bar{\epsilon}, \rho_o, \epsilon_1}$ for all $x \in X_o$, 
which proves \eqref{eqtoprove}. 

Next, to show that 
\begin{align} \label{eqproofstep}
B^{v}_{\cl(K_{\bar{\epsilon}, \rho_o, \epsilon_1})}(x)>0 \qquad \forall x \in X_u, 
\end{align}
we start noting that 
$B_{\cl(K_{\bar{\epsilon}, \rho_o, \epsilon_1})}(y)>0$ for all $y \in \mathbb{R}^n \backslash \cl(K_{\bar{\epsilon}, \rho_o, \epsilon_1}). $
Next, using Lemma \ref{lem6} and \eqref{eqpropadd},  we conclude that $ \cl(K_{\bar{\epsilon},\delta}) \cap \cl(X_u) = \emptyset. $
Hence,  
$ x \in \mathbb{R}^n \backslash \cl(K_{\bar{\epsilon}, \delta})$ for all $x \in X_u$. 
Moreover, using \eqref{eqpropadd1}-\eqref{eqpropadd2},  we conclude that
$ x+\rho_o(x) v \notin \cl(K_{\bar{\epsilon}, \rho_o,\epsilon_1}) $ for all $x \in \mathbb{R}^n \backslash \cl(K_{\bar{\epsilon}, \delta})$.   Therefore, \eqref{eqproofstep} is satisfied and $B^{v}_{K_{\bar{\epsilon}, \rho_o,\epsilon_1}}$ is a barrier function candidate with respect to 
$(X_o,X_u)$. Now, since the function $\rho_2$ is  continuous,  the upper and lower semicontinuity properties of $B^v_{\cl(K_{\bar{\epsilon}, \rho_o,\epsilon_1})}$ are inherited from the those  of $B_{\cl(K_{\bar{\epsilon}, \rho_o,\epsilon_1})}$. 
Indeed, for each $x \in K_{vo}$, 
$y := x + \rho_2(x) v \in \mathbb{R}^n \backslash \text{int} (K_{\bar{\epsilon},\rho_o,\epsilon_1})$. 
According to item \ref{item:B3--+}, the function $B_{\cl(K_{\bar{\epsilon},\rho_o,\epsilon_1})}$ is upper semicontinuous at $y$. Hence, $B^v_{\cl(K_{\bar{\epsilon},\rho_o,\epsilon_1})}$ is upper semicontinuous on $K_{vo}$. Similarly, for each $x \in K_{vi}$,  $y := x + \rho_2(x) v \in  \cl(K_{\bar{\epsilon},\rho_o,\epsilon_1})$. 
Using item \ref{item:B4--+}, the function $B_{\cl(K_{\bar{\epsilon},\rho_o,\epsilon_1})}$ is lower semicontinuous at $y$.  Hence, $B^v_{\cl(K_{\bar{\epsilon},\rho_o,\epsilon_1})}$ is lower semicontinuous on $K_{vi}$.

Finally,   using Lemma \ref{lem8}, we conclude the existence of a smooth function $\rho_2 : \mathbb{R}^n \rightarrow \mathbb{R}_{>0}$ satisfying \eqref{eqrho}, such that \ref{item:star} holds. That is, for this choice of the function $\rho_2$, for each $v \in \mathbb{B}$, for each $\phi$ solution to $\Sigma$ starting from $x_o \in \hat{U}_1$ and remaining in $\hat{U}_1$, and given $(t,t') \in \dom \phi \times \dom \phi$ with $t' \geq t$, we conclude the existence of $\psi : [t,t'] \rightarrow \mathbb{R}^n$ solution to $\Sigma_{\rho_1}$ such that $ \psi(s) := \phi(s) + \rho_2(\phi(s)) v$ for all $s \in [t,t']$.  Thus, the function $\psi$ is solution to $\Sigma_{\epsilon_2}$.   Next, using the last item in Lemma \ref{lem6+},  we  conclude that $ \psi(s) := \phi(s) + \rho_2(\phi(s)) v \in U_1$ for all $s \in [t,t']$.   Finally, using \ref{item:B5--+}, we conclude that
$$ B_{\cl(K_{\bar{\epsilon},\rho_o,\epsilon_1})}
(\psi(t')) - B_{\cl(K_{\bar{\epsilon},\rho_o,\epsilon_1})}(\psi(t)) \leq - (t'-t), $$
which implies that \ref{item:Bhhat5} is satisfied.

\subsection{Proof of Step 5}
Finally,  we consider  the function $B : \mathbb{R}^n \rightarrow \mathbb{R}$ introduced in \eqref{eqBcertificate}  and we show  the following properties:
\begin{enumerate} [label={S5\arabic*)},leftmargin=*]
\item \label{item:B2-+}  $B$ is a barrier candidate with respect to $(X_o,X_u)$. 
\item \label{item:B3-+}  $B$ is continuously differentiable.
\item \label{item:B4-+} 
For each $x \in \hat{U}_1$,
$\langle \nabla B(x), \eta \rangle 
\leq -1$ for all $\eta \in F(x)$.

\item \label{item:B5-+} 
$\partial K \subset  \text{int}(\hat{U}_1)$, where $K := \{ x \in \mathbb{R}^n : B(x) \leq 0 \}$.
\end{enumerate}

To prove \ref{item:B2-+}, 
we use \ref{item:Bhat1} to conclude that, for each $v \in \mathbb{B}$, the function $B^{v}_{\cl(K_{\bar{\epsilon},\rho_o,\epsilon_1})}$ is a barrier function candidate with respect to $(X_o, X_u)$. Hence, according to \eqref{eqBcertificate}, it follows that $B$ is also a barrier function candidate. Indeed, for each $x \in X_o$, we have 
$$B^{v}_{\cl(K_{\bar{\epsilon},\rho_o,\epsilon_1})}(x) \leq 0 \qquad \forall v \in \mathbb{B}. $$
Hence, since $\Psi$ is nonnegative, we conclude that
$$ \int_{\mathbb{R}^n} B^{v}_{\cl(K_{\bar{\epsilon},\rho_o,\epsilon_1})}(x) \Psi(v) dv \leq 0 \qquad  \forall x \in X_o.  $$ 
Furthermore, for each $x \in X_u$, we have 
$$ B^{v}_{\cl(K_{\bar{\epsilon},\rho_o,\epsilon_1})}(x) > 0 \qquad 
\forall v \in \mathbb{B}. $$   
Hence, using the fact that $\psi$ is nonnegative and there exist points in $\mathbb{B}$ where $\psi$ not zero, we conclude that 
$$ \int_{\mathbb{R}^n} B^{v}_{\cl(K_{\bar{\epsilon},\rho_o,\epsilon_1})}(x) \Psi(v) dv > 0 \qquad  \forall x \in X_u. $$ 

To prove \ref{item:B3-+}, we use Lemma \ref{lemsmoothing} under the fact that $B_{\cl(K_{\bar{\epsilon}, \rho_o, \epsilon_1})}$ is locally bounded and,  at every $x \in \mathbb{R}^n$,   either upper or lower semicontinuous.   

To prove \ref{item:B4-+}, we consider a solution $\phi$ to $\Sigma$ such that $\phi(\dom(\phi))\subset \hat{U}_1$.   For each $(t,t') \in \dom \phi \times \dom \phi$, we have 
\begin{equation*}
\begin{aligned}  
B(\phi(t')) & = \int_{\mathbb{B}} B^{v}_{\cl(K_{\bar{\epsilon},\rho_o,\epsilon_1})}(\phi(t')) \Psi(v) dv
\\ &
\leq \int_{\mathbb{B}} [B^{v}_{\cl(K_{\bar{\epsilon},\rho_o,\epsilon_1})}(\phi(t)) - (t'-t)] \Psi(v) dv 
\\ &
= \int_{\mathbb{B}} B^{v}_{\cl(K_{\bar{\epsilon},\rho_o,\epsilon_1})}(\phi(t)) \Psi(v) dv - (t'-t) = B(\phi(t))-(t'-t),
\end{aligned}
\end{equation*}
where the first inequality is obtained  using \ref{item:Bhhat5}. Next, using Lemma \ref{lemA9bis+}, \ref{item:B4-+} follows.

To prove \ref{item:B5-+}, we use \eqref{eqrho} and \eqref{Kinter1+} to conclude that
$x+ \rho_2(x) \mathbb{B} \subset x+ \rho_o(x) \mathbb{B} \subset K_{\bar{\epsilon}, \rho_o, \epsilon_1}$
for all $x \in K_{\bar{\epsilon}}.$ 
Hence, for each $v \in \mathbb{B}$, we have 
$$ B^{v}_{\cl(K_{\bar{\epsilon}, \rho_o, \epsilon_1})}(x)  = B_{\cl(K_{\bar{\epsilon}, \rho_o, \epsilon_1})}(x + \rho_2(x) v) \leq 0 \qquad  \forall x \in K_{\bar{\epsilon}}. $$
The latter implies that 
\begin{align} \label{eqtouse} 
B(x) \leq 0 \qquad \forall x \in K_{\bar{\epsilon}}. 
\end{align}

Next,  using \eqref{eqpropadd1}-\eqref{eqpropadd2} and \eqref{eqrho},  we have
$$
x + \rho_2(x) \mathbb{B}  \subset
x + \rho_o(x) \mathbb{B} \subset \mathbb{R}^n \backslash 
\cl(K_{\bar{\epsilon}, \rho_o, \epsilon_1}) \qquad \forall x \in \mathbb{R}^n \backslash \text{int}(K_{\bar{\epsilon},\delta}).  
$$
Hence, for each 
$v \in \mathbb{B}$, 
$$ B^{v}_{\cl(K_{\bar{\epsilon}, \rho_o, \epsilon_1})}(x) = B_{\cl(K_{\bar{\epsilon}, \rho_o, \epsilon_1})}(x + \rho_2(x)v)  > 0 \qquad \forall x \in \mathbb{R}^n \backslash \text{int}(K_{\bar{\epsilon}, \delta}).  $$
The latter yields
\begin{align} \label{eqtouse1}
B(x) > 0 \qquad \forall x \in \mathbb{R}^n \backslash \text{int}(K_{\bar{\epsilon}, \delta}).
\end{align}

Using \eqref{eqtouse} and \eqref{eqtouse1}, we conclude that
$K := \{ x \in \mathbb{R}^n : B(x) \leq 0 \} \subset K_{\bar{\epsilon},\delta}$
and 
$K_{\bar{\epsilon}} \subset K$.
Hence, \ref{item:B5-+} follows using \eqref{eqpropadd}.

\section{Conclusion}
In this paper, we establish the equivalence between robust safety and the existence of a smooth barrier certificate,  in the context of continuous-time systems modeled by differential inclusions.   Our result requires only continuity of the set-valued dynamics and empty intersection between the closures of the initial and unsafe sets. We relax most of the assumptions used in existing literature such as boundedness of the safety region, smoothness of the system's dynamics,  and uniqueness of solutions. 
 Future works pertain to address the considered problem in the more general context of hybrid systems.

\appendix

\section*{Appendix}

Given $T \geq 0$,  $M \geq 0$,  and $x \in \mathbb{R}^n$,  we introduce set of solutions 
\begin{equation*}
\begin{aligned}
 \hspace{-0.4cm}
  \mathcal{A}_{\Sigma}  (T,M, x) :=  
\left\{ \phi \in \mathcal{S}^T_{\Sigma}(x) : 
 \esssup_{t \in [0,T]} |\dot{\phi}(t)| \leq M  \right\}, 
\end{aligned}
\end{equation*}
where  $\mathcal{S}^T_{\Sigma}(x)$ is the set of solutions to 
$\Sigma$ starting from $x$ whose domain is $[0,T]$. 

\begin{remark} \label{remreach}
\blue{Note that the set $ \mathcal{A}_{\Sigma} (T,M, x)$ coincides with to the set of solutions to $\Sigma'$ starting from $x$ and defined on $[0,T]$,  where 
$\Sigma'$ is defined as 
\begin{align} \label{eqSigma'} 
\Sigma' :  \dot{x}  \in F'(x) := F(x) \cap M \mathbb{B} \qquad x \in \mathbb{R}^n.   
\end{align}
In other words,  we have that $\mathcal{A}_{\Sigma} (T,M, x) = \mathcal{S}^T_{\Sigma'}(x)$.  

Under Assumption \ref{ass1},  we conclude that $F'$ is upper semi-continuous with compact and convex images on $\mathbb{R}^n$,  but can be empty on a subset $E \subset \mathbb{R}^n$,  which is necessarily open since the graph of $F$ is closed.
Hence,  the set $ \mathcal{A}_{\Sigma} (T,M, x)$ is not guaranteed to be nonempty. }
\end{remark}

\blue{ In the following lemma,  we recall three different consequences of Assumption \ref{ass1} on the set of solutions.   The first and the third properties are studied in   \cite[Theorem 1, Page 103]{aubin2012differential} under a slightly different set of assumptions.  The second property can be directly proved using  \cite[Lemma 4]{refId0}.  }

\begin{lemma} \label{lem1-}
Consider system $\Sigma$ such that  Assumption \ref{ass1} holds and let $T$,  $M \geq 0$.  Then,   
\begin{itemize}
\item The  map $x \mapsto \mathcal{A}_{\Sigma}(T,M, x)$  is outer semicontinuous and locally bounded in the metric of uniform convergence.  
\end{itemize}
As a consequence,  given  $U \subset   \mathbb{R}^n$ such that  $R_\Sigma(T, U)$ is bounded,  we can show that 
\begin{itemize}
\item  The  map $x \mapsto \mathcal{S}^{T}_{\Sigma}(x)$  is both outer and upper semicontinuous on $U$ in the metric of uniform convergence.  

\item  $R_{\Sigma}$ and $R^b_{\Sigma}$ are outer semicontinuous on $[0,T] \times U$.
\end{itemize}
\end{lemma}
'
\begin{proof}
\underline{Proof of the first item:} We start noting that 
\begin{align} \label{eqprflem}
\phi ([0,T]) \subset x + T M \mathbb{B} \qquad \forall  \phi \in \mathcal{A}_\Sigma(T,M,x).   
\end{align} 
As a consequence,  the map 
$x \mapsto \mathcal{A}_\Sigma(T,M,x)$ 
is locally bounded.  

\ifitsdraft

 Next,  to show outer semicontinuity, we pick a sequence of initial conditions 
$\{ x_i \}^{\infty}_{i=1}$ that converges to $x_o$, and a sequence of solutions $\{ \phi_i \}^{\infty}_{i=1}$,  with $\phi_i \in 
\mathcal{A}_\Sigma(T,M,x_i)$,  that converges to a function  
$\phi$ (in the sense of uniform convergence over  $[0,T]$).  
Using a consequence of the Ascoli-Arzel{\'a} and the Alaoglu Theorems,  which is established in \cite[Theorem 4, Page 13]{aubin2012differential},  we conclude that $\phi$ is absolutely continuous on $[0,T]$ and that $\{ \dot{\phi}_i \}^{\infty}_{i=1}$ converges weakly to $\dot{\phi}$; that is,  for any $y \in \mathcal{L}^1([0,T])$,  $$\lim_{i \rightarrow \infty} \int^T_{0} \dot{\phi}_i(t) y(t) dt =  \int^T_{0} \dot{\phi}(t) y(t) dt. $$
Now,  we propose to show that $\esssup \{ |\dot{\phi}(t)|  : t \in [0,T]  \}  \leq M$ using contradiction.  That is,  we let $I \subset [0,T]$ be a non-null-measure interval such that  $|\dot{\phi}(t)| > M$ for all $t \in I$.  
Furthermore,  we note that, for any $y \in \mathcal{L}^1([0,T])$,  
$$ \lim_{i \rightarrow \infty} \int^T_{0} \dot{\phi}_i(t) y(t) dt =  \int^T_{0} \dot{\phi}(t) y(t) dt \leq M \int^T_{0} |y(t)| dt.  $$
So, for 
$$y(t) :=  \left\{  
\begin{matrix} 
\phi(t)/|\phi(t)| & \text{if} ~ t \in I ~ \text{and} ~ \phi(t) \neq 0
\\
 0 & \text{otherwise},  
\end{matrix}
 \right.
$$
we obtain
$$ \int_I \dot{\phi}(t) y(t) dt \leq M \int_I |y(t)| dt = M \int_I 1 dt $$
and at the same time 
$$  \int_I \dot{\phi}(t) y(t) dt =  \int_I |\dot{\phi}(t)| dt > M \int_I  1 dt,  $$
which yields to a contradiction.  Finally,  using the Convergence Theorem in \cite[Theorem 1, Page 60]{aubin2012differential} while noting that upper semicontinuous maps with compact images are upper hemicontinuous
\cite[Corollary 2.4.1 and Definition 2.4.2]{Aubin:1991:VT:120830},   we conclude that  
$$ (\phi(t), \dot{\phi}(t)) \in \gph(F) \quad  \text{for almost all}  ~ t \in [0,T].   $$  
Hence,  $\phi$ is a solution to $\Sigma$ and it belongs to $\mathcal{A}(T,M,x)$. 

\else
\blue{In view of Remark \ref{remreach},   we know that $\mathcal{A}_{\Sigma} (T,M, x) = \mathcal{S}^T_{\Sigma'}(x)$, where $\Sigma'$ is defined in \eqref{eqSigma'}. 
As a result,  we can introduce the constrained system  
$$ \Sigma_c : \dot{x} \in F'(x) \qquad x \in C := \mathbb{R}^n \backslash E,  $$
where $E$ is the open set  on which $F'$ is empty. 

Since a solution $\phi \in \mathcal{S}^T_{\Sigma'}(x)$ starting from $x \in C$ cannot leave the closed set $C$,  we conclude that
\begin{align*}
\mathcal{A}_{\Sigma}(T,M,x) & = \mathcal{S}^T_{\Sigma'}(x) = \mathcal{S}^T_{\Sigma_c}(x) \qquad \forall x \in C,  
\\
\mathcal{A}_{\Sigma}(T,M,x) & = \emptyset \qquad  \qquad \qquad \qquad \quad \forall x \in E = \mathbb{R}^n \backslash C.  
\end{align*}
As a result,  using \cite[Lemma 5.27]{goebel2012hybrid} under local boundedness of $x \mapsto \mathcal{S}^T_{\Sigma_c}(x)$, 
 we can verify  that the map $x \mapsto \mathcal{S}^T_{\Sigma_c}(x)$ is outer semi-continuous on $C$,  which completes the proof.  }

\fi

\underline{Proof of the second item:} 
\ifitsdraft

We note that since $R_\Sigma(T, U)$ is bounded and $F$ is locally bounded,  then,  by letting
\begin{align} \label{eqMM}
M :=  \sup \{ |F(x)|  : x \in R_\Sigma(T,U)  \},  
\end{align}
we conclude that 
$$ \mathcal{A}_{\Sigma}(T,M,x) =  \mathcal{S}^T_\Sigma(x) \qquad \forall x \in U.  $$   
As a result,  using the first item,  more precisely,  by exploiting the upper semicontinuity of $x \mapsto \mathcal{A}_{\Sigma}(T,M,x)$,   we conclude 
that, for each $\epsilon > 0$,  there exists $\delta > 0$ such that, for each $x_o \in (x + \delta \mathbb{B}) \cap U$, we have
$ \mathcal{A}_{\Sigma}(T,M,x_o) \subset \mathcal{A}_{\Sigma}(T,M,x) + \epsilon \mathbb{B},  $
where $\mathbb{B}$ is given by 
$ \mathbb{B} := \{ \zeta \in \mathcal{C}([0,T]) :   |\zeta|_\infty \leq 1   \}. $

\else 
\blue{ The same way,  using \cite[Lemma 5.27]{goebel2012hybrid}  under boundedness of $\mathcal{S}^T_{\Sigma}(U)$,  we can verify  that the map 
 $x \mapsto \mathcal{S}^T_{\Sigma}(x)$ is outer semi-continuous on $U$.   This allows us,  in view of Remark \ref{remplus},  to conclude its upper semicontinuity on $U$.  } 
 
 \fi

\underline{Proof of the third item:}  
We start showing outer semicontinuity of $R^b_{\Sigma}$. 
Indeed,  let $(t,x) \in [0,T] \times U$ and let two sequences $\left\{ (t_{i},x_{i}) \right\}^{\infty}_{i=0} \subset [0,T] \times U$ and $\left\{ y_i \right\}^{\infty}_{i=0} \subset \mathbb{R}^n$ such that $\lim_{i \rightarrow \infty} (t_{i}, x_{i}) = (t,x)$, $y_i \in R^b_{\Sigma}(t_{i},x_{i})$, and $\lim_{i \rightarrow \infty} y_i = y \in \mathbb{R}^n$.   
Outer semicontinuity of $R^b_{\Sigma}$ at $(t,x)$ follows if we show that $y \in R^b_{\Sigma}(t,x)$.   To this end,  we consider a sequence of  solutions $\left\{\phi_i \right\}^{\infty}_{i=0}$  such that $\phi_i \in \mathcal{S}^{T}_{\Sigma}(x_{i})$ and $y_i = \phi_i(t_{i})$ for all $i \in \{0,1,...\}$.  Now,  since $\left\{\phi_i \right\}^{\infty}_{i=0}$ is uniformly bounded,   by passing to an adequate subsequence,  we conclude the existence of a continuous function $\phi : [0,T] \rightarrow \mathbb{R}^n$ such that $\lim_{i \rightarrow \infty} \phi_i(t) = \phi(t)$ for all $t \in [0,T]$; hence,  $\phi(0) = x$ and $y = \lim_{i \rightarrow \infty} \phi_i(t_{i}) = \phi(t)$.   Finally,  using outer semicontinuity of $x \mapsto \mathcal{S}^T_\Sigma(x)$ on $U$,   we conclude that $\phi \in  \mathcal{S}^T_\Sigma(x)$ and thus $y \in R^b_{\Sigma}(t,x)$.

Now,  to show outer semicontinuity of $R$, we consider two sequences $\left\{(t_{i},x_{i})\right\}^{\infty}_{i=0} \subset [0,T] \times U$ and 
$\left\{ y_i \right\}^{\infty}_{i=0} \subset \mathbb{R}^n$ such that 
$\lim_{i \rightarrow \infty} (t_{i}, x_{i}) = (t,x)$, 
$y_i \in R(t_{i},x_{i})$,  and  $\lim_{i \rightarrow \infty} y_i = y \in \mathbb{R}^n$. 
Outer semicontinuity of $R$ at $(t,x)$ follows if we show that $y \in R(t,x)$.  
Having $y_i \in R(t_{i},x_{i})$, for each $i \in \mathbb{N}$, implies the existence of $t'_i \in [0,t_{i}]$ such that $y_i \in R^b(t'_{i},x_{i})$, for each $i \in \mathbb{N}$. By passing to an adequate subsequence, we conclude the existence of $t' \in [0,t]$ such that $t' = \lim_{i \rightarrow \infty} t'_i$. Hence,  since $R^b$ is outer semicontinuous, we conclude that 
$y \in R^b(t',x) \subset R(t,x)$. 
\end{proof}

\blue{We next present a  useful consequence of having $F$ locally bounded.   \ifitsdraft \else A proof can be found in \cite{maghenem2022converse}.  \fi }

\begin{lemma}  \label{lemprepre}
\blue{ Consider the differential inclusion $\Sigma$ and assume that the map $F$ therein is locally bounded.  
Then,  for each compact set $K \subset \mathbb{R}^n$,   there exist $b>0$ and $\bar{T}>0$ such that the sets
$R_{\Sigma} (\bar{T},  (K + b \mathbb{B}))$
and $R_{\Sigma^-} (\bar{T},  (K + b \mathbb{B}))$
are bounded.  }
\end{lemma}

 \ifitsdraft 

\begin{proof}
Since $F$ is locally bounded,  we conclude the existence of $L$ and $M >0$ such that
$   F( K + L \mathbb{B}) \subset M \mathbb{B}.   $
Now,  by taking $b := L/2$ and $\bar{T} := L/(2M)$,  it follows that 
$ F( K + (b + \bar{T} M) \mathbb{B}) \subset M \mathbb{B}.   $
This implies that the solutions to $\Sigma$ starting from $K + b \mathbb{B}$, over the window of time $[0,\bar{T}]$ or $[-\bar{T},0]$,  cannot leave the set $ K + (b + \bar{T} M) \mathbb{B}$.  As a result, we obtain 
$$ R_{\Sigma} (\bar{T},  (K + b \mathbb{B}) ) \subset  K + (b + \bar{T} M) \mathbb{B} \quad \text{and} \quad 
 R_{\Sigma^-} (\bar{T},  (K + b \mathbb{B}) ) \subset  K + (b + \bar{T} M) \mathbb{B}.  $$
Hence,  $R_{\Sigma} (\bar{T},  (K + b \mathbb{B}))$
and
$R_{\Sigma^-} (\bar{T},  (K + b \mathbb{B}))$
are bounded.      
\end{proof}

\fi

\blue{The following Lemma follows from the combination of    \cite[Theorem 5.2.1 and Lemma 5.1.2]{Aubin:1991:VT:120830}. 
\ifitsdraft
\else
 A detailed proof can be found in \cite{maghenem2022converse}.  
 \fi }

\begin{lemma} \label{lemA14}
Consider system $\Sigma$ in \eqref{eq.1} such that Assumption \ref{ass1} holds.  A closed set $K \subset \mathbb{R}^n$ is forward invariant for $\Sigma$ if 
\begin{align} \label{eqtengApp} 
F(x) \subset C_K(y) \qquad \forall  x \in \mathbb{R}^n \backslash K, \quad \forall y \in \text{Proj}_{K}(x).
\end{align}
\end{lemma}

\ifitsdraft

\begin{proof}
\blue{We start using \cite[Theorem 5.2.1]{Aubin:1991:VT:120830} to conclude that the set
$K$ is forward invariant if
\begin{equation}  
\label{eq.econe}
\begin{aligned} 
F(x) \subset E_K(x) & \qquad \forall x \in \mathbb{R}^n \backslash K,
\end{aligned}
\end{equation}
where 
$$E_K(x) := \left\{ v\in \mathbb{R}^n :\liminf_{h \rightarrow 0^+} \frac{|x+hv|_K - |x|_K}{h} \leq 0  \right\}.  $$
Next,  to complete the proof,  we use \cite[Lemma 5.1.2]{Aubin:1991:VT:120830} to conclude that
$$  C_K(y) \subset E_K(x) \qquad \forall y \in \text{Proj}_K(x), \quad \forall x \in \mathbb{R}^n \backslash K.   $$
Hence,   \eqref{eq.econe} is verified under  \eqref{eqtengApp}.}
\end{proof} 

\fi

\textcolor{blue}{ The following lemma can be deduced from \cite[Theorem 6.3]{clarke2008nonsmooth}.  
Although formulated in the nonsmooth setting and for locally-Lipschitz dynamics,  the same proof applies to our case.  \ifitsdraft
\else A detailed proof is in  \cite{maghenem2022converse}.  \fi }

\begin{lemma} \label{lemA9bis+}
Consider system $\Sigma$ in \eqref{eq.1} such that Assumption  \ref{ass3} holds.  Consider an open set  $O \subset \mathbb{R}^n$  and a continuously-differentiable function $B : \mathbb{R}^n \rightarrow \mathbb{R}$ such that, 
along each solution $\phi$ satisfying $\phi(\dom \phi) \subset O$,  the map $t \mapsto B(\phi(t))$ is nonincreasing.  Then, 
$$ \langle \nabla B(x) , \eta  \rangle \leq 0 \qquad  \forall \eta \in F(x),  \quad \forall x \in O.  $$ 
\end{lemma}

\ifitsdraft

\begin{proof} 
Let $x_o \in O$ and $v_o \in F(x_o)$.  Since $F$ is continuous and has closed and convex images,  using Michael's selection theorem \cite{michael1956continuous},  we conclude the existence of a continuous selection $v : U(x_o) \rightarrow \mathbb{R}^n $ such that $v(x) \in F(x)$ for all $x \in U(x_o)$ with $v(x_o) = v_o$.  
Next,  using \cite[Proposition 3.4.2]{Aubin:1991:VT:120830},  we conclude the existence of a nontrivial continuously differentiable solution $\phi$ starting from $x_o$ solution to the system $\dot{x} = v(x)$; thus,  $\phi$ is also solution to $\Sigma$.  Furthermore, we consider a sequence  $\left\{ t_i \right\}^{\infty}_{i=0} \subset \dom \phi$ such that $\lim_{i \rightarrow \infty} t_i = 0$.  Note that
$
\frac{d}{dt} (B(\phi(t)))|_{t=0}  = \langle \nabla B(x_o),v(x_o) \rangle 
 = \lim_{t_i \rightarrow 0} \frac{B(\phi(t_i)) - B(\phi(0)) }{t_i}  \leq 0. 
$
\end{proof}

\fi

\blue{ The remaining lemmas are deduced from  \cite{goebel2012hybrid},  where they are  formulated for the general context of hybrid inclusions.   
\ifitsdraft
\else
The proofs can be found in \cite{maghenem2022converse}.  \fi }

\begin{lemma} [Lemma 7.37.  \cite{goebel2012hybrid}] \label{lem8}
Consider system $\Sigma$ in \eqref{eq.1} such that Assumption  \ref{ass3} holds.  Then, for each $\rho_1 \in \mathcal{C}_+$, there exists a smooth function $\rho_2 : \mathbb{R}^n \rightarrow \mathbb{R}_{>0}$ satisfying  
\begin{align} \label{eqrhorho0}
\rho_2(x) \leq \rho_1(x) \qquad \forall x \in \mathbb{R}^n
\end{align}
such that the following property holds. 
\begin{enumerate} [label={P\ref{lem8}\arabic*)},leftmargin=*]
\item \label{item:star} 
For each $\phi \in \mathcal{S}_{\Sigma}(x_o)$,  for each $v \in \mathbb{B}$,  and for each $t \in \dom \phi \cap \mathbb{R}_{\geq 0}$,  the function  $\psi: [0,t] \rightarrow  \mathbb{R}^n$ given by $\psi(s) := \phi(s) + \rho_2(\phi(s)) v$  is solution to $\Sigma_{\rho_1}$. 
\end{enumerate}
\end{lemma}
 
\ifitsdraft
\begin{proof}
Let $\phi \in \mathcal{S}_{\Sigma}(x_o)$,  $t \in \dom \phi \cap \mathbb{R}_{\geq 0}$,    and  $v \in \mathbb{B}$.  Given $\rho_2 : \mathbb{R}^n \rightarrow \mathbb{R}_{>0}$ continuous,  the function 
$\psi$ satisfies  
\begin{align*}
\dot{\psi}(s)  = \dot{\phi}(s) + \nabla \rho_2(\phi(s)) \dot{\phi}(s) v  \subset  F(\phi(s)) + |\nabla \rho_2(\phi(s))| f(\phi(s)) \mathbb{B} \quad  \text{for a. a.} ~ s \in [0,t],
\end{align*}
where $f(x) :=  \sup \{ |\zeta| : \zeta \in F(x) \}$,  which is  continuous in our case.   Now,   given 
$\rho : \mathbb{R}^n \rightarrow \mathbb{R}_{>0}$ continuous,  we show the existence of $\rho_2 : \mathbb{R}^n \rightarrow \mathbb{R}_{>0}$ continuous
such that 
\begin{align} \label{eqTCK}
 \rho_2(x) \leq \rho(x) ~ \text{and} ~  |\nabla \rho_2(x)|  \leq  \frac{\rho(x)}{  f(x) +1}  \quad \forall x \in \mathbb{R}^n.  
\end{align}
Under \eqref{eqTCK}, we conclude that 
\begin{align} \label{eqIntermuse}
\dot{\psi}(s)   \subset  F(\phi(s)) + \rho(\phi(s)) \mathbb{B} \quad \text{for a. a.} ~ s \in [0,t].  
\end{align}
Now,  to verify \eqref{eqTCK},  we partition $\mathbb{R}^n$ using a locally finite cover $\{K_i\}^\infty_{i=1}$ with $\cl(K_i)$ compact, and subordinate to this
cover a smooth partition of unity $\{ \psi_i \}^{\infty}_{i=1}$.   Finally,   we let 
$$ \rho_2(x) :=  \sum^\infty_{i=1} \frac{2^{1-i} a_i}{\max_{z \in K_i} \max \{ \psi_i(z),  |\nabla \psi_i(z)|  \} }  \psi_i(x), $$
where $a_i \in (0,1)$ such that $a_i \leq \rho(x)$ for all $x \in K_i$ and
$ a_i \sup_{z \in K_i }  f(z)  \leq 1$.

To complete the proof,  we use Lemma \ref{lemA333} in the Appendix to conclude that given
  $\rho_1 \in \mathcal{C}_+$,  we can find $\rho \in \mathcal{C}_+$  such that 
\begin{equation}
\label{eqinterPr}
\begin{aligned}
\rho(x) & \leq \rho_1(x) \qquad \forall x \in \mathbb{R}^n,
\\
F(x) + \rho(x) \mathbb{B} &  \subset F(y) + \rho_1(y) \mathbb{B}  \qquad \forall y \in x + \rho(x) \mathbb{B}. 
\end{aligned}
\end{equation}
As a result,  since 
$$ \psi(s) \in \phi(s) + \rho_2(\phi(s)) \mathbb{B} \subset \phi(s) + \rho(\phi(s)) \mathbb{B} \quad  \forall s \in [0,t], 
$$    
 applying \eqref{eqinterPr},  we conclude that 
$$ F(\phi(s)) + \rho(\phi(s)) \mathbb{B}   \subset F(\psi(s)) + \rho_1(\psi(s)) \mathbb{B} \quad \forall s \in [0,t],  $$
and using \eqref{eqIntermuse},  we obtain 
$$ \dot{\psi}(s)  \subset  F(\psi(s)) + \rho_1(\psi(s)) \mathbb{B} \quad \text{for almost all} ~ s \in [0,t]. $$
The latter implies that $\psi$ is solution to $\Sigma_{\rho_1}$. 
\end{proof}
  
  \fi

\begin{lemma} [Section 7.6.0.3.  in \cite{goebel2012hybrid}] \label{lemsmoothing}
Let $B_o : \mathbb{R}^n \rightarrow \mathbb{R}$ be  locally bounded.  Assume that,  at every $x \in \mathbb{R}^n$,  $B_o$ is either upper or lower semicontinuous.   Then,  the function 
$B : \mathbb{R}^n \rightarrow \mathbb{R}$ 
given by
\begin{align} \label{eqBinteg} 
B(x) := \int_{\mathbb{R}^n} B_{o}(x + \rho_o(x) v) \Psi(v) dv,  
\end{align}
where $\rho_o : \mathbb{R}^n \rightarrow \mathbb{R}_{>0}$ and $\Psi : \mathbb{R}^n \rightarrow [0,1]$ are smooth functions such that
\begin{equation}
\label{eqPsi} 
\begin{aligned} 
\Psi(v) = 0  \quad \forall v \in \mathbb{R}^n \backslash \mathbb{B} \quad \text{and} \quad
\int_{\mathbb{R}^n}  \Psi(v) dv = 1, 
\end{aligned}
\end{equation}
 is continuously differentiable. 
\end{lemma}

\ifitsdraft

\begin{proof}
 Consider the change of coordinate 
$ w := x + \rho_o(x) v. $
Furthermore, we introduce the closed set  
$ \mathbb{B}_{\rho_o}(x) := x + \rho_o(x) \mathbb{B}. $ 
Hence, we obtain 
$$B(x)= \int_{\mathbb{R}^n} B_{o}(w) \Psi\left(\frac{w-x}{\rho_o(x)} \right) \frac{dw}{\rho_o(x)}.$$
Next, we let 
$$ g(x):= \int_{\mathbb{R}^n} f(x,w) dw  :=  \int_{\mathbb{R}^n} B_{o}(w) \Psi\left(\frac{w-x}{\rho_o(x)} \right) dw. $$

We will show that $g$ is continuously differentiable,  which would imply that $B$ is continuously differentiable,  since $\rho_o$ is smooth. 
\begin{itemize}
\item[a.] Note that 
$ w \mapsto f(x,w)$ 
is $\mathcal{L}^1$ since 
$ w \mapsto \Psi\left(\frac{w-x}{\rho_o(x)} \right) $
is null outside the bounded set  $\mathbb{B}_{\rho_o}(x)$. 

\item[b.] Note that the map $ x \mapsto \Psi \left( \frac{w-x}{\rho_o(x)} \right) $ is smooth and null outside $\mathbb{B}_{\rho_o}(x)$.  Hence,  $ w \mapsto  \nabla_x f(x,w) \in \mathcal{L}^1$. 
 
 \item[c.] Since   the map $ x \mapsto  \nabla_x \left( \Psi\left(\frac{w-x}{\rho_o(x)} \right) \right)$ is smooth and null outside the bounded set $\mathbb{B}_{\rho_o}(x)$, we conclude the existence of a positive constant $k$ such that
 $$ \sup_{y \in \mathbb{B}_{\rho_o}(x)} \displaystyle \left\lvert \nabla_x \left( \Psi\left(\frac{w-y}{\rho_o(y)} \right) \right) \displaystyle \right\rvert \leq k. $$
 Hence, 
 $ |\nabla_x f(x,w)| \leq  b(w), $
  where $b : \mathbb{R}^n \rightarrow \mathbb{R}_{\geq 0}$ is the upper semicontinuous function given by
 $$ b(w) := 
 \left\{ 
 \begin{matrix}
 k |B_{o}(w)| & \text{if} ~ w \in \mathbb{B}_{\rho_o}(x) 
 \\ 
 0 & \text{otherwise}.
 \end{matrix}
 \right. $$
Being null outside the  set  $\mathbb{B}_{\rho_o}(x)$,    $w \mapsto b(w)$ is $\mathcal{L}_1$. 
\end{itemize} 

As a result,   using  \cite[Lemma 7.38]{goebel2012hybrid},  we conclude that $g$ is differentiable and 
$$ \nabla_x g(x) = \int_{\mathbb{R}^n} B_{o}(w) \nabla_x \left(\Psi\left(\frac{w-x}{\rho_o(x)} \right) \right)^\top dw,  $$
which is continuous.  
\end{proof}

\begin{lemma}  \label{lemA333}
Let $F: \mathbb{R}^n \rightrightarrows \mathbb{R}^n$ such that Assumption \ref{ass3} holds.   Then,  for each   $\rho_1 \in \mathcal{C}_+$, there exists  $\rho \in \mathcal{C}_+$ such that
\begin{align} \label{eqcontic}
F(x) + \rho(x) \mathbb{B}  \subset F(y) +  \rho_1 (y)\mathbb{B} \quad 
\forall y \in x + \rho(x) \mathbb{B}. 
\end{align}
\end{lemma}

\begin{proof}
We grid $\mathbb{R}^n$ using a sequence of nonempty compact subsets $\{I_i\}_{i=1}^{N} \subset \mathbb{R}^n$,  where $N \in \{1,2,3, \dots ,\infty \}$,  such that,  for each $i \in \{1,2,...,N\}$,  there exists 
$\mathcal{N}_i \subset \{1,2,...,N \}$  finite such that 
$ I_i \cap I_j = \emptyset$  for all  $j \notin \mathcal{N}_i$,   and $ \mbox{int}(I_i \cap I_j) = \emptyset$ for all $j \in \mathcal{N}_i\backslash \{i\}$.
Since $F$ is continuous; thus,  uniformly continuous  on each $I_i + \mathbb{B}$,  
$i \in \{1,2,...,N\}$ \footnote{  
\blue{The map $F: K \rightrightarrows \mathbb{R}^n$ is \textit{uniformly continuous} if, for each $\epsilon > 0$, there exists 
$\delta > 0$ such that, for each $x \in K$, we have  
$$ 
|F(x_1) - F(x_2) |_H \leq \epsilon \qquad  \forall x_1, x_2 \in x + \delta \mathbb{B}.  
 $$  
When $K$ is compact,  using the same arguments as in the  single-valued case,  continuity and uniform continuity become equivalent. }},  
  we conclude that,  for each $\varepsilon_i > 0$,  there exists 
$\delta_i \in (0,1]$ such that 
$ F(x) + \delta_i  \mathbb{B}  \subset F(y)+ \varepsilon_i \mathbb{B}$  for all $x, y \in I_i + \mathbb{B}$  such that $ |y - x| \leq  \delta_i$.    In particular,  since $\delta_i \leq 1$,  we conclude that
$ F(x) + \delta_i  \mathbb{B}  \subset F(y)+ \varepsilon_i \mathbb{B}$ for all $x \in I_i$,   for all $y \in x +  \delta_i \mathbb{B}$.     
Now,  if we let $\varepsilon_i := \min_{y \in I_i} \rho_1(y)$,  we conclude the existence of $\delta_i > 0$ such that 
$ F(x) + \delta_i \mathbb{B}  \subset F(y)+ \rho_1(y) \mathbb{B}$ for all $x \in I_i$,  for all $y \in x +  \delta_i \mathbb{B}$.   Next, we introduce the function $\delta :\mathbb{R}^n \to \mathbb{R}_{>0}$ given by
$\delta(x) :=\text{min}\{\delta_i: i\in \{1,2,\dots , N\} \;\ \text{such that} \;\ x\in I_i\} $.  By definition,  $\delta$ is lower semi-continuous; hence,  $-\delta$ is upper semicontinuous.  
Using   \cite[Theorem 1]{katvetov1951real},  we conclude the existence of $\rho \in \mathcal{C}_+$ and satisfying 
$- \rho(x) \geq -\delta(x)$ for all $x \in \mathbb{R}^n$.
Hence, \eqref{eqcontic} follows. 
\end{proof}

\fi

\balance

\bibliographystyle{unsrt}      
\bibliography{biblio.bib}

\end{document}